\documentclass[reqno,11pt]{amsart}
\usepackage{cases}
\usepackage{mathrsfs}
\usepackage[T1]{fontenc}
\usepackage{mathrsfs}
\usepackage{amsmath,latexsym,amssymb,amsfonts,amsbsy, amsthm}
\usepackage{bm}
\usepackage[usenames]{color}
\usepackage{xspace,colortbl}
\usepackage{epsfig}
\usepackage{graphicx}
\usepackage{subfigure}
\usepackage{amsmath,amsfonts,amsthm,amssymb,amscd}
\input amssym.def 
\input amssym.tex
\usepackage{color}
\usepackage[title,titletoc,toc]{appendix}
\usepackage{bigints}
\allowdisplaybreaks[4]

\newtheorem {remark1}{Remark}[section]
\newtheorem{theorem1}{Theorem}[section]

\newtheorem{definition1}{Definition}[section]

\numberwithin{equation}{section}
	\title{On a  rod-Kadomtsev-Petviashvili  shallow water equation in two dimensions}

\author[Zhao]{Tiantian Zhao}
\address{Tiantian Zhao\newline
	School of Mathematics and Statistics, Huazhong University of Science and Technology, Wuhan, Hubei 430074, China}
\email{zhaotiantian@hust.edu.cn}

\author[Yan]{Kai Yan\textsuperscript{*}}
\address{Kai Yan (Corresponding auther)\newline
	School of Mathematics and Statistics, Huazhong University of Science and Technology, Wuhan, Hubei 430074, China}
\email{kaiyan@hust.edu.cn}
\begin{document}
	\renewcommand{\thefootnote}{*}\footnotetext{Corresponding author}\renewcommand{\thefootnote}{\arabic{footnote}}
	\begin{abstract}
		%The Camassa-Holm-Kadomstev-Petviashvili (CH-KP) equation ia a two dimensional Camass-Holm equation-type with weakly transverse effect for the horizontal velocity component and captrues stronger nonlinear effects that the classical dispersive integrable equation li
		In this paper, we derive  a new  two-dimensional  rod-Kadomtsev-Petviashvili (rod-KP) equation from the incompressible and irrotational  three-dimensional Euler equation under the shallow water scaling. We  establish the local well-posedness of the Cauchy problem in a suitable  Sobolev space and derive  the blow-up criterion for strong solutions via the  energy method. Then, by exploring two appropriate conservation laws of the rod-KP equation, we are able to construct its global strong solution when the physical dimensionless parameter $\sigma=0$, and on the other hand  produce the finite-time blow up solutions under some certain conditions when $\sigma \neq 0$. Furthermore, we  present a uniqueness continuation property for the solutions. Finally, we investigate the existence of traveling-wave solutions in order to highlight  the influence of weak transverse effects on wave stability, and we  also exhibit   the  symmetry of solitary waves in the propagation direction. 
	\end{abstract}

		\maketitle
%	\author[Liu]{Yue Liu}	\address{Yue Liu \newline	Department of Mathematics, University of Texas at Arlington, Arlington, TX 76019, USA}	\email{yliu@uta.edu}
	%\author{Tiantian Zhao} \address{Tiantian Zhao \newline	School of Mathematics and Statistics\\ Huazhong University of Science and Technology, Wuhan 430074,  China}\email{m202270029@hust.edu.cn}
	
	%\author{Kai Yan} \address{Kai Yan (Corresponding author)\newline School of Mathematics and Statistics\\ Huazhong University of Science and Technology, Wuhan 430074,  China} \email{kaiyan@hust.edu.cn}

% -----------------------------------------


\vskip 0.2cm

%%%%%%%%%%%%%%%%%%%%%%%%%%%%%%%%%%%%%%%%%%%%%%%%%%%%%%%%%%%%%%
%%%%%%%%%%%%%%%%%%%%%%%%%%%%%%%%%%%%%%%%%%%%%%%%%%%%%%%%%%%%%
\renewcommand{\theequation}{\thesection.\arabic{equation}}
\setcounter{equation}{0}
	
	\section{Introduction}
	Recently, the Camassa-Holm-Kadomtsev-Petviashvili (CH-KP) equation has been derived  as a model in the context of incompressible and irrotational shallow water wave theory, taking the form 
	\begin{align}
		(u_t -u_{txx} + \kappa u_x + 3uu_x - \left( 2u_xu_{xx} +uu_{xxx} \right) )_x + u_{yy} =0,  \label{CH-KP eqaution}
	\end{align}
where $u=u(t,x,y) $ denotes  the horizontal velocity field. 
	The CH-KP equation (\ref{CH-KP eqaution}) is  obtained under the Camassa-Holm (CH) scaling $\mu \ll 1$ and $\varepsilon = O (\sqrt{\mu})$ \cite{Gui-Liu-Luo-Yin AIM 2021}. Here,  $\varepsilon =\frac{a}{h_0}$ and $\mu =\frac{h_0^2}{\lambda^2}$ are two small dimensionless parameters measuring nonlinearity and shallowness, respectively, where $a$ is the typical amplitude, $\lambda$ is the wavelength of the surface wave, and $h_0$ is the typical depth of the water.   CH-KP equation (\ref{CH-KP eqaution})  generalizes the CH equation in the same way as the Kadomtsev-Petviashvili  (KP) equation generalizes the classical Korteweg-de Vries (KdV) eqaution.  
	By capturing both strong unidirectional dynamics
	and weak transverse modulation within a single framework, the CH-KP model offers a powerful tool for analyzing long-crested, moderately nonlinear shallow-water waves, with direct
	relevance to geophysical phenomena such as tsunami evolution, tidal bores, and nearshore
	wave-current interactions.
	
	 The formal Hamiltonian structure of the CH-KP equation (\ref{CH-KP eqaution}) was established in \cite{Gui-Liu-Luo-Yin AIM 2021}. Gui {\it et. al} also studied 
	the   local well-posedness in suitable Sobolev space, the formation of singularities and existence of traveling wave solutions.  
	Later on, the transverse spectral stability of small periodic traveling waves of the CH-KP equation (\ref{CH-KP eqaution}) was proved in \cite{Geyer-Liu-Pelinovsky JMPA 2024}. In  a related context, the transverse spectral stability of small periodic traveling waves of the $b$-KP equation as a two-dimensional generalization of the $b$-family of CH equation was demonstrated in \cite{Chen-Fan-Wang-Xu MA 2024}.
	
	The CH-KP equation (\ref{CH-KP eqaution}) admits a reduction to the following KP-\uppercase\expandafter{\romannumeral2} equation \cite{Kadomtsev-Petviashvili SPD 1970}, 
	\begin{align}
		\left(  u_t + uu_x + u_{xxx} \right)_x + u_{yy} = 0 , \label{KP-II equation}
	\end{align}
	under the slow-modulation approximation of small-amplitude perturbations about a constant  background, together with an appropriate asymptotic truncation \cite{Geyer-Liu-Pelinovsky JMPA 2024,Mizumachi MAMS 2015}. The KP-\uppercase\expandafter{\romannumeral2}  equation derived as a model for propagation of the weakly transverse water water waves in a long wave regime, and generalizes the classical KdV equation with transversal effects. On the analytical side,   the most important point is  that  this equation is a completely integrable model \cite{Ablowitz-Clarkson}. The well-posedness of the KP-\uppercase\expandafter{\romannumeral2} equation has been studied extensively. By combining the Fourier restriction norm method with  the conservation laws, the global well-posedness has been extablished \cite{Bourgain GFA 1993}. The Cauchy problem of the KP-\uppercase\expandafter{\romannumeral2} equation is  globally well-posed in $H^s(\mathbb{R} \times \mathbb{T})$ for $s>0$ \cite{Molinet-Saut-Tzvetkov AIHCANL 2011}. In \cite{Geyer-Liu-Pelinovsky JMPA 2024}, the small data globally well-posedness  and scattering reslut are presented in the scale-invariant Sobolev space $\dot{H}^{-\frac{1}{2},0}(\mathbb{R}^2)$. 
	
 By suppressing the $y$-dependence in the CH-KP equation (\ref{CH-KP eqaution}), one obtains the  celebrated  CH equation  \cite{Camassa-Holm PRL 1993}:
	\begin{align}
		u_t -u_{txx} + \kappa u_x + 3uu_x- \left( 2u_xu_{xx} + uu_{xxx}\right) =0,
		\label{CH equation}
	\end{align}
a  unidirectional shallow-water model notable for its orbitally stable peaked solitons (peakons)\cite{Camassa-Holm PRL 1993,Constantin-Strauss CPAM 2000} and finite-time wave breaking \cite{Constantin-Escher AM 1998}.
 The CH equation also models the propagation of axially symmetric waves in hyperelastic rods \cite{Dai Act M 1998} and possesses  a bi-Hamiltonian structure \cite{Fuchssteiner-Fokas PD 1981}, as well as being  completely integrable \cite{Camassa-Holm PRL 1993}. It also admits a geometric interpretation as a  geodesic flow on the diffeomorphism group of the circle \cite{Constantin-Kolev CMH 2003}. The  Cauchy problem for the CH equation (\ref{CH equation}) has been extensively analyzed, with results on well-posedness, blow-up, and global existence \cite{Constantin AIF 2000,Constantin-Escher AM 1998,Constantin-Escher CPAM 1998,Danchin JDE 2003}.  On the other hand, the global weak solutions  to the CH  equation (\ref{CH equation}) is investigated in  \cite{Xin-Zhang CPAM 2000}. 
 %The advantage of the CH equation in comparison with the KdV equation lies in the fact that the CH equation not only has orbitally stable peakon solutions \cite{Camassa-Holm PRL 1993,Constantin-Strauss CPAM 2000}, but also models wave breaking (namely, the wave remains bounded while its slope becomes unbounded in finite time \cite{Whitham 1974}.

The CH equation also admits many  generalizations,  one of them  is the following rod equation: 
\begin{align}
	u_t -u_{txx} +\kappa u_x + 3uu_x -\sigma \left( 2u_xu_{xx} + uu_{xxx} \right) =0,\label{Rod equation}
\end{align}
which was derived by Dai in 1998 \cite{Dai Act M 1998}.
This model  describes  the propagation of nonlinear waves in  cylindrical hyper-elastic rods, where  $u(t,x)$ represent the radial strech relatiove to a presented state, and the physical parameter $\sigma$ is a real dimensionless constant.  
When $\sigma =1$, the rod equation (\ref{Rod equation}) recovers the standard  CH equation(\ref{CH equation}). For $\sigma=0$, it  reduces  to the  celebrated Benjamin-Bona-Mahony (BBM) equation, which models  surface waves in a canal. It is worth pointing out that the solutions of the BBM equation existly globally, implying that wave breaking does not occur in this case.  For $\sigma\neq 0$, the finite-time blow-up solutions for the rod equation (\ref{Rod equation}) has been established in \cite{Brandloese CMP 2014,Brandolese-Cortez JFA 2014}.  Moreover, smooth solitary waves of the rod equation (\ref{Rod equation}) have been constructed at least for some $\sigma$; see \cite{Dai Act M 1998,Lenells DCDS 2006}.

It is worth noting that the CH-KP equation (\ref{CH-KP eqaution}) can be interpreted as a two-dimensional CH-type model that incorporates weak transverse effects into the horizontal velocity component.  This naturally raises the following

{\bf Question:} {\it Does the rod  equation (\ref{Rod equation})   adimit an analogous higher-dimensional extension, similar in spirit to the CH-KP equation, that can more effectively capture the complexity of fluid dynamical phenomena?   }

One of the goals of the present paper is to formally derive a  two-dimensional rod 	equation from the Euler  equation. This new model equation named as the rod-Kadomtsev-Petviashvili (rod-KP, in short) equation,  namely 
	\begin{equation}\label{Rod-KP equation form-1}
		\left\{\begin{array}{ll}
			\left( m_t + \kappa u_x+ \sigma \left( 2mu_x +um_x \right)   + 3(1-\sigma) uu_x \right)_x + u_{yy}  =0  ,& \quad t>0, (x,y) \in  \mathbb{R}^2,\\
			m= u-u_{xx} ,& \quad  t>0, (x,y)\in \mathbb{R}^2 ,
		\end{array}\right.
	\end{equation}
 where $u=u(t,x,y)$ denotes the horizontal velocity (primarily in the $x$-direction but with weak $y$-dependence).  
 It is noted that this two-dimensional equation (\ref{Rod-KP equation form-1}) has two conserved quantities
 \begin{equation}
 	E (u):= \frac{1}{2} \int_{\mathbb{R}^2} \left( u^2+ u_x^2\right)  dxdy , \label{Rod-KP conserved quantity}
 \end{equation}
 and 
 \begin{align}
 	F_{\sigma}(u) :=\frac{1}{2} \int_{\mathbb{R}^2}  \left( \kappa u^2 + u^3 + \sigma uu_x^2 +\left(  \partial_x^{-1} \partial_y u\right)^2  \right) dxdy. \label{Rod-KP conserved quantity-1}
 \end{align} 
 
 Similar to the KP and CH-KP equations, the rod-KP equation involoves the nonlocal term $\partial_x^{-1} u_{yy}$, which aries from the asymptotic description of transverse effects in weakly  dispersive shallow-water wave flows. Physically, this operator reflects the balance between longitudinal dispersion and lateral modulation of the wave field. Analytically, however, its nonlocal nature complicates the estimation of remainder terms and necessitates refined techniques to establish the validity of the model.

Among the many asymptotic models for shallow-water dynamics, the Green-Naghdi (GN) equations \cite{Green-Naghdi JFM 1976}-also known as the Serre \cite{Serre HB 1953} or Su-Gardner equations \cite{Su-Gardner JMP 1969}-are the best known
and most widely used, particularly in coastal oceanography for the numerical simulation of surface wave propagation. Unlike the classical shallow-water equations, the GN equations capture fully nonlinear effects and allow for waves of non-negligible amplitude, while providing higher-order dispersive corrections. The derivation of  (\ref{Rod-KP equation form-1}) from a two-dimensional
GN-type system (hereafter referred to as the GN-KP-type system), itself obtained from the three-dimensional, incompressible, irrotational Euler equations under the standard shallow-water scaling.  The structure of (\ref{Rod-KP equation form-1}) exhibits three important limiting cases.  Suppressing the $y$-dependence  reduces the model to the rod equation (\ref{Rod equation}); setting $\sigma=1$ recovers the CH-KP equation (\ref{CH-KP eqaution}); and 
  when $\sigma =0$, the rod-KP equation (\ref{Rod-KP equation form-1}) reduces to the following Benjamin-Bona-Mahony-Kadomtsev-Petviashvili (BBM-KP, in short) equation, namely
\begin{align}
	\left(  u_t - u_{txx} + 3uu_x + \kappa u_x  \right)_x + u_{yy} =0. \label{BBM-KP equation}
\end{align}

%The Cauchy problem for  (\ref{Rod-KP equation form-1}) is also an interesting issue  investigated in this paper and a local well-posedness result (Theorem \ref{LWP theorem}) is obtained in a suitable Sobolev space by the energy method. 
For convenience, let us  rewrite (\ref{Rod-KP equation form-1}) with the initial data $u_0$ as the following nonlocal form:   
\begin{equation}
	\left\{
	\begin{aligned}
		&u_t + \sigma uu_x + p_x \ast \left( \frac{3-\sigma}{2}u^2 +\frac{\sigma}{2} u_x^2+ \kappa u  \right)  + p \ast \partial_x^{-1} u_{yy} =0,  &&  t>0, (x,y) \in \mathbb{R}^2, \\
		& u|_{t=0} = u_0(x,y) , &&(x,y) \in \mathbb{R}^2,
	\end{aligned}
	\right
	.
	\label{the nonlocal weak form}
\end{equation}
where $p(x) :=\frac{1}{2} e^{-|x|}$ and $\partial_x^{-1} u(x,y) := \mathscr{F}^{-1} \left( (i\xi)^{-1} \mathscr{F} (u) (\xi, \eta) \right)$, operators $\mathscr{F}$ and $\mathscr{F}^{-1}$ denote the Fourier transform and inverse of the Fourier transform in terms of variables $x$, $y$, rspectively. By the structure of  (\ref{the nonlocal weak form}), we introduce the following the Hilbert space $X^s$ for $s>0$: 
\begin{align}
	X^s (\mathbb{R}^2) := \left\lbrace  u \in H^s(\mathbb{R}^2) |  u_x \in H^s(\mathbb{R}^2) , \partial_x^{-1} u \in H^s(\mathbb{R}^2) \right\rbrace \nonumber 
\end{align}
equipped with the norm $\Vert u \Vert_{X^s}:= \left(   \Vert u \Vert_{H^s}^2 +\Vert \partial_x^{-1} u \Vert_{H^s}^2+ \Vert u_x \Vert_{H^s}^2\right)^{\frac{1}{2}}$, $\forall u \in X^s$,   and  the inner product
\begin{align}
	\left(  u, v\right)_{X^s} :=  (u,v)_{H^s} +(\partial_x^{-1} u, \partial_x^{-1} v) +  (u_x, v_x)_{H^s}, \quad \forall u, v\in X^s(\mathbb{R}^2). \nonumber 
\end{align}

Now we are in  the position to state the local well-posedness result.
\begin{theorem1}\label{LWP theorem}
	Let $\sigma \in \mathbb{R}$ and  $u_0 \in X^s(\mathbb{R}^2)$ with $s\geq 2$. Then there exists a maximal existence time  $T^{\ast} = T^{\ast} (\Vert u_0 \Vert_{X^s}) >0$ and a unique strong solution $u$ to the Cauchy problem (\ref{the nonlocal weak form})  such that 
	\begin{align}
		u\in C([0,T^{\ast} ); X^s(\mathbb{R}^2)) \cap C^1 ([0,T^{\ast}); X^{s-2} (\mathbb{R}^2)). \nonumber 
	\end{align}
	Moreover, the solution $u$ depends continuously on the initial value $u_0$, i.e., the mapping $u_0 \mapsto u : \, X^s(\mathbb{R}^2) \to C([0,T^{\ast}); X^s(\mathbb{R}^2)) \cap C^1 ([0,T^{\ast}); X^{s-2} (\mathbb{R}^2))$ is continuous. In addition, the conservation laws $E(u)$ and $F_{\sigma}(u)$ in (\ref{Rod-KP conserved quantity})-(\ref{Rod-KP conserved quantity-1}) are independent of the existence time $t>0$. 
\end{theorem1}
 By the argument in the proof of Theorem \ref{LWP theorem} in Section \ref{Local Well-posedness}, we may derive the following blow-up criterion for rod-KP equation (\ref{the nonlocal weak form}) in the case $\sigma \neq 0$. 
\begin{theorem1}[Blow-up criterion for $\sigma \neq 0$]\label{Blow-up criterion}
	Let $\sigma \neq 0$,  $u_0 \in X^s(\mathbb{R}^2)$ with $s\geq 2$ and $u$ be the corresponding solution to (\ref{the nonlocal weak form}). Assume $T^{\ast}$ is the maximal time of existence. Then 
	\begin{align}
		T^{\ast} < \infty \quad \Rightarrow \quad  \int_0^{T^{\ast}} \Vert \sigma \nabla u (\tau) \Vert_{L^{\infty}(\mathbb{R}^2)} d \tau = \infty. \label{blow-up criterion eq-1}
	\end{align}
\end{theorem1}
\begin{remark1}
For the proof of Theorem \ref{Blow-up criterion} in Section \ref{Local Well-posedness}, we provide a detailed analysis according to the regularity index $s$. More precisely,  via the Sobolev embedding, we treat the cases $s=2$, $s=3$ and $s>3$ step by step, and the remaining case $2<s<3$ are then handled by the interpolation and density arguments. It is noted to point out that in contrast to work in \cite{Gui-Liu-Luo-Yin AIM 2021} (for $\sigma =1$) where only the case $s=2$ is proved, we here give a complete proof. 

\end{remark1}

	%, which completes the proof of Remark  \ref{blow-up criterion remark}. 

While for the case $\sigma =0$, we establish the following blow-up criterion which is different from the case $\sigma  \neq 0$.
\begin{theorem1}[Blow-up criterion for $\sigma =0$]\label{BBM blow-up criterion}
	Let  $\sigma =0$,  $u_0 \in X^s(\mathbb{R}^2)$ with $s \geq 2$ and $u$ be the corresponding solution to (\ref{the nonlocal weak form}). Assume $T^{\ast}$ is the maximal time of existence. Then 
	\begin{align}
		T^{\ast} < \infty \quad \Rightarrow \quad \int_0^{T^{\ast}} \Vert u(\tau ) \Vert_{L^{\infty}(\mathbb{R}^2)} d\tau =\infty. \nonumber 
	\end{align}
\end{theorem1}
\begin{remark1}\label{blow-up criterion remark sigma0}
	Under the same assumption inTheorem \ref{BBM blow-up criterion}, the solution $u$ blows up in finite time $T^{\ast} < +\infty$ if and only if 
	\begin{align}
		\limsup\limits_{t \to T^{\ast}} \Vert u (\tau ) \Vert_{L^{\infty}(\mathbb{R}^2)} =+\infty. \label{blow-up criterion remark sigma0 eq}
	\end{align}
%Indeed,  if    $T^{\ast} <+\infty$ and (\ref{blow-up criterion remark sigma0 eq}) is not valid, then there exists $ M>0$ such that $\Vert u(t) \Vert_{L^{\infty} (\mathbb{R}^2)} \leq M$, $ \forall t\in [0,T^{\ast})$. Thus, we have  $ \int_0^{T^{\ast}} \Vert u(\tau ) \Vert_{L^{\infty}} d \tau \leq MT^{\ast}<+\infty $. Therefore, Theorem \ref{BBM blow-up criterion} implies that the maximal existence $T^{\ast} =\infty$, which contradicts the assumption that $T^{\ast} <+\infty$. Conversely, the Sobolev embedding theorem $H^s  (\mathbb{R}^2) \hookrightarrow L^{\infty}(\mathbb{R}^2) $ with $s>1$ implies that  if (\ref{blow-up criterion remark sigma0 eq}) holds, the corresponding solution $u$ blows up in finite time. 
\end{remark1}
\begin{remark1}
	The local well-posedness result (Theorem \ref{LWP theorem}) and blow-up criterion (Theorems \ref{Blow-up criterion}-\ref{BBM blow-up criterion}) also work for the Sobolev space $X^s(\Omega_1 \times \Omega_2)$, $s\geq 2$, where $\Omega_1$ and $\Omega_2$ are either the whole line $\mathbb{R}$ or the unit circle $\mathbb{T} = \mathbb{R} / \mathbb{Z}$. 
\end{remark1}
%Motivated by  the global well-posedness results for the KP equation (\ref{KP-II equation}) and BBM equation, we rely on Remark \ref{blow-up criterion remark sigma0} as a key ingredient.  By combining Sobolev estimates and the structural properties of the BBM-KP equation (\ref{BBM-KP equation}), we prove that the Cauchy problem of BBM-KP admits  a unique global solution, for all $u_0 \in X^s(\mathbb{R}^2)$, $s\geq 2$. In this argument, the coservation law $F_{0} (u)$ in (\ref{Rod-KP conserved quantity-1}) plays an indispensable role.

It is known that one of the more challenging aspects of the Cauchy problem for nonlinear evolution equation is the possible global existence or finite-time blow-up of solutions. With Theorems \ref{Blow-up criterion}-\ref{BBM blow-up criterion} and Remark \ref{blow-up criterion remark sigma0} in hand, we are able to construct the global and blow-up data for the rod-KP equation (\ref{the nonlocal weak form}).  First of all, by combining the two conservation laws $E(u)$ and $F_{\sigma} (u)$ with $\sigma =0$ in (\ref{Rod-KP conserved quantity})-(\ref{Rod-KP conserved quantity-1}), and a Sobolev inequality for $\Vert u\Vert_{L^3(\mathbb{R}^2)}$ (see Lemma \ref{Sobolev inequality lemma}- (\ref{Sobolev inequality lemma-4})): 
\begin{align}
	\Vert u \Vert_{L^3}^3 \lesssim \Vert u \Vert_{L^2}^{\frac{3}{2}} \Vert u_x \Vert_{L^2} \Vert \partial_x^{-1} u_y \Vert_{L^2}^{\frac{1}{2}}, \nonumber 
\end{align}
one gets a bound for $\Vert u_y \Vert_{L^2(\mathbb{R}^2)}$ (see (\ref{BBM global result proof eq-1})). This along with the special structure of rod-KP equation (\ref{the nonlocal weak form}) with $\sigma =0$ (i.e., the BBM-KP equation (\ref{BBM-KP equation})), and another Sobolev inequality for $\Vert u \Vert_{L^{\infty}(\mathbb{R}^2)}$ (see Lemma \ref{Sobolev inequality lemma}-(\ref{Sobolev inequality lemma-1})):
\begin{align}
	\Vert u \Vert_{L^{\infty}}^4 \lesssim \Vert u \Vert_{L^2} \Vert u_x \Vert_{L^2} \Vert u_y \Vert_{L^2} \Vert u_{xy} \Vert_{L^2}, \nonumber 
\end{align}
give rise to a bound for $\Vert  u \Vert_{L^{\infty}(\mathbb{R}^2)}$ within the lifespan. So, (\ref{blow-up criterion remark sigma0 eq}) implies the global existence of the solution to (\ref{BBM-KP equation}), and hence, the wave breaking does not occur in the case of $\sigma =0$. More precisely, we have the following result. 
\begin{theorem1}[Global data for $\sigma=0$]\label{BBM global result}
	Assume $\sigma =0$, and $u_0 \in X^s(\mathbb{R}^2)$, $s\geq 2$. Then the Cauchy problem (\ref{the nonlocal weak form}) has a  unique global strong solution $u\in C([0,+\infty); X^s(\mathbb{R}^2)) \cap C^1([0,+\infty); X^{s-2}(\mathbb{R}^2))$. 
\end{theorem1}

On the other hand, for the case $\sigma \neq 0$, the Cauchy problem (\ref{the nonlocal weak form}) admits finite-time blow-up strong solutions under some certain conditions.  Indeed, we track the dynamics of the effective blow-up quantity $B(t):= \int_{\mathbb{R}} \sigma u_xdy$ at $x=0$. After a careful estimates, one gets a Riccati-type differential inquality for 
\begin{align}
	\frac{d B}{dt} \leq -A_1 B^2(t) + A_2 \quad \text{with $A_1>0$ and $A_2 \geq 0$,} \nonumber 
\end{align}
which ensures finite-time blow-up. More precisely, we have the following result. 
\begin{theorem1}[Blow-up data for $\sigma \neq 0$]\label{blow-up data-1}
	Fix $\phi \in H^2(\mathbb{R})$ with $\phi \geq 0$ and $\int_{\mathbb{R}} \phi dy=1$. Suppose that $\sigma \neq 0$, $u_0 \in X^s(\mathbb{R}^2)$ with $s\geq 2$. If  the corresponding solution $u$ satisfies that $u(t,x,y) =-u(t,-x,y)$, for all $t \geq 0$ and $ (x,y) \in \mathbb{R}^2$, and the initial data $u_0$ satisfies
	\begin{align}
		\int_{\mathbb{R}} \sigma u_{0,x}(0,y) \phi(y) dy < -\sqrt{2C_3}, \nonumber 
	\end{align}
	where 
	\begin{equation}
		C_3:= \left\{ 
		\begin{aligned}
			&  \frac{\sigma (\sigma-3)}{2} E(u_0) \Vert \phi\Vert_{L^{\infty}},  && \quad\text{$\sigma>4$ or $\sigma<0$}, \\
			&0, & &\quad 0<\sigma \leq 4. 
		\end{aligned}
		\right.\nonumber 
	\end{equation}
	Then the corresponding strong solution $u$ blows up in finite time $T^{\ast}$ such that 
	\begin{equation}
		T^{\ast} \leq T_1 := \left\{ 
		\begin{aligned}
			& \frac{1}{\sqrt{2C_3}}\ln \left( \frac{\int_{\mathbb{R}} \sigma u_{0,x}(0,y) \phi (y) dy -\sqrt{2C_3} }{\int_{\mathbb{R}} \sigma u_{0,x}(0,y) \phi (y) dy +\sqrt{2C_3}} \right) , && \quad \text{$\sigma>4$ or $\sigma<0$}, \\
			& -\frac{2}{\int_{\mathbb{R}} \sigma u_{0,x} (0,y)\phi(y) dy}, && \quad 0<\sigma \leq 4. 
		\end{aligned}
		\right.\nonumber 
	\end{equation}
\end{theorem1}

Next, inspired by the argument in study of CH equation (\ref{CH equation}) in \cite{Linares-Ponce PAMS 2020}, we obtain the following unique continuation property for rod-KP equation (\ref{the nonlocal weak form}).
\begin{theorem1}\label{UCP Rod}
	Assume $\kappa=0$ and  $0\leq \sigma \leq 3$. Let $u$ be a  solution of (\ref{the nonlocal weak form}) in $C([0,T^{\ast}); X^s(\mathbb{R}^2)) \cap C^1([0,T^{\ast}); X^{s-2} (\mathbb{R}^2))$, $s\geq 2$ with the maximal existence time $T^{\ast}$.  If there is an open set $\Omega \subset [0,T^{\ast}) \times \mathbb{R}$ such that $u(t,x,y)=0$, $\forall (t,x) \in \Omega$, $y \in \mathbb{R}$. Then  $u=0$  in $[0,T^{\ast}) \times \mathbb{R}^2$. 
\end{theorem1}
\begin{remark1}
	We point out that  Theorem \ref{UCP Rod}  recovers  the corresponding result in \cite{Gui-Liu-Luo-Yin AIM 2021} for the case $\sigma=1$. On the other hand, if we suppress the $y$-dependence in the rod-KP equation (\ref{the nonlocal weak form}), Theorem \ref{UCP Rod} generalizes the corresponding result in \cite{Linares-Ponce PAMS 2020}. 
\end{remark1}
% \begin{remark1} 	It is obvious that the travelling wave solution of the form $u(x-ct,y)$ to the CH-KP equation (\ref{CH-KP eqaution}) is also symmetric in the variable $x$.  \end{remark1}
  
The properties of localized travelling waves play an   important role in the study of nonlinear dispersive equations.  
If we consider the travelling wave solutions of the form $u(x-ct,y)$, where $c\in \mathbb{R}$ denotes the  wave speed, then such solutions satisfy the following equation: 
\begin{align}
	\left( -c+\kappa \right) u_x - cu_{xxx} + 3uu_x -\sigma \left( 2u_xu_{xx} + uu_{xxx} \right)  + \partial_x^{-1} u_{yy} =0. \label{travelling wave solution eq-1}
\end{align} 
By studying the structure of the equation (\ref{travelling wave solution eq-1}), we obtain the following rigidity theorem which  reveals that the travelling wave solution is trivial when the wave speed $c> \min\left\lbrace 0, \kappa \right\rbrace$.
\begin{theorem1}\label{Travelling wave solution rigidity theorem}
	Let $\sigma \in \mathbb{R}$ and $u$ be the solution to the equation (\ref{travelling wave solution eq-1}).  If $c> \min\left\lbrace 0, \kappa \right\rbrace $, then $u=0$. 
\end{theorem1}
\begin{remark1}
	When $\sigma=1$,  Theorem \ref{Travelling wave solution rigidity theorem} also holds true  for the CH-KP eqaution (\ref{CH-KP eqaution}). In addition, by suppressing the $y$-dependence in the rod-KP equation (\ref{Rod-KP equation form-1}), we get  a similar property of the travelling wave solution  to the rod equation (\ref{Rod equation}). 
\end{remark1}
Then, we focus our attention on  the symmetry of the solution to equation  (\ref{travelling wave solution eq-1}). Motivated by the argument in study of the rotation-modified KP equation in \cite{Esfahani-Levandosky PAMS 2019}, the following result is established. 
\begin{theorem1}\label{Symmetry result}
	If  $0 \leq \sigma \leq 3$, $ 0 < c < \kappa$ and $u$ is a solution to the equation (\ref{travelling wave solution eq-1}), then $u$ must be  symmetric in $x$-variable. 
\end{theorem1}
Existence and nonexistence of the peaked solitary waves to the rod-KP model is an interesting issue investigated in this paper. More precisely,  we have the following result.
\begin{theorem1}\label{Peaked solitary wave}
	The rod-KP equation (\ref{Rod-KP equation form-1}) has a global weak solution in the peak form of $u(t,x,y) = ce^{-|x+\beta y -ct|}$ for some constant $\beta\in \mathbb{R}$ if and only if there holds that $\sigma=1$ and $\beta^2 +\kappa  =0$. 
\end{theorem1}
\begin{remark1}
	Similar to the $1$-D case for the rod equation (\ref{Rod equation}) in \cite{Lenells DCDS 2006},  $\sigma=1$ is the only case that   the equation admits a peaked solitary wave with a choice of different parameter $\kappa$ and wave speed $c$. 
\end{remark1}  
The remainder of the paper is organized as follows. In Section \ref{Derivation}, we formally derive  the rod-KP model (\ref{Rod-KP equation form-1})  from the incompressible and irrotational Euler equation.  Section \ref{Local Well-posedness} is devoted to establishing the local well-posedness  (Theorem \ref{LWP theorem}) and deriving the blow-up criterion (Theorems \ref{Blow-up criterion}-\ref{BBM blow-up criterion}) for the Cauchy problem (\ref{the nonlocal weak form}). In Section \ref{Global BBM}, we focus on the global existence of solutions to the BBM-KP equation (\ref{BBM-KP equation}) (Theorem \ref{BBM global result}). In Section \ref{Blow-up rod}, we  construct the blow-up solutions to rod-KP equation (\ref{the nonlocal weak form}) with $\sigma \neq 0$ (Theorem \ref{blow-up data-1}).  Next, we present the uniqueness continuation property (Theorem \ref{UCP Rod}) of the solutions to the rod-KP equation (\ref{the nonlocal weak form})  in Section \ref{Liouville-type property of the solution}.  Finally, in Section \ref{Travelling wave  solutions}, by proving Theorems \ref{Travelling wave solution rigidity theorem}-\ref{Peaked solitary wave},  we describe  the travelling-wave solutions to the rod-KP equation (\ref{the nonlocal weak form})  and discuss the symmetry of its  solitary wave. 
\section{Derivation of the rod-KP model in shallow water}
\label{Derivation}
In this section, we focus our attention on the model of the rod-KP equation (\ref{Rod-KP equation form-1}) derived from the GN-KP type system, which is considered as the model from the incompressible and irrotational three-dimensional Euler equations under certain assumptions.    
Firstly, we present the formal derivation of the GN-KP-type system. 
 Attention is given here is the so-called long-wave limit.

 In this setting, it is assumed that water flows are incompressible and inviscid with a constant density $\varrho$ and no surface tension,  and the interface between the air and the water is a free surface. Then such a motion of water flow occupying a domain $\Omega_t$ in $\mathbb{R}^3$ under the influence of the gravity $ g $ can be described by the  Euler equations, {\it viz.}
\begin{equation*} \label{R-Euler}
	\begin{cases}
		&\vec{u}_t+\left( \vec{u}\cdot\nabla \right) \vec{u} =-{1\over\varrho}\nabla P +\vec{g},\quad x\in \Omega_t,\\
		&\nabla\cdot \vec{u}=0, \quad x\in \Omega_t,\\
		&\vec{u}|_{t=0}=\vec{u_0}, \quad x\in \Omega_0,
	\end{cases}
\end{equation*}
where $ \vec u = (u, v, w )^T $ is  the fluid velocity, $ P(t,x,y,z)$ is the pressure in the fluid, $ \vec{g} = ( 0, 0, -g )^T$ with $g \approx 9.8 m/s^2$ the constant gravitational acceleration at the Earth's surface. We consider here waves at the surface of water with a flat bed, and assume that $\Omega_t=\{(x, y, z): 0<z<h_0+\zeta(t, x, y)\}$, where $h_0$ is the typical depth of the water and $\zeta(t, x, y)$ measures the deviation from the average level. The motion of inviscid irrotational fluid in the region $0 < z < h_0 + \zeta(t,x,y)$ with a constant density $\varrho$ is described by the Euler equations in the form
\begin{equation*} \label{f-plane}
	\begin{cases}
		u_t + uu_x+ vu_y + wu_z = -\frac{1}{\varrho}P_x, \\
		v_t + uv_x+ vv_y + wv_z = -\frac{1}{\varrho}P_y, \\
		w_t + uw_x + vw_y + ww_z  = -\frac{1}{\varrho}P_z - g,
	\end{cases}
\end{equation*}
the incompressibility of the fluid,
\begin{equation*}\label{incom-1}
	u_x + v_y + w_z = 0,
\end{equation*}
and the irrotational condition,
\begin{equation*}\label{irrot-1}
	(w_y-v_z, u_z-w_x, v_x-u_y)^T = (0,0,0)^T.
\end{equation*}
The pressure is written as
\begin{equation*}
	P(t, x,z) = P_a + \varrho g(h_0 - z) + p(t, x, y, z),
\end{equation*}
where $P_a$ is the constant atmosphere pressure, and $p$ is a pressure variable measuring the hydrostatic pressure distribution.

The dynamic condition posed on the surface $z = h_0 + \zeta$ yields $P = P_a$. Then there appears that
\begin{equation*}\label{perssure-1}
	p = \varrho g \zeta.
\end{equation*}
Meanwhile, the kinematic condition on the surface is given by
\begin{equation*}\label{KC-1}
	w = \zeta_t + u\zeta_x + v\zeta_y, \quad \mbox{when} \quad z = h_0 + \zeta(t, x, y).
\end{equation*}
Finally, we pose "no-flow" condition at the flat bottom $z = 0$, that is,
\begin{equation*}\label{bottom-1}
	w|_{z=0} = 0.
\end{equation*}

According to the magnitude of the physical quantities, we introduce dimensionless quantities as follows
\begin{equation*}
	x \rightarrow \lambda x,\quad y \rightarrow \lambda y,\quad z \rightarrow h_0 z,\quad \zeta \rightarrow a \zeta,\quad t \rightarrow \frac{\lambda}{\sqrt{gh_0}}t,
\end{equation*}
which implies
\begin{equation*}
	u \rightarrow \sqrt{gh_0}u,\quad v \rightarrow \sqrt{gh_0}v,\quad  w \rightarrow \sqrt{\mu gh_0} w,\quad p \rightarrow \varrho g h_0 p.
\end{equation*}

Furthermore, considering whenever $\varepsilon \rightarrow 0$,
\begin{equation*}
	u \rightarrow 0,\quad  v \rightarrow 0, \quad w \rightarrow 0,\quad p \rightarrow 0,
\end{equation*}
that is, $u, w$ and $p$ are proportional to the wave amplitude so that
we require a scaling
\begin{equation*}\label{rescall-2}
	u \rightarrow \varepsilon u,\quad v \rightarrow \varepsilon v,\quad w \rightarrow \varepsilon w,\quad p \rightarrow \varepsilon p.
\end{equation*}
Therefore the governing equations become
\begin{align}
	& u_t + \varepsilon(uu_x+v u_y+wu_z) = - p_x &\text{in}&\quad 0 < z < 1+\varepsilon \zeta(t, x, y), \label{1} \\
	& v_t + \varepsilon(uv_x+vv_y+wv_z) = - p_y &\text{in}&\quad 0 < z < 1+\varepsilon \zeta(t, x, y),  \label{2}\\
	& \mu \{w_t + \varepsilon (uw_x +vw_y+ ww_z)\}  = -p_z &\text{in} & \quad 0 < z < 1+\varepsilon \zeta(t, x, y), \label{3}\\
	& u_x+v_y + w_z = 0 &\text{in} & \quad 0 < z < 1+\varepsilon \zeta(t, x, y),\label{4}\\
	& \mu w_y-v_z  = 0, \, u_z - \mu w_x = 0, \, v_x - u_y = 0 &\text{in} & \quad 0 < z < 1+\varepsilon \zeta(t, x, y),\label{5} \\
	& p = \zeta &\text{on}& \quad z = 1 + \varepsilon \zeta(t, x, y),\label{6}\\
	& w = \zeta_t + \varepsilon (u \zeta_x+v\zeta_y) &\text{on} & \quad z = 1 + \varepsilon \zeta(t, x, y),\label{7}\\
	& w = 0 & \text{on}& \quad z = 0.\label{8}
\end{align}
Similar to the classical Green-Naghdi equations, let $\bar{u}$ be the $x$-direction average velocity and $\bar{v}$ be the $y$-direction average velocity, i.e.
\begin{align}
	&\bar{u}(t,x,y) = \frac{1}{h}\int_0^h u(t,x,y,z)~dz,\label{u-x-velocity}\\
	&\bar{v}(t,x,y) = \frac{1}{h}\int_0^h v(t,x,y,z)~dz,\label{u-y-velocity}
\end{align}
where $h =h(t,x,y)= 1 + \varepsilon \zeta(t,x,y)$. Times $h$ both sides of \eqref{u-x-velocity} and \eqref{u-y-velocity}, then differentiate the first one by $x$ and the second one by $y$, it follows
\begin{align}
	(h \bar{u})_x = \int_0^h u_x~dz +\varepsilon \zeta_x u_s, \nonumber%\label{h-u-x-velocity}, 
	\\
	(h \bar{v})_y = \int_0^h v_y~dz +\varepsilon \zeta_y v_s, %\label{h-v-x-velocity}
	\nonumber
\end{align}
where $u_s = u(t,x,y,z)|_{z = h}$ and $v_s = uv(t,x,y,z)|_{z = h}$.
Then the above equation combining with \eqref{4}, \eqref{7} and \eqref{8} gives rise to the first equation of the Green-Naghdi-KP equations
\begin{equation}
	\zeta_t + (h \bar{u})_x + (h \bar{v})_y= 0.\label{GN1}
\end{equation}

%%%%%%%%%%%%%%%%%%%%%%%%%%%%%%%%%%%%%%%%%%%%%%%%%%%%%%%%%%%%%%%%%%%%%%%%%%%%%%%%%%%%%%%%%%%%%%%%%%%%%%%%%%%%%%%%%%%%%%%%%%%%%
%%%%%%%%%%%%%%%%%%%%%%%%%%%%%%%%%%%%%%%%%%%%%%%%%%%%%%%%%%%%%%%%%%%%%%%%%%%%%%%%%%%%%%%%%%%%%%%%%%%%%%%%%%%%%%%%%%%%%%%%%%%%%

To derive the second equation of the Green-Naghdi-KP equations, let
\begin{align}
	u(t,x,y,z)=u_0(t,x,y,z)+\mu u_1(t,x,y,z)+O(\mu^2), \nonumber\\ %\label{u-u0-u1-mu},\\
	v(t,x,y,z)=v_0(t,x,y,z)+\mu v_1(t,x,y,z)+O(\mu^2),\nonumber % \label{v-v0-v1-mu},
\end{align}
here we require the shallowness parameter $\mu\ll 1$, but without any assumption on $\varepsilon$.
Considering $\mu \rightarrow 0$, the first and second equations in \eqref{5} imply $u_0$ and $v_0$ are functions independent of $z$, respectively. We can say $u_0=u_0(t,x,y)$ and $v_0=v_0(t,x,y)$. Then it follows from \eqref{4} and \eqref{5} that
\begin{align}
	\mu u_{xx} +\mu u_{yy} +u_{zz}=0, \nonumber \\%\label{mu-u-xx}, \\
	\mu v_{xx} +\mu v_{yy} +v_{zz}=0. \nonumber%\label{mu-v-xx}.
\end{align}
From the coefficients of $O(\mu^0)$ and $O(\mu^1)$, we have
\begin{align}
	 u_{0,zz}=0,~~v_{0,zz}=0, \nonumber %\label{mu0-u-v-xx}
	 \end{align} 
 and 
 \begin{align}
	u_{0,xx} +u_{0,yy} +u_{1,zz}=0, ~~v_{0,xx} +v_{0,yy} +v_{1,zz}=0.\nonumber %\label{mu1-u-v-xx}
\end{align}
Therefore, it follows from the equation of $\mu^1$ order that
\begin{align}
	u=u_0-\mu\frac{z^2}{2}(u_{0,xx}+u_{0,yy})+\mu z\Psi+O(\mu^2)  \label{mu1-u},\\
	v=v_0-\mu\frac{z^2}{2}(v_{0,xx}+v_{0,yy})+\mu z \Phi+O(\mu^2) \label{mu1-v},
\end{align}
where $\Psi=\Psi(t,x,y)$ and $\Phi=\Phi(t,x,y)$ are arbitrary functions.

Then, by using the identity $v_x=u_y$ in \eqref{5}, it infers that
\begin{equation}%\label{uy-vx}
	\nonumber
	\begin{split}
		v_{0,x}-\mu\frac{z^2}{2}(v_{0,xxx}+v_{0,xyy})+\mu z \Phi_x
		=u_{0,y}-\mu\frac{z^2}{2}(u_{0,xxy}+u_{0,yyy})+\mu z \Psi_y+O(\mu^2).
	\end{split}
\end{equation}
By comparing terms of order $\mu$, one can see $v_{0,x}=u_{0,y}$ and $\Phi_x=\Psi_y$.

Applying \eqref{4} and \eqref{8} to the identity $w=w|_{z=0}+\int_0^z w_{z'}~dz'$, there then appears that
\begin{equation}%\label{w-u-v}
	\nonumber
	\begin{split}
		w=-z(u_{0,x}+v_{0,y})+\mu\frac{z^3}{6}(u_{0,xxx}+u_{0,xyy}+v_{0,xxy}+v_{0,yyy})
		-\mu\frac{z^2}{2}(\Phi_x+\Psi_y)+O(\mu^2),
	\end{split}
\end{equation}
which along with \eqref{mu1-u}, \eqref{mu1-v} and \eqref{5} implies
\begin{equation*}
	\Phi(t,x,y)=\Psi(t,x,y)=0.
\end{equation*}
Therefore, we obtain
\begin{align}
	u=u_0-\mu\frac{z^2}{2}(u_{0,xx}+u_{0,yy})+O(\mu^2)  \label{mu1-u-1},\\
	v=v_0-\mu\frac{z^2}{2}(v_{0,xx}+v_{0,yy})+O(\mu^2) \label{mu1-v-1},
\end{align}
and
\begin{equation}\label{w-u-v-1}
	\begin{split}
		w=-z(u_{0,x}+v_{0,y})+\mu\frac{z^3}{6}(u_{0,xxx}+u_{0,xyy}+v_{0,xxy}+v_{0,yyy})
		+O(\mu^2).
	\end{split}
\end{equation}

Substituting \eqref{mu1-u-1}, \eqref{mu1-v-1} and \eqref{w-u-v-1} into \eqref{3}, it then follows from \eqref{6} and the identity $p=p|_{z=h}-\int_z^h p_{z'}~dz'$ that
\begin{equation}%\label{p-u-v}
	\nonumber
	\begin{split}
		p=\zeta-\mu\frac{h^2-z^2}{2}(u_{0,xt}&+v_{0,yt}+\varepsilon u_0(u_{0,xx}+v_{0,xy})\\
		&+\varepsilon v_0(u_{0,xy}+v_{0,yy})-\varepsilon (u_{0,x}+v_{0,y})^2)
		+O(\mu^2).
	\end{split}
\end{equation}
Then taking account of the first equation \eqref{1}, there appears that
\begin{equation}%\label{u-t-v-eta}
	\nonumber
	\begin{split}
		u_{0,t}+\varepsilon (u_0u_{0,x}&+v_0v_{0,x})+\zeta_x=\mu\frac{h^2}{2}\Big{(}u_{0,xt}+v_{0,yt}+\varepsilon u_0(u_{0,xx}+v_{0,xy})\\
		&+ \varepsilon v_0(u_{0,xy}+v_{0,yy})- \varepsilon (u_{0,x}+v_{0,y})^2\Big{)}_x
		+\mu hh_x\Big{(}u_{0,xt}+v_{0,yt}\\&+\varepsilon u_0(u_{0,xx}+v_{0,xy})
		+ \varepsilon v_0(u_{0,xy}+v_{0,yy})- \varepsilon (u_{0,x}+v_{0,y})^2\Big{)}
		+O(\mu^2).
	\end{split}
\end{equation}
It is noted that since the first equation \eqref{1} and the second equation \eqref{2} are equivalent, the final substitution of either will result in the expression.  
Recall the definition of $\bar{u}$ in \eqref{u-x-velocity} and $\bar{v}$ in \eqref{u-y-velocity},
it then follows from the expressions \eqref {mu1-u-1} and \eqref{mu1-v-1} that
\begin{align}
	\bar{u}=u_0-\mu\frac{h^2}{6}(u_{0,xx}+u_{0,yy})+O(\mu^2) ,\nonumber\\% \label{mu1-u-bar},\\
	\bar{v}=v_0-\mu\frac{h^2}{6}(v_{0,xx}+v_{0,yy})+O(\mu^2). \nonumber % \label{mu1-v-bar}.
\end{align}
Hence, a direct computation shows that
\begin{equation}%\label{u-bar-t-eta-x}
	\nonumber
	\begin{split}
		\bar{u}_t+\zeta_x&+\varepsilon (\bar{u} \bar{u}_x +\bar{u}_y \partial_x^{-1}\bar{u}_y)
		=\frac{\mu}{3(1+\varepsilon \zeta)} \Big{[} (1+\varepsilon \zeta)^3\Big{(}\bar{u}_{xt}+\partial_x^{-1}\bar{u}_{yyt} \\
		&+ \varepsilon \bar{u}(\bar{u}_{xx}+\bar{u}_{yy})+  \varepsilon \partial_x^{-1}\bar{u}_y(\bar{u}_{xy}+\partial_x^{-1}\bar{u}_{yyy}) -\varepsilon (\bar{u}_x+\partial_x^{-1}\bar{u}_{yy})^2\Big{)}  \Big{]}_x\\
		&-\varepsilon\mu\frac{hh_x}{3}\partial_x^{-1}\bar{u}_y(\bar{u}_{xy}+\partial_x^{-1}\bar{u}_{yyy})
		+\varepsilon\mu\frac{hh_y}{3}\partial_x^{-1}\bar{u}_y(\bar{u}_{xx}+\bar{u}_{yy})
		+O(\mu^2),
	\end{split}
\end{equation}
where use has been made of the fact that $\bar{v}_x=\bar{u}_y+O(\mu^2)$, i.e. $\bar{v}=\partial_x^{-1}\bar{u}_y+O(\mu^2)$.

It is noticed that
\begin{equation}%\label{u-bar-tyy}
	\nonumber
	\bar{u}_t+\zeta_x+\varepsilon (\bar{u} \bar{u}_x +\bar{u}_y \partial_x^{-1}\bar{u}_y)=O(\mu),
\end{equation}
which infers that
\begin{equation}%\label{u-bar-tyy-partial}
	\nonumber
	\partial_x^{-1}\bar{u}_{yyt}=-\zeta_{yy}-\varepsilon \partial_x^{-1}(\bar{u} \bar{u}_x)_{yy} -\varepsilon \partial_x^{-1}(\bar{u}_y \partial_x^{-1}\bar{u}_y)_{yy}+O(\mu).
\end{equation}
Then combining \eqref{GN1} and using the variable transformation $y \rightarrow  \varepsilon^{\frac{1}{2}} y$, we have the following GN-KP type equations with the error terms,
\begin{equation}%\label{GN-eq-1-y}
	\nonumber
	\left\{\begin{aligned}
		&\zeta_t+((1+\varepsilon \zeta) u)_x+\varepsilon \partial_x^{-1} u_{yy}
		+\varepsilon^2\left( \zeta \partial_x^{-1} u_y\right)_y=O(\varepsilon^2\mu^2), \\
		&u_t +\zeta_x +\varepsilon u u_x + \varepsilon^2 u_y \partial_x^{-1} u_y=\frac{\mu}{3(1+\varepsilon \zeta)}\left[( 1 + \varepsilon \zeta ) ^ { 3 } (u_{xt}+\varepsilon u u_{xx} -\varepsilon u_{x}^2- \varepsilon \zeta_{yy})\right]_x    \\
		&\quad+\frac{\varepsilon^2\mu}{3(1+\varepsilon \zeta)}\left[( 1 + \varepsilon \zeta ) ^ { 3 } \left(  u u_{yy}-\frac{1}{2}   \left(u^2\right)_{yy}-  \varepsilon \partial_x^{-1}\left(u_y\partial_x^{-1} u_y \right)_{yy}+ u_{x y}\partial_x^{-1} u_y \right.\right. \\
		&\quad \left.\left.+ \varepsilon \partial_x^{-1} u_y \partial_x^{-1} u_{yyy}-2 u_x \partial_x^{-1} u_{yy}- \varepsilon \left( \partial_x^{-1} u_{yy}\right)^2\right)\right]_x \\
		&\quad-\frac{1}{6}  \varepsilon^2 \mu\left((1+\varepsilon \zeta)^2\right)_x u_{x y}\partial_x^{-1} u_y -\frac{1}{6}  \varepsilon^3 \mu\left((1+\varepsilon \zeta)^2\right)_x \partial_x^{-1} u_y\partial_x^{-1} u_{yyy}\\
		&\quad+\frac{1}{6} \varepsilon^2 \mu\left((1+\varepsilon \zeta)^2\right)_y u_{x x}\partial_x^{-1} u_y  +\frac{1}{6}  \varepsilon^3 \mu\left((1+\varepsilon \zeta)^2\right)_y u_{yy} \partial_x^{-1} u_y
		+O(\mu^2).
	\end{aligned}\right.
\end{equation}
For simplicity, $u$ is used instead of $\bar{u}$ here and in the following article. Thus, we have the following system 
\begin{equation}
	\left\{ 
	\begin{aligned}
		&  \zeta_t + ((1+\varepsilon \zeta) u)_x + \varepsilon \partial_x^{-1} u_{yy} + \varepsilon^2 (\zeta \partial_x^{-1} u_y)_y = O(\varepsilon^2\mu^2), \\
		& u_t -\frac{\mu}{3} u_{txx} + \zeta_x + \varepsilon uu_x + \varepsilon^2 u_y \partial_x^{-1} u_y  -\varepsilon \mu \zeta_x u_{tx} -\frac{\mu}{3}u_{txx} -\frac{\mu\varepsilon}{3} \left(  uu_{xx} - u_x^2 -\zeta_{yy} \right)_x 
	\nonumber \\
	&\quad 	- \frac{2\varepsilon \mu}{3} \zeta u_{txx}  = O(\varepsilon^2  \mu, \mu^2) .
	\end{aligned}
	\right.
\end{equation}
Then, we introduce the transformation 
\begin{align}
	\tau = \varepsilon t, \quad \xi = x-Vt, \quad u(\tau, \xi,y ) = u(t,x,y ),\nonumber 
\end{align}
where $V$ will be determined later. This appears that
\begin{align}
	\frac{\partial}{\partial t} = \varepsilon \frac{\partial}{\partial \tau} -  V \frac{\partial}{\partial \xi } , \quad \frac{\partial}{\partial x} =\frac{\partial}{\partial \xi} , \nonumber 
\end{align}
which implies
\begin{equation}
	\left\{ 
	\begin{aligned}
		&  \varepsilon \zeta_{\tau}  -V \zeta_{\xi} + u_{\xi} + \varepsilon (\zeta u)_{\xi}  + \varepsilon \partial_{\xi}^{-1} u_{yy} + \varepsilon^2 (\zeta \partial_{\xi}^{-1} u_y)_y = O(\varepsilon^2\mu^2) ,\\
		& \varepsilon u_{\tau} -V u_{\xi} -\frac{\varepsilon \mu}{3} u_{\tau \xi\xi} + \frac{\mu V}{3} u_{\xi\xi\xi} + \zeta_{\xi} + \varepsilon uu_{\xi} +\varepsilon^2 u_y \partial_{\xi}^{-1} u_y +\varepsilon \mu \zeta_{\xi} u_{\xi\xi}  \\
		& \quad -\frac{\varepsilon\mu}{3} \left( uu_{\xi\xi\xi} -u_{\xi} u_{\xi\xi} -\zeta_{\xi yy} \right) +\frac{2 \varepsilon \mu }{3} \zeta u_{\xi\xi\xi} = O(\varepsilon^2 \mu, \mu^2) . 
	\end{aligned}
	\right. \label{derivation eq-2}
\end{equation}
The solution $(\zeta, u)$ is expected to be expanded as the following form
\begin{align}
	& \zeta (\tau , \xi) = \zeta_0 (\tau , \xi, y)  + \varepsilon \zeta_1 (\tau, \xi, y) +O(\varepsilon^2) ,\nonumber \\
	& u (\tau, \xi) = u_0(\tau, \xi, y) + \varepsilon u_1 (\tau, \xi, y) +O(\varepsilon^2) . \nonumber 
\end{align}
Then equation (\ref{derivation eq-2})  is converted to 
\begin{equation}
	\left\{ 
	\begin{aligned}
		&\varepsilon \zeta_{0,\tau}+ \varepsilon^2 \zeta_{1, \tau}  -V\zeta_{0,\xi} - \varepsilon V \zeta_{1,\xi} + u_{0,\xi}  +\varepsilon u_{1, \xi}   +\varepsilon \left(  \zeta_0 u_0\right)_{\xi} +\varepsilon  \partial_{\xi}^{-1} u_{0,yy} =  O(\varepsilon^2\mu,\mu^2), \\
		&  \varepsilon u_{0,\tau} + \varepsilon^2 u_{1, \tau} -Vu_{0,\xi}  - \varepsilon V u_{1,\xi}-\frac{\varepsilon\mu}{3} u_{0, \tau\xi\xi} + \frac{\mu V}{3} u_{0, \xi\xi\xi}  + \frac{\varepsilon \mu V}{3} u_{1, \xi\xi\xi} + \zeta_{0,\xi}  + \varepsilon \zeta_{1, \xi} +\varepsilon u_0 u_{0,\xi}  \\
		& \quad +\varepsilon \mu \zeta_{0,\xi} u_{0,\xi\xi}  -\frac{\varepsilon\mu}{3} \left( u_0u_{0,\xi\xi\xi} -u_{0,\xi} u_{0,\xi\xi} -\zeta_{0,\xi yy} \right) +\frac{2 \varepsilon \mu }{3} \zeta_0 u_{0,\xi\xi\xi} = O(\varepsilon^2 \mu,\mu^2).
	\end{aligned} 
	\right.\nonumber
\end{equation}
From  the coefficients of $O(\varepsilon^0)$,  we have 
\begin{align}
	\left\{ 
	\begin{aligned}
		& -V \zeta_{0,\xi} + u_{0, \xi} =0,\\
		& \zeta_{0,\xi} + Vu_{0,\xi} =0. 
	\end{aligned}
	\right.\nonumber %\label{derivation eq-10}
\end{align}
For the existence of the non-trivial of the solution $(\zeta_0, u_0)$, the determinant of the coefficient matrix is required to be zero, that is $V^2=1$. Let $V=1$, one could derive the relation $u_0= \zeta_0$.  Next, we consider 
the coeffieients of $O(\varepsilon^1)$,  it is inferred from
\begin{align}
	& \zeta_{1, \xi} - u_{1, \xi} = \zeta_{0,\tau}  + (\zeta_0u_0)_{\xi} + \partial_{\xi}^{-1}u_{0,yy}, \nonumber \\
	& u_{1, \xi} - \zeta_{1, \xi} = u_{0,\tau}  + u_0u_{0,\xi}. \nonumber 
\end{align}
From the above inequalities, we obtain
\begin{align}
	\zeta_1 - u_1 = \frac{1}{4} u_0^2 + \frac{1}{2} \partial_{\xi}^{-2} u_{0,yy}.\nonumber
\end{align}  
Thus, we have     
\begin{align}
	\zeta =  \zeta_0 + \varepsilon \zeta_1 + O (\varepsilon^2)  
	= u + \frac{1}{4} \varepsilon u^2 + \frac{1}{2} \varepsilon \partial_{\xi}^{-2} u_{yy}   + O (\varepsilon^2)  . \label{derivation eq-3}
\end{align}
Taking (\ref{derivation eq-3}) into (\ref{derivation eq-2}), we have 
\begin{align}
	&  \varepsilon u_{\tau}  -\frac{\varepsilon \mu}{3} u_{\tau \xi\xi} + \frac{3\varepsilon}{2} uu_{\xi}  +\frac{\varepsilon}{2} \partial_{\xi}^{-1} u_{yy} +\frac{\mu }{3} u_{\xi\xi\xi} + \varepsilon^2 u_y \partial_{\xi}^{-1} u_y \nonumber \\
	 & \quad +\frac{\varepsilon\mu}{3} \left( uu_{\xi\xi} + \frac{3}{2} u_{\xi}^2 + u_{yy} \right)_x   = O(\varepsilon^2 \mu, \mu^2, \varepsilon^3) . \nonumber
	\end{align} 
Back to the original variables and set $u(t,x,y) = u(\tau, \xi,y)+\frac{2(1-\kappa)}{3\varepsilon } $, which satisfies the following eqution
\begin{align}
	 u_t -\frac{\mu}{3} u_{txx} +\kappa u_x  +\frac{3\varepsilon}{2}  uu_x +\frac{\varepsilon}{2} \partial_x^{-1} u_{yy} + \varepsilon^2 u_y \partial_x^{-1} u_y =O(\varepsilon \mu, \mu^2, \varepsilon^3) \nonumber
\end{align}
Now we define $m= u-\frac{\mu}{3} u_{xx}$, for some constant $\sigma$, we have 
\begin{align}
	\frac{3}{2}  \varepsilon uu_x= \frac{\sigma}{2}  \varepsilon \left( 2mu_x + um_x \right) +\frac{3}{2} (1-\sigma) \varepsilon uu_x + O (\varepsilon \mu) ,\nonumber
\end{align}
which gives rise to 
\begin{align}
	 m_t + \kappa u_x + \frac{\sigma}{2} \varepsilon \left( 2mu_x + um_x  \right)  +\frac{3}{2} \left( 1-\sigma \right) \varepsilon uu_x +\frac{\varepsilon}{2} \partial_x^{-1} u_{yy} + \varepsilon^2 u_y \partial_x^{-1} u_y = O (\varepsilon \mu, \mu^2, \varepsilon^3). \nonumber
\end{align}
If we consider the CH regime $\varepsilon= O(\mu^{\frac{1}{2}})$, then $ \varepsilon^3= O(\mu^{\frac{3}{2}})$ and $ \varepsilon \mu =O (\mu^{\frac{3}{2}})$. Now we focus on the model equation with ignoring terms of order $O (\varepsilon^2, \mu)$,
\begin{align}
	 \left( m_t + \kappa u_x +  \frac{\sigma}{2} \varepsilon \left( 2mu_x + um_x  \right)  +\frac{3}{2} \left( 1-\sigma \right) \varepsilon uu_x\right)_x  +\frac{\varepsilon}{2}  u_{yy} = O (\varepsilon^2, \mu)\nonumber
\end{align}
Then applying the scaling 
\begin{align}
	u\to \frac{2}{\varepsilon} u , \quad x\to \sqrt{\frac{\mu}{3}} x, \quad t \to \sqrt{\frac{\mu}{3}} t, \quad  y \to \sqrt{\frac{\varepsilon\mu}{6}}y,\nonumber
\end{align} 
and ignoring terms of oreder $O(\varepsilon^2,\mu)$,
the following rod-KP equation is established, 
\begin{align}
	\left( m_t +\kappa u_x + \sigma \left( 2mu_x + um_x\right) + 3(1-\sigma) uu_x   \right)_x + u_{yy} =0, \label{Rod-KP equation form}
\end{align}
where the momentum density $m= u-u_{xx}$, and  $u(t,x,y)$ is the verage fluid velocity in the $x$-direction and $y$-direction, and $u$ satisfy the far field conditions $u\to 0$ as $|x| \to \infty$.    The real dimensionless constant $\sigma$ is a  parameter which shows a balance between nonlinear steepening and amplification in fluid convection due to stretching. According to the value of $\sigma$, the derived rod-KP equation exhibits two special models. When $\sigma=0$, the rod-KP equation (\ref{Rod-KP equation form}) reduces the BBM-KP equation (\ref{BBM-KP equation}); when $\sigma =1$, the rod-KP equation (\ref{Rod-KP equation form}) recovers the CH-KP equation (\ref{CH-KP eqaution}).

%Let $\Omega_{\lambda} = \left\lbrace (x,y) \in \mathbb{R}^2 , x\geq \lambda  \right\rbrace $,  $\Omega_{\lambda}' =  \left\lbrace (x,y) \in \mathbb{R}^2 , x<  \lambda  \right\rbrace$, $x_{\lambda} = 2\lambda -x$ and $ \varphi_{\lambda} = \varphi (x_{\lambda}, y)$.  We have the following lemma.  It is easy from (\ref{symmetry eq-4}) and the parity of $K$ in $x$ to see that 	\begin{align 	\varphi (x,y) =&  \left( \int_{\Omega_{\lambda}}+ \int_{\Omega_{\lambda}'} \right) K(x-x_1,y-y_1)  f(\varphi) (x_1,y_1) dx_1 dy_1 \nonumber \\ 		= & \int_{\Omega_{\lambda}} \left( K(x-x_1,y-y_1) f(\varphi) (x_1,y_1) + K(x_{\lambda} -x_1,y-y_1) f(\varphi_{\lambda } ) (x_1,y_1) \right) dx_1dy_1, \label{symmetry lemma_1  eq-1} 	\end{align} and  \begin{align 	\varphi_{\lambda} (x,y) =& \varphi (2\lambda-x,y) = \left( K \ast f(\varphi) \right)(2\lambda-x,y) \nonumber \\ 	=  & \left( \int_{\Omega_{\lambda}} + \int_{\Omega_{\lambda}'} \right) K(2\lambda -x -x_1,y-y_1) f(\varphi) (x_1,y_1) dx_1dy_1\nonumber \\ = &\int_{\Omega_{\lambda}} K(x_{\lambda}-x_1,y-y_1) f(\varphi) (x_1,y_1) + K(x_1-x,y-y_1) f(\varphi_{\lambda}) (x_1,y_1)dx_1dy_1.\label{symmetry lemma_1  eq-2} \end{align} It is inferred that  \begin{align} & \varphi (x,y) - \varphi_{\lambda}(x,y)\nonumber \ 	 = & \int_{\Omega_{\lambda}} \left( K(x-x_1,y-y_1) -K(x_{\lambda} -x_1,y-y_1) \right) \left(  f(\varphi)(x_1,y_1) - f(\varphi_{\lambda}) (x_1,y_1) \right) dx_1 dy_1.  \end{align} \end{proof}
\section{Local well-posedness and blow-up criterion}\label{Local Well-posedness}
	\newtheorem {remark2}{Remark}[section]
	\newtheorem{theorem2}{Theorem}[section]
	\newtheorem{lemma2}{Lemma}[section]
	\newtheorem{proposition2}{Proposition}[section]
	\newtheorem{corollary2}{Corollary}[section]
	In this section, we consider the issues of local well-posedness (Theorem \ref{LWP theorem}) and blow-up criterion  (Theorems \ref{Blow-up criterion}-\ref{BBM blow-up criterion}) for  the Cauchy problem (\ref{the nonlocal weak form}). 
\subsection{Local well-posedness}
 In order to prove Theorem \ref{LWP theorem}, we need the following useful lemmas
\begin{lemma2}[\cite{Moser ASNSP 1966}]\label{Moser-type Estimate} 
	Let $\Lambda^s \triangleq (1-\Delta_{x,y})^{\frac{s}{2}}$ with $s>0$, $p_1, \, p_2 \in [2,+\infty)$, $q_1, \, q_2 \in (2,+\infty]$, $\frac{1}{p_1} + \frac{1}{q_1}  = \frac{1}{p_2} + \frac{1}{q_2} =\frac{1}{2} $. Then there holds 
	\begin{equation}
		\Vert fg \Vert_{H^s(\mathbb{R}^2)} \leq C \left(  \Vert \Lambda^s f \Vert_{L^{p_1} (\mathbb{R}^2)} \Vert g \Vert_{L^{q_1}(\mathbb{R}^2)}   + \Vert \Lambda^s g \Vert_{L^{p_2}(\mathbb{R}^2)} \Vert f \Vert_{L^{q_2}(\mathbb{R}^2)} \right) ,
		\label{Moser-type estimate}
	\end{equation}
where the constant $C$ is independent of $f$ and $g$.
\end{lemma2}

\renewcommand{\theenumi}{\roman{enumi}}
\renewcommand{\labelenumi}{(\theenumi)}
\begin{lemma2}[\cite{Kato-Ponce CPAM 1988}]\label{Commutator estimates}
	Let $\Lambda^s \triangleq  (1-\Delta_{x,y})^{\frac{s}{2}}$ with $s>0$. Then we have the following commutator estimates:
	\begin{enumerate}
		\item $\Vert [\Lambda^s, f] g \Vert_{L^2(\mathbb{R}^2)} \leq C \left(  \Vert \Lambda^s f \Vert_{L^2(\mathbb{R}^2)} \Vert g \Vert_{L^{\infty}(\mathbb{R}^2)} + \Vert \nabla f \Vert_{L^{\infty}( \mathbb{R}^2)} \Vert \Lambda^{s-1} g \Vert_{L^2(\mathbb{R}^2)}       \right) $;
		\item $\Vert [\Lambda^s, f] g \Vert_{L^2(\mathbb{R}^2)} \leq  C \Vert \nabla f \Vert_{H^{q_0}(\mathbb{R}^2)}  \Vert g  \Vert_{H^{s-1}(\mathbb{R}^2)} $, $ \forall q_0 >1$, $ 0 \leq s \leq q_0 +1$, 
	\end{enumerate}
where all the constants $C$ are independent of $f$ and $g$.
\end{lemma2}
%\begin{lemma2}[\cite{Gui-Liu-Luo-Yin AIM 2021}]	For $s\geq 2$, the space $X^s(\mathbb{R}^2)$ is continuously embedded in $Lip (\mathbb{R}^2)$.  \end{lemma2}
\begin{lemma2}[\cite{Gui-Liu-Luo-Yin AIM 2021,Molinet-Saut-Tzvetkov CMP 2007}]\label{Sobolev inequality lemma}
	Let $u \in X^s(\mathbb{R}^2)$, $ s\geq 2$. Then the following estimates hold:
	\begin{enumerate}
		\item \label{Sobolev inequality lemma-2}  $\Vert \nabla u \Vert_{L^{\infty}(\mathbb{R}^2)}^2 \leq  C \Vert u \Vert_{X^s(\mathbb{R}^2)}^2$,
		\item \label{Sobolev inequality lemma-3}$\Vert u \Vert_{L^{\infty}(\mathbb{R}^2)}\leq C\left(\Vert u \Vert_{L^{\infty}(\mathbb{R}^2)} + \Vert u_x \Vert_{L^2(\mathbb{R}^2)} + \Vert u_y \Vert_{L^{\infty}(\mathbb{R}^2)} \right) $,
		\item\label{Sobolev inequality lemma-1} $\Vert u \Vert_{L^{\infty}(\mathbb{R}^2)}^2 \leq C \Vert u \Vert_{L^2(\mathbb{R}^2)}^{\frac{1}{2}}  \Vert u_x \Vert_{L^2(\mathbb{R}^2)}^{\frac{1}{2}} \Vert u_y \Vert_{L^2(\mathbb{R}^2)}^{\frac{1}{2}} \Vert u_{xy} \Vert_{L^2(\mathbb{R}^2)}^{\frac{1}{2}} $,
		\item \label{Sobolev inequality lemma-4}$\Vert u \Vert_{L^3(\mathbb{R}^2)}^3 \leq C \Vert u \Vert_{L^2(\mathbb{R}^2)}^{\frac{3}{2}} \Vert u_x \Vert_{L^2(\mathbb{R}^2)} \Vert \partial_x^{-1} u_y \Vert_{L^2(\mathbb{R}^2)}^{\frac{1}{2}}$,
	\end{enumerate}
for some positive constant $C$ independent of $u$. 
\end{lemma2}
\begin{lemma2}\label{Sobolev embedding remark-2}
	Let $ u, \, v \in H^1(\mathbb{R}^2)$, the following  inequality holds  true:
	\begin{align}
		\Vert uv \Vert_{L^2(\mathbb{R}^2)} \leq \Vert u \Vert_{L^2(\mathbb{R}^2)}^{\frac{1}{2}} \Vert u_x \Vert_{L^2(\mathbb{R}^2)}^{\frac{1}{2}}\Vert v \Vert_{L^2(\mathbb{R}^2)}^{\frac{1}{2}}
		\Vert v_y \Vert_{L^2(\mathbb{R}^2)}^{\frac{1}{2}}. 
		\label{Sobolev embedding remark-2 eq}
	\end{align}
\end{lemma2}
	\begin{proof}
		For all $u, \, v\in H^1(\mathbb{R}^2)$, the following estimate  holds:
		\begin{align}
			\Vert uv \Vert_{L^2} \leq \Vert u \Vert_{L_x^{\infty} L_y^2} \Vert v \Vert_{L_x^2 L_y^{\infty}}. 
			\label{Sobolev embedding remark-2 proof eq}
		\end{align}
		Thanks to $u , \, v\in H^1(\mathbb{R}^2)$, we get  \[ u^2(x,y) = 2 \int_{-\infty}^x u(x',y)u_x (x' ,y) dx', \quad v^2(x,y) =2 \int_{-\infty}^y v (x,y') v_y(x,y') dy'. \] 
		Thus,  by the Cauchy-Schwartz inequality, we have 
		\begin{align}
			\Vert u \Vert_{L_x^{\infty}}^2 \leq 2  \Vert u \Vert_{L_x^2} \Vert u_x \Vert_{L_x^2},  \quad 
			\Vert v\Vert_{L_x^2}^2  \leq  2 \Vert v \Vert_{L^2_{x,y} } \Vert v_y \Vert_{L_{x,y}^2},\nonumber 
		\end{align}
		which along with (\ref{Sobolev embedding remark-2 proof eq})   completes the proof of Lemma \ref{Sobolev embedding remark-2}.
	\end{proof}
%\begin{lemma2}\label{complex interpolation} Suppose that $s_1, s_2\in \mathbb{R}$, we have the following complex interpolation:  	\begin{align} 	\Vert f \Vert_{H^{\theta s_1 + (1-\theta) s_2 }} \leq \Vert f \Vert_{H^{s_1}}^{\theta} \Vert f \Vert_{H^{s_2}}^{1-\theta}, \quad \forall  u \in H^{s_1} \cap H^{s_2}, \quad  \forall \theta \in [0,1]. \nonumber  	\end{align}\end{lemma2}

%With these lemmas in hand, we are now in the position to prove Theorem \ref{LWP theorem}

\begin{proof}[Proof of Theorem \ref{LWP theorem}]
	Firstly, applying the operator $\Lambda^s \partial_x^{-1}$ to the equation (\ref{Rod-KP equation form-1}), we have 
	\begin{align}
		\partial_t \Lambda^s u - \partial_t \Lambda^su_{xx} + \left( \kappa \Lambda u + \frac{3}{2} \Lambda^s \left(  u^2 \right) - \frac{\sigma}{2}  \Lambda^s (u_x^2)  - \sigma\Lambda^s (uu_{xx})  \right)_x + \Lambda^s \partial_x^{-1} u_{yy} =0 . \label{Rod-KP LWP eq-1}
	\end{align}
	Applying the operator $\partial_x^{-1}$ to (\ref{Rod-KP LWP eq-1}),  we have 
	 \begin{align} 
	 		\partial_t\Lambda^s \partial_x^{-1} u - \partial_t \Lambda^s u_x +\kappa \Lambda^s u + \frac{3}{2} \Lambda^s (u^2) -\sigma \Lambda^s (uu_x)_x +\frac{\sigma }{2} \Lambda^s (u_x^2) + \Lambda^s \partial_x^{-2} u_{yy} =0. \label{Rod-KP LWP eq-2}	
	 \end{align}
	Then we take the $L^2$ inner-product between (\ref{Rod-KP LWP eq-2}) and $\Lambda^s \partial_x^{-1} u$ to get  
		\begin{align} 
			& \frac{d}{dt} \frac{1}{2} \left( \Vert \partial_x^{-1} u \Vert_{H^s}^2 + \Vert u \Vert_{H^s}^2  \right)  	\nonumber \\	= & \kappa \int_{\mathbb{R}^2} \Lambda^s u \Lambda^s \partial_x^{-1} u dxdy - \frac{3}{2} \int_{\mathbb{R}^2} \Lambda^s (u^2) \Lambda^s \partial_x^{-1}u dxdy + \sigma \int_{\mathbb{R}^2} \Lambda^s (uu_x)_x \Lambda^s \partial_x^{-1} u dxdy \nonumber \\ 	& -\frac{\sigma}{2} \int_{\mathbb{R}^2}  \Lambda^s(u_x^2) \Lambda^s \partial_x^{-1} u dxdy  -\int_{\mathbb{R}^2} \Lambda^s \partial_x^{-2} u_{yy}\Lambda^s  \partial_x^{-1}  dxdy .\label{Rod-KP LWP eq-3}
			\end{align}
	Integrating by parts, we have  
		\begin{align} 
				\int_{\mathbb{R}^2} \Lambda^s u \Lambda^s \partial_x^{-1} u dxdy =0,	\label{Rod-KP LWP eq-4} 
				\end{align} 
			and 
				\begin{align} 
					\int_{\mathbb{R}^2} \Lambda^s \partial_x^{-2} u_{yy}\Lambda^s  \partial_x^{-1}  dxdy = - \int_{\mathbb{R}^2} \Lambda^s \partial_x^{-2} u_{y}\Lambda^s  \partial_x^{-1} u_y  dxdy =0 .  	\label{Rod-KP LWP eq-4-1} 
					\end{align} 	
				Thanks to the H\"{o}lder inequality, we have
	\begin{align} 	
		\left| \int_{\mathbb{R}^2}\Lambda^s u \Lambda^s \partial_x^{-1} u dxdy  \right| \leq \Vert u \Vert_{H^s} \Vert \partial_x^{-1} u\Vert_{H^s}, \label{Rod-KP LWP eq-5}	
	\end{align}
\begin{align}
		\left|  \int_{\mathbb{R}^2} \Lambda^s (u^2) \Lambda^s \partial_x^{-1} u dxdy \right| \leq \Vert u^2 \Vert_{H^s} \Vert \partial_x^{-1} u \Vert_{H^s} \leq C \Vert u \Vert_{L^{\infty}} \Vert u \Vert_{H^s} \Vert\partial_x^{-1} \Vert_{H^s} , \label{Rod-KP LWP eq-5.1} 	
	\end{align}
\begin{align} 	
	\left| 	\int_{\mathbb{R}^2} \Lambda^s (uu_x)_x \Lambda^s \partial_x^{-1} u dxdy\right| =&  \frac{1}{2} \left| \int_{\mathbb{R}^2} \Lambda^s (u^2) \Lambda^s u_x dxdy \right| \leq \frac{1}{2} \Vert u^2 \Vert_{H^s} \Vert u_x \Vert_{H^s} \nonumber \\	\leq & \frac{1}{2} C \Vert u \Vert_{L^{\infty}} \Vert u \Vert_{H^s} \Vert u_x \Vert_{H^s},  \label{Rod-KP LWP eq-6} 
	\end{align}
 and 
\begin{align}	
	\left| 	\int_{\mathbb{R}^2} \Lambda^s (u_x^2) \Lambda^s \partial_x^{-1} udxdy \right| \leq \Vert u_x^2 \Vert_{H^s} \Vert \partial_x^{-1} u \Vert_{H^s} \leq C\Vert u_x\Vert_{L^{\infty}} \Vert u_x \Vert_{H^s} \Vert \partial_x^{-1} u \Vert_{H^s}. \label{Rod-KP LWP eq-7}
	\end{align}
 Plugging (\ref{Rod-KP LWP eq-4})-(\ref{Rod-KP LWP eq-7}) into  (\ref{Rod-KP LWP eq-3}) to yield that for $s\geq 2$ 
 \begin{align} 
 	& \frac{d}{dt} \left( \Vert \partial_x^{-1} u \Vert_{H^s}^2 + \Vert u \Vert_{H^s}^2 \right)  	\nonumber\\ 	\leq & C \left( \Vert u\Vert_{L^{\infty}} +\Vert u_x \Vert_{L^{\infty}} \right)  \left( \Vert \partial_x^{-1} u \Vert_{H^s}^2 + \Vert u \Vert_{H^s}^2 + \Vert u_x \Vert_{H^s}^2 \right) . \label{Rod-KP LWP eq-8}	
 \end{align}
	In order to get estimate of $\Vert \Lambda^s u  \Vert_{L^2}^2$, we first consider  the case $s>2$, thake the $L^2$ inner-product between (\ref{Rod-KP LWP eq-1}) and $\Lambda^s u$ to give 
	\begin{align}
		&\frac{1}{2} \frac{d}{dt} \left(  \Vert u\Vert_{H^s}^2 + \Vert u_x \Vert_{H^s}^2 \right) \nonumber \\
		= & - \kappa \int_{\mathbb{R}^2} \Lambda^s u_x \Lambda u dxdy +\frac{3}{2} \int_{\mathbb{R}^2} \Lambda^s (u^2) \Lambda^s u_x dxdy -\frac{\sigma}{2} \int_{\mathbb{R}^2} \Lambda^s (u_x^2) \Lambda^s u_x dxdy \nonumber \\
		& -\sigma \int_{\mathbb{R}^2} \Lambda^s (uu_{xx}) \Lambda^s u_x dxdy -\int_{\mathbb{R}^2} \partial_x^{-1} u_{yy} \Lambda^s udxdy . \label{Rod-KP LWP eq-9}
	\end{align}
	Thanks to integration by parts, we get
	\begin{align}
		\int_{\mathbb{R}^2} \Lambda^s u_x \Lambda u dxdy=0,\label{Rod-KP LWP eq-10-1}
	\end{align}
and  
	\begin{align}
		\int_{\mathbb{R}^2} \Lambda^s \partial_x^{-1} u_{yy} \Lambda^s u dxdy = - \int_{\mathbb{R}^2} \Lambda^s (\partial_x^{-1} u_y)  \Lambda^s u_y dxdy =0. \label{Rod-KP LWP eq-10}
	\end{align}
	On the other hand, one can see
	\begin{align}
		\left|  \int_{\mathbb{R}^2} \Lambda^s (u^2) \Lambda^s u_x dxdy \right| \leq \Vert u^2 \Vert_{H^s} \Vert u_x \Vert_{H^s} \leq C\Vert u\Vert_{L^{\infty}} \Vert u \Vert_{H^s} \Vert u_x \Vert_{H^s} ,  \label{Rod-KP LWP eq-11}
	\end{align}
	and 
	\begin{align}
		\left| \int_{\mathbb{R}^2} \Lambda^s (u_x^2) \Lambda^s u_x dxdy  \right| \leq \Vert u_x^2 \Vert_{H^s} \Vert u_x \Vert_{H^s} \leq C\Vert u_x \Vert_{L^{\infty}}  \Vert u_x \Vert_{H^s}^2 .  \label{Rod-KP LWP eq-12}
	\end{align}
	In view of Lemma \ref{Commutator estimates}, we get
	\begin{align}
		& \left| 	\int_{\mathbb{R}^2} \Lambda^s (uu_{xx}) \Lambda^s u_x dxdy \right| \nonumber \\
		= &  \left|-  \frac{1}{2} \int_{\mathbb{R}^2}  u_x \left( \Lambda^s u_x\right)^2 dxdy + \int_{\mathbb{R}^2} [\Lambda^s, u] u_{xx} \Lambda^s u_x dxdy \right| \nonumber \\
		\leq & \frac{1}{2} \Vert u \Vert_{L^{\infty}} \Vert u_x \Vert_{H^s}^2 + C \left(\Vert \nabla u \Vert_{L^{\infty}} \Vert u_{xx} \Vert_{H^{s-1}} + \Vert u_{xx} \Vert_{L^{\infty}} \Vert u \Vert_{H^s}  \right) \Vert u_x \Vert_{H^s}. \label{Rod-KP LWP eq-13}
	\end{align}
	 Therefore, plugging (\ref{Rod-KP LWP eq-10})-(\ref{Rod-KP LWP eq-13}), we obtain
	\begin{align}
		& \frac{d}{dt} \left( \Vert u\Vert_{H^s}^2 + \Vert u_x \Vert_{H^s}^2  \right) \nonumber \\
		\leq &  C \left(  \Vert u\Vert_{L^{\infty}}  + \Vert \nabla u \Vert_{L^{\infty}} +\Vert u_{xx} \Vert_{L^{\infty}}  \right) \left( \Vert u \Vert_{H^s}^2+ \Vert u_x \Vert_{H^s}^2 \right) . \label{Rod-KP LWP eq-14}
	\end{align}
Combining (\ref{Rod-KP LWP eq-14}) with (\ref{Rod-KP LWP eq-8}) implies 
\begin{align}
	\frac{d}{dt} \left( 2\Vert u \Vert_{H^s}^2 + \Vert \partial_x^{-1} u \Vert_{H^s}^2 + \Vert u_x \Vert_{H^s}^2 \right) \leq&  C \left(\Vert u \Vert_{L^{\infty}} + \Vert u_x \Vert_{L^{\infty}} + \Vert u_{xx} \Vert_{L^{\infty}} \right) \nonumber \\
& 	\times \left( \Vert u \Vert_{H^s}^2 + \Vert \partial_x^{-1} u \Vert_{H^s}^2 + \Vert u_x \Vert_{H^s}^2 \right). \label{Rod-KP LWP eq-14.1}
\end{align}
	Next we consider the case $s=2$, acting the operators $\nabla$, $ \nabla \partial_x$ and $\nabla \partial_x^2$ to (\ref{the nonlocal weak form}) respectively gives rise to 
	\begin{align}
		& 	\partial_t \nabla u + \sigma \left( u \nabla u_x + u_x \nabla u\right) + p_x \ast \left(\left(3-\sigma\right) u \nabla u + \sigma u_x \nabla u_x  +\kappa \nabla u  \right)  + p \ast  \partial_x^{-1} \nabla u_{yy}=0, \label{Rod-KP LWP case s=2 eq-1} \\
		&\partial_t \nabla u_x + \sigma \left(  u\nabla u_{xx} + u_x \nabla u_x + u_{xx} \nabla u \right) 
		-(3-\sigma) u \nabla u-\kappa \nabla u \nonumber \\
		&\quad \quad \quad \quad + p \ast \left(\left( 3-\sigma \right) u\nabla u + \sigma u_x \nabla u_x  +\kappa \nabla u \right) + p \ast \nabla u_{yy} =0, \label{Rod-KP LWP case s=2 eq-2}\\
		& \partial_t \nabla u_{xx} +\sigma\left(  u\nabla u_{xxx} + 2u_x \nabla u_{xx} + 2u_{xx} \nabla u_x + u_{xxx} \nabla u  \right)  - (3-\sigma) \left( u_x \nabla u + u\nabla u_x  \right) -\kappa \nabla u_x  \nonumber \\
		& \quad \quad \quad \quad + p_x  \ast \left(\left( 3-\sigma \right) u\nabla u + \sigma u_x \nabla u_x+\kappa \nabla u_x  \right) + p \ast \nabla u_{xyy}=0. \label{Rod-KP LWP case s=2 eq-3} 
	\end{align}
	Taking the $L^2$ inner-product between three equations in (\ref{Rod-KP LWP case s=2 eq-1})-(\ref{Rod-KP LWP case s=2 eq-3}) and $\nabla u$, $ \nabla u_x$, $\nabla u_{xx}$ respectively, and summating them, we get
	\begin{align}
		\frac{1}{2} \frac{d}{dt} \left(  \Vert \nabla u \Vert_{L^2}^2 + \Vert \nabla u_x \Vert_{L^2}^2 \right) = \sum\limits_{i=1}^4 I_i,\label{Rod-KP LWP case s=2 eq-4}
	\end{align}
	\begin{align}
		\frac{1}{2} \frac{d}{dt} \left(  \Vert \nabla u \Vert_{L^2}^2 + \Vert \nabla u_x \Vert_{L^2}^2  + \Vert \nabla u_{xx}\Vert_{L^2}^2\right) =\sum\limits_{i=1}^7 I_i,\label{Rod-KP LWP case s=2 eq-5}
	\end{align}
	where 
	\begin{align}
		I_1 :=&  -\sigma \int_{\mathbb{R}^2}  \left( u \nabla u_x + u_x \nabla u  \right)  \cdot \nabla u  dxdy , \nonumber \\
		I_2 := & \int_{\mathbb{R}^2} -\left( \sigma \left(  u\nabla u_{xx} + u_x \nabla u_x + u_{xx} \nabla u \right) 
		-(3-\sigma) u \nabla u  - \kappa  \nabla u \right)\cdot  \nabla u_x    dxdy, \nonumber\\
		I_3 := & -\int_{\mathbb{R}^2} p_x \ast \left(\left(3-\sigma\right) u \nabla u + \gamma u_x \nabla u_x +\kappa u   \right) \cdot \nabla u  \nonumber \\
		& \quad \quad + p \ast \left(\left(3-\sigma\right) u \nabla u + \gamma u_x \nabla u_x  + \kappa u   \right) \cdot \nabla u_x dxdy ,\nonumber \\  
		I_4:=&  - \int_{\mathbb{R}^2} p \ast \partial_x^{-1} \nabla u_{yy}  \cdot \nabla u + p \ast \nabla u_{yy}  \cdot \nabla u_x dxdy,	 \nonumber \\
		I_5:= & -\int_{\mathbb{R}^2} (\sigma\left(  u\nabla u_{xxx} + 2u_x \nabla u_{xx} + 2u_{xx} \nabla u_x + u_{xxx} \nabla u  \right) \nonumber \\
		&\quad \quad  - (3-\sigma) \left( u_x \nabla u + u\nabla u_x  \right) -\kappa \nabla u_x ) \cdot \nabla u_{xx} dxdy,  \nonumber 
		\\
		I_6:= & -\int_{\mathbb{R}^2} p_x  \ast \left(\left( 3-\sigma \right) u\nabla u + \sigma u_x \nabla u_x + \kappa \nabla u  \right) \cdot \nabla u_{xx}  dxdy,  \nonumber \\
		I_7:= & -\int_{\mathbb{R}^2} p\ast  \nabla u_{xyy} \cdot \nabla u_{xx}dxdy . \nonumber 
	\end{align}
	Thanks to integration by parts, we get $I_3=0$, $I_4=0$, $I_7=0$, and 
	\begin{align}
		\left| I_1  \right| \leq \frac{|\sigma|}{2} \Vert u_x \Vert_{L^{\infty}} \Vert \nabla u \Vert_{L^2}^2,
		\label{Rod-KP LWP case s=2 eq-6}
	\end{align}
	\begin{align}
		\left| I_2  \right|  \leq & \frac{|\sigma|}{2} \Vert u_x \Vert_{L^{\infty}} \Vert \nabla u_x \Vert_{L^2}^2 + |\sigma |\Vert \nabla u \Vert_{L^{\infty}} \Vert \nabla u_x \Vert_{L^2}^2 + |3-\sigma| \Vert \nabla u \Vert_{L^{\infty}} \Vert u \Vert_{L^2} \Vert \nabla u_x \Vert_{L^2}  \nonumber \\
		\leq & C \left(  \Vert u_x \Vert_{L^{\infty}} +\Vert \nabla u \Vert_{L^{\infty}} \right) \left(  \Vert u \Vert_{L^2}^2 + \Vert \nabla u_x \Vert_{L^2}^2\right) .\label{Rod-KP LWP case s=2 eq-7}
	\end{align}
	It is noticed that 
	\begin{align}
		\left|  I_5 \right| \leq & \frac{3|\sigma|}{2} \Vert u_x \Vert_{L^{\infty}} \Vert \nabla u_{xx} \Vert_{L^2}^2 + 2|\sigma | \Vert u_{xx} \nabla u_x \Vert_{L^2} \Vert \nabla u_{xx} \Vert_{L^2} \nonumber \\
		& + \left| \sigma \right| \Vert \nabla u \Vert_{L^{\infty}} \Vert \nabla u \Vert_{L^2} \Vert \nabla u_{xx} \Vert_{L^2}  + (3-\sigma) \Vert \nabla u \Vert_{L^{\infty}}  \Vert \nabla u \Vert_{L^2} \Vert \nabla u_{xx} \Vert_{L^2}  \nonumber \\
		& +\frac{\left| 3-\sigma\right| }{2} \Vert u_x\Vert_{L^{\infty}} \Vert \nabla u_x\Vert_{L^2}^2  . \label{Rod-KP LWP case s=2 eq-8}
	\end{align}
	Thanks to Lemma \ref{Sobolev embedding remark-2}, we have 
	\begin{align}
		\Vert u_{xx} \nabla u_x \Vert_{L^2} \leq \Vert u_{xx} \Vert_{L^2}^{\frac{1}{2}} \Vert u_{xxy} \Vert_{L^2}^{\frac{1}{2}} \Vert \nabla u_x \Vert_{L^2}^{\frac{1}{2}} \Vert \nabla u_{xx} \Vert_{L^2}^{\frac{1}{2}}, \nonumber 
	\end{align}
	which along with 
	\begin{align}
		\left|  I_5 \right| \leq C \left( \Vert \nabla u \Vert_{L^{\infty}} + \Vert \nabla u_x \Vert_{L^2} \right)  \left( \Vert  \nabla u \Vert_{L^2}^2 + \Vert \nabla u_x \Vert_{L^2}^2 + \Vert \nabla u_x \Vert_{L^2}^2 \right) .\label{Rod-KP LWP case s=2 eq-9}
	\end{align} 
	In view of the H\"{o}lder inequality, we have 
	\begin{align}
		\left|  I_6 \right| \leq \left| 3-\sigma \right|  \Vert \nabla u \Vert_{L^{\infty}} \Vert u \Vert_{L^2} \Vert \nabla u_{xx} \Vert_{L^2} + \left| \sigma \right|  \Vert u_x \Vert_{L^{\infty}} \Vert \nabla u_x \Vert_{L^2} \Vert u_{xx} \Vert_{L^2}.
		\label{Rod-KP LWP case s=2 eq-10}
	\end{align}
	Plugging (\ref{Rod-KP LWP case s=2 eq-6}) and (\ref{Rod-KP LWP case s=2 eq-7}) into (\ref{Rod-KP LWP case s=2 eq-4}), we have 
	\begin{align}
	\frac{d}{dt} \left( \Vert \nabla u \Vert_{L^2}^2 + \Vert \nabla u_x \Vert_{L^2}^2  \right) 
		\leq   C  \Vert \nabla u \Vert_{L^{\infty}} \left(  \Vert u \Vert_{L^2}^2 + \Vert \nabla u \Vert_{L^2}^2 + \Vert \nabla u_{xx} \Vert_{L^2}^2 \right) ,	\label{Rod-KP LWP case s=2 eq-11}
	\end{align}
	and 
	\begin{align}
		& \frac{d}{dt} \left( \Vert \nabla u \Vert_{L^2}^2 + \Vert \nabla u_x \Vert_{L^2}^2 + \Vert \nabla u_{xx} \Vert_{L^2}^2   \right)  \nonumber \\
		\leq  & C \left( \Vert \nabla u \Vert_{L^{\infty}} + \Vert \nabla u_x \Vert_{L^2} \right) \left(  \Vert u \Vert_{L^2}^2 + \Vert \nabla u \Vert_{L^2}^2 +\Vert \nabla u_x \Vert_{L^2}^2 + \Vert \nabla u_{xx} \Vert_{L^2}^2 \right) . \label{Rod-KP LWP case s=2 eq-12}
	\end{align}
	Applying the operators $\partial_y^2$ and $\partial_x \partial_y^2$ to (\ref{the nonlocal weak form}) respectively gives rise to 
	\begin{align}
		& 	\partial_t u_{yy} +\sigma\left(  u u_{xyy} + 2u_y u_{xy} + u_x u_{xy}\right)  \nonumber \\
		&\quad \quad + p_x \ast \left( \left(3-\sigma \right) \left(uu_{yy} + u_y^2  \right)  +\sigma \left(  u_xu_{xyy} + u_{xy}^2 \right)+ \kappa u_{yy} \right)  + p \ast \partial_x^{-1} u_{yyyy} =0, \label{Rod-KP LWP case s=2 eq-13} \\
		& \partial_t u_{xyy} + \sigma \left( uu_{xxyy} +u_{xy}^2 + 2u_y u_{xyy} + u_{xx} u_{xy}  + u_xu_{xxy} \right)  -(3-\sigma) ( uu_{yy} + u_y^2)  -\kappa u_{yy}  \nonumber \\
		& \quad \quad + p  \ast \left( \left(3-\sigma \right) \left(uu_{yy} + u_y^2  \right)  + \sigma \left(  u_xu_{xyy} + u_{xy}^2 \right) + \kappa u_{yy} \right)  + p \ast  u_{yyyy} =0. \label{Rod-KP LWP case s=2 eq-14} 
	\end{align}
	Taking the $L^2$ inner product between two equations in (\ref{Rod-KP LWP case s=2 eq-13}), (\ref{Rod-KP LWP case s=2 eq-14}) and $u_{yy}$, $u_{xyy}$ respectively, and summating them, we get 
	\begin{align}
		\frac{1}{2} \frac{d}{dt} \left( \Vert u_{yy} \Vert_{L^2}^2 + \Vert u_{xyy} \Vert_{L^2}^2  \right) =\sum\limits_{i=1}^4 J_i, \label{Rod-KP LWP case s=2 eq-15} 
	\end{align}
	with 
	\begin{align}
		& J_1 := - \int_{\mathbb{R}^2} \sigma \left( uu_{xyy} + 2u_y u_{xy} + u_x u_{xy} \right) u_{yy} dxdy , \nonumber \\
		& J_2 :=-  \int_{\mathbb{R}^2}  \sigma \left( uu_{xxyy} +u_{xy}^2 + 2u_y u_{xyy} + u_{xx} u_{xy}  + u_xu_{xxy} \right)  u_{xyy} \nonumber \\
		& \quad \quad \quad \quad \quad -(3-\sigma) ( uu_{yy} + u_y^2) u_{xyy} - \kappa u_{yy} u_{xyy} dxdy \nonumber \\
		& J_3:= -\int_{\mathbb{R}^2}  p_x \ast \left( \left(3-\sigma \right) \left(uu_{yy} + u_y^2  \right)  +\sigma\left(  u_xu_{xyy} + u_{xy}^2 \right) +\kappa u_{yy} \right)  u_{yy}  \nonumber \\
		& \quad \quad \quad \quad \quad + p \ast \left( \left(3-\sigma\right) \left(uu_{yy} + u_y^2  \right)  + \sigma \left(  u_xu_{xyy} + u_{xy}^2 \right)  + \kappa u_{yy} \right)  u_{x yy} dxdy ,\nonumber \\
		& 	J_4:= -\int_{\mathbb{R}^2}p \ast  \partial_x^{-1} u_{yyyy} u_{yy} + p \ast  u_{yyyy} u_{xyy} dxdy . \nonumber 
	\end{align}
	Applying integration by parts yields, we have $J_3=J_4=0$, and 
	\begin{align}
		\left|  I_1 \right| \leq \frac{\left| \sigma\right| }{2} \Vert u_x\Vert_{L^{\infty}} \Vert u_{yy} \Vert_{L^2}^2 + 3 \left| \sigma\right| \Vert u_{xy} \Vert_{L^2} \Vert u_{yy} \Vert_{L^2} , \label{Rod-KP LWP case s=2 eq-16} 
	\end{align}
	\begin{align}
		\left| I_2 \right| \leq & |\sigma| ( \frac{1}{2} \Vert u_x \Vert_{L^{\infty}} \Vert u_{xyy} \Vert_{L^2}+\Vert u_{xy}^2 \Vert_{L^2} + 2 \Vert u_y \Vert_{L^{\infty}} \Vert u_{xyy} \Vert_{L^2} \nonumber \\
		& +\Vert u_{xx} u_{xy} \Vert_{L^2}  + \Vert u_x \Vert_{L^{\infty}} \Vert u_{xxy} \Vert_{L^2}  +  |3-\sigma| \Vert u_y \Vert_{L^{\infty}} \Vert u_y \Vert_{L^2} )\Vert u_{xyy} \Vert_{L^2} \nonumber \\
		& + \frac{|3-\sigma|}{2} \Vert u_x\Vert_{L^{\infty}} \Vert u_{yy} \Vert_{L^2}^2. \label{Rod-KP LWP case s=2 eq-17}
	\end{align}
	Applying Lemma \ref{Sobolev embedding remark-2}, we have 
	\begin{align}
		\Vert u_{xy}^2 \Vert_{L^2} \leq C  \Vert u_{xy} \Vert_{L^2} \Vert u_{xxy} \Vert_{L^2}^{\frac{1}{2}} \Vert u_{xyy} \Vert_{\frac{1}{2}}, 
		\label{Rod-KP LWP case s=2 eq-18}
	\end{align} 
	and 
	\begin{align}
		\Vert u_{xx} u_{xy} \Vert_{L^2} \leq C \Vert u_{xx} \Vert_{L^2}^{\frac{1}{2}} \Vert u_{xy} \Vert_{L^2}^{\frac{1}{2}} \Vert u_{xxx} \Vert_{L^2}^{\frac{1}{2}} \Vert u_{xyy} \Vert_{L^2}^{\frac{1}{2}}, \label{Rod-KP LWP case s=2 eq-19}
	\end{align}
	which along with 
	\begin{align}
		\left| I_2 \right|  \leq  C \left(  \Vert \nabla u \Vert_{L^{\infty}} + \Vert \nabla u_x \Vert_{L^2} \right) \left(  \Vert \nabla u \Vert_{L^2}^2  + \Vert \nabla u_{xx}\Vert_{L^2}^2  + \Vert u_{yy}\Vert_{L^2}^2 + \Vert u_{xyy} \Vert_{L^2}^2\right) . 
		\label{Rod-KP LWP case s=2 eq-20}
	\end{align}
	Therefore, it follows from (\ref{Rod-KP LWP case s=2 eq-16}) and (\ref{Rod-KP LWP case s=2 eq-20}) that 
	\begin{align}
		& \frac{d}{dt} \left(  \Vert u_{yy} \Vert_{L^2}^2 + \Vert u_{xyy} \Vert_{L^2}^2 \right) \nonumber \\
		\leq & \left(  \Vert \nabla  u \Vert_{L^{\infty}} + \Vert  \nabla u_x \Vert_{L^2} \right) \left(\Vert \nabla u \Vert_{L^2}^2  + \Vert \nabla u_x \Vert_{L^2}^2 + \Vert \nabla u_{xx}\Vert_{L^2}^2  + \Vert u_{yy}\Vert_{L^2}^2 + \Vert u_{xyy} \Vert_{L^2}^2\right). 
		\label{Rod-KP LWP case s=2 eq-21}
	\end{align}
	Combining (\ref{Rod-KP conserved quantity}),  (\ref{Rod-KP LWP case s=2 eq-12}) and  (\ref{Rod-KP LWP case s=2 eq-21}),  there appears that 
	\begin{align}
		& \frac{d}{dt} \left(  \Vert ( u, u_x ) \Vert_{H^1}^2 + \Vert( \nabla u_{xx} , u_{yy},u_{xyy} ) \Vert_{L^2}^2  + \Vert (u, \partial_x^{-1} u ) \Vert_{H^2}^2 \right) \nonumber \\
		\leq  & C \left( \Vert (u, \nabla u ) \Vert_{L^{\infty}} + \Vert (\nabla u_x, \nabla u_{xx})\Vert_{L^2}  \right) \nonumber \\
		& \quad \quad \quad \quad \times  \left(  \Vert ( u, u_x ) \Vert_{H^1}^2 + \Vert( \nabla u_{xx} , u_{yy},u_{xyy} ) \Vert_{L^2}^2 + 
	\Vert (u, \partial_x^{-1} u ) \Vert_{H^2}^2  \right) . 	\label{Rod-KP LWP case s=2 eq-22}
	\end{align}
	Thus, we have 
	\begin{align}
		\frac{d}{dt} \left(  \Vert ( u, u_x ) \Vert_{H^1}^2 + \Vert( \nabla u_{xx} , u_{yy},u_{xyy} ) \Vert_{L^2}^2+ \Vert (u, \partial_x^{-1} ) \Vert_{H^2}^2  \right) \leq C \Vert u \Vert_{X^2}^3 .
		\label{Rod-KP LWP case s=2 eq-23}
	\end{align}
	Since $\Vert u \Vert_{X^2}^2 \leq \Vert (u,u_x ) \Vert_{H^1}^2 + \Vert( \nabla u_{xx} , u_{yy},u_{xyy} ) \Vert_{L^2}^2+ \Vert (u, \partial_x^{-1} u ) \Vert_{H^s}^2  \leq 2\Vert u \Vert_{X^2}^2$, due to Gronwall's inequality, one may find from (\ref{Rod-KP LWP case s=2 eq-23}) that 
	\begin{align}
		\sup\limits_{\tau \in [0,t]} \Vert u(\tau ) \Vert_{X^s}^2 \leq  2 \Vert u_0 \Vert_{X^s}^2  e^{C \int_0^t \Vert u (\tau) \Vert_{X^s} d\tau} , \quad \text{with $s\geq 2$}.
	\end{align}
	Taking $T>0$ such that $4C\Vert u_0 \Vert_{X^s} T <1$, we get
	\begin{align}
		\sup\limits_{\tau \in [0,t]}  \Vert u (\tau)\Vert_{X^s} \leq 2\Vert u_0 \Vert_{X^s}, \quad \forall \, t \in (0,T].\label{LWP eq-22}
	\end{align}
	It then follows that for some positive constant $C$ 
	\begin{align}
		\sup\limits_{\tau \in [0,t]} \Vert u_t(\tau) \Vert_{X^{s-2}} \leq C\Vert u_0 \Vert_{X^s}, \quad \forall \, t \in (0,T].\label{LWP eq-23}
	\end{align}
With these a priori estimates we may use the classical Firedrichs regularization mathod to constuct the approximate solutions to (\ref{the nonlocal weak form}), and then apply the compactness argument to get the local existence of a solution $u$ in $C([0,T ] ; X^s (\mathbb{R}^2)) \cap C^1 ([0,T] ; X^{s-2} (\mathbb{R}^2) )$ to equation (\ref{the nonlocal weak form}) with estimate (\ref{LWP eq-22}) and (\ref{LWP eq-23}). 

	Our attention is now turned to uniqueness of the solution. Assume that $u^{(1)}$ and $u^{(2)}$ are two solutions to (\ref{Rod-KP equation form-1}) with the same initial data, we denote that $\delta u = u^{(1)} - u^{(2)}$.  From (\ref{Rod-KP equation form-1}), we have
	\begin{equation}
		\left\{ 
		\begin{aligned}
			& \delta u_t - \delta u_{txx} +R_1 + R_2 +R_3 =0,\nonumber \\
			& \delta u|_{t=0} =0,
		\end{aligned}
		\right.\label{Rod-KP LWP eq-15}
	\end{equation}
	where $R_1 = \frac{3}{2} \left(  \left(  u^{(1)} + u^{(2)} \right) \delta u  \right)_x $, $R_2 =\sigma \left( \frac{1}{2} \left(  u^{(1)}_x + u^{(2)}_x \right)\delta u_x + u^{(1)} \delta u_{xx} + \delta u  u_{xx}^{(2)}\right)_x$, $R_3 =\delta u_x +  \partial_x^{-1} \delta u_{yy}$. 
	We first take the $L^2$-inner product between the first equation of (\ref{Rod-KP LWP eq-15}) and $\delta u$ to get 
	\begin{align}
		\frac{d}{dt} \frac{1}{2} \left(  \Vert \delta u \Vert_{L^2}^2 + \Vert \delta u_x \Vert_{L^2}^2 \right) = -\int_{\mathbb{R}^2} R_1 \delta u dxdy 
		- \int_{\mathbb{R}^2} R_2 \delta u dxdy -\int_{\mathbb{R}^2} R_3 \delta u dxdy. \label{Rod-KP LWP eq-16}
	\end{align}
	Integrating by parts yields 
$\int_{\mathbb{R}^2} R_3 \delta u dxdy  =0$. 
	Thanks to the H\"{o}lder inequality, we have 
	\begin{align}
		\left| -\int_{\mathbb{R}^2} R_1 \delta u dxdy\right|  = & \left|  \frac{3}{2} \int_{\mathbb{R}^2} (u^{(1)} + u^{(2)} ) \delta u \delta u_x dxdy \right|  \nonumber \\
		\leq  & \frac{3}{2} \Vert (u^{(1)} , u^{(2)}) \Vert_{L^{\infty}} \Vert  \delta u \Vert_{L^2} \Vert \delta u_x \Vert_{L^2},\label{Rod-KP LWP eq-18}
	\end{align}
	\begin{align}
		\left|  \int_{\mathbb{R}^2} \left(  u_x^{(1)} + u_x^{(2)} \right)  \delta u_x \delta u_x \right| \leq \Vert (u^{(1)}_x , u^{(2)}_x) \Vert_{L^{\infty}} \Vert \delta u_x \Vert_{L^2}^2, 
		\label{Rod-KP LWP eq-19}
	\end{align}
	and 
	\begin{align}
		& 	\left| \int_{\mathbb{R}^2} u^{(1)} \delta u_{xx} \delta u_x  + \delta u \delta u_x u^{(2)}_{xx} dxdy  \right|  \nonumber \\
		\leq &\left(  \frac{1}{2} \Vert u_x^{(1)} \Vert_{L^{\infty}} \Vert \delta u_x \Vert_{L^2} + \Vert  u_{xx}^{(2)}  \delta u \Vert_{L^2}\right)  \Vert  \delta u_x \Vert_{L^2} \nonumber \\
		\leq & C  \left(  \Vert u_x^{(1)} \Vert_{L^{\infty}} \Vert \delta u_x \Vert_{L^2} + \Vert  u_{xx}^{(2)}   \Vert_{L^2}^{\frac{1}{2}} \Vert u_{xxy}^{(2)} \Vert_{L^2}^{\frac{1}{2}} \Vert \delta u \Vert_{L^2}^{\frac{1}{2}} \Vert \delta u_x \Vert_{L^2}^{\frac{1}{2}}\right)  \Vert  \delta u_x \Vert_{L^2} .
		\label{Rod-KP LWP eq-20}
	\end{align}
	From (\ref{Rod-KP LWP eq-19}) and (\ref{Rod-KP LWP eq-20}), we have 
	\begin{align}
		\left|  \int_{\mathbb{R}^2} R_2 \delta u dxdy \right| \leq \frac{|\sigma|}{2} \left( 2 \Vert (u^{(1)}, u^{(2)}) \Vert_{L^{\infty}}  + \Vert u_{xx}^{(2)} \Vert_{L^2} + \Vert u_{xxy} \Vert_{L^2} \right)  \left( \Vert \delta u \Vert_{L^2}^2 + \Vert \delta u_x\Vert_{L^2}^2  \right) . \label{Rod-KP LWP eq-21}
	\end{align}
	Plugging  (\ref{Rod-KP LWP eq-18}) and (\ref{Rod-KP LWP eq-21}) into (\ref{Rod-KP LWP eq-16}) gives rise to 
	\begin{align}
		\frac{d}{dt} \left(  \Vert \delta u \Vert_{L^2}^2 + \Vert \delta u_x \Vert_{L^2}^2 \right) \leq C \Vert (u^{(1)}, u^{(2)} ) \Vert_{X^2} \left(  \Vert \delta u \Vert_{L^2}^2 + \Vert \delta u_x \Vert_{L^2}^2 \right) . \label{Rod-KP LWP eq-22}
	\end{align}
	Using Gronwall's inequality, we deduce that $\delta u=0$, that is $ u^{(1)} = u^{(2)}$.  This ends the proof of the uniqueness. Moreover, thanks to (\ref{LWP eq-22}) and Fatou's lemma, we may verify that the solution $u$ depends continuously on the initial value $u_0 \in X^s $, and this complete the proof of Theorem   \ref{LWP theorem}. 
\end{proof}
\subsection{Blow-up criterion}
We now turn to the proof of the blow-up criterion  for  (\ref{the nonlocal weak form}). The argument follows the same spirit as that of Theorem \ref{LWP theorem}, but additional techniques require a more detailed presentation. 
\begin{proof}[Proof of Theorem \ref{Blow-up criterion}]
 We will prove the theorem by considering different cases of  the regular index $s$ $(s\geq 2)$.	To this end, we divide the proof into four steps.\\
	{\bf Step 1: $\mathbf{s=2}$.} From (\ref{Rod-KP conserved quantity}) and  (\ref{Rod-KP LWP case s=2 eq-11}), we have 
	\begin{align}
		\frac{d}{dt} \left( \Vert \nabla u(\tau) \Vert_{L^2}^2 + \Vert \nabla u_x (\tau)\Vert_{L^2}^2  \right)   \leq & \left( 4|\sigma| + |3-\sigma|\right)  \Vert \nabla u \Vert_{L^{\infty}} \nonumber \\
		& \times  \left( \sqrt{2E(u_0)} + \Vert \nabla u(\tau) \Vert_{L^2}^2 + \Vert \nabla u_x (\tau)\Vert_{L^2}^2 \right) .
		\label{Rod-KP blow-up criterion eq-1}
	\end{align}
	Applying Gronwall's inequality to (\ref{Rod-KP blow-up criterion eq-1}) yields that  for all $t \in [0,T^{\ast}_{u_0})$,
	\begin{align}
		& \sup\limits_{\tau \in[0,t]} \left( \Vert \nabla u(\tau) \Vert_{L^2}^2 + \Vert \nabla u_x (\tau)\Vert_{L^2}^2  \right)  \nonumber \\
		\leq & \left( \sqrt{2E(u_0)}+  \Vert \nabla u_0 \Vert_{L^2}^2 + \Vert \nabla u_{0,x} \Vert_{L^2}^2 \right) e^{\left( 4|\sigma| + |3-\sigma|\right) \int_0^t  \Vert \nabla u(\tau) \Vert_{L^{\infty}} d \tau } -\sqrt{2E(u_0)}.\label{Rod-KP blow-up criterion eq-2}
	\end{align}
	Therefore, if $\int_0^{T_{u_0}^{\ast}} \Vert \nabla u(\tau ) \Vert_{L^{\infty}} d\tau$ is finite, then for all $t \in [0,T_{u_0}^{\ast})$, 
	\begin{align}
		\left( \Vert \nabla u(\tau) \Vert_{L^2}^2 + \Vert \nabla u_x (\tau)\Vert_{L^2}^2  \right)  \leq M_1 (T_{u_0}^{\ast}), \label{Rod-KP blow-up criterion eq-3}
	\end{align}
	where
	\begin{align}
		M_1 (T_{u_0}^{\ast}) :=  \left( \sqrt{2E(u_0)}+  \Vert \nabla u_0  \Vert_{L^2}^2 + \Vert \nabla u_{0,x} \Vert_{L^2}^2 \right) e^{\left( 4|\sigma| + |3-\sigma|\right) \int_0^t  \Vert \nabla u(\tau) \Vert_{L^{\infty}} d \tau } -\sqrt{2E(u_0)}.\nonumber 
	\end{align}
	It is inferred from (\ref{Rod-KP LWP case s=2 eq-12}) and (\ref{Rod-KP blow-up criterion eq-3})  that 
	\begin{align}
		& \frac{d}{dt} \left( \Vert \nabla u(\tau) \Vert_{L^2}^2 + \Vert \nabla u_x (\tau)\Vert_{L^2}^2  + \Vert \nabla u_{xx} (\tau)  \Vert_{L^2}^2 \right)  \nonumber \\
		\leq  & C \left(  \Vert \nabla u \Vert_{L^{\infty}} + \sqrt{M_1(T_{u_0}^{\ast})} \right) \left( \sqrt{2E(u_0)} + \Vert \nabla u(\tau) \Vert_{L^2}^2 + \Vert \nabla u_x (\tau)\Vert_{L^2}^2  + \Vert \nabla u_{xx} (\tau)  \Vert_{L^2}^2 \right) .
		\label{Rod-KP blow-up criterion eq-5}
	\end{align}
	Applying Gronwall's inequality to (\ref{Rod-KP blow-up criterion eq-5}), we get 
	\begin{align}
		& \sup\limits_{\tau \in [0,t]} \left( 	\Vert \nabla u(\tau) \Vert_{L^2}^2 + \Vert \nabla u_x (\tau)\Vert_{L^2}^2  + \Vert \nabla u_{xx} (\tau)  \Vert_{L^2}^2 \right)  \nonumber \\
		\leq & \left(	\sqrt{E(u_0)} + \Vert \nabla u_0 \Vert_{L^2}^2 + \Vert \nabla u_{0,x}\Vert_{L^2}^2  + \Vert \nabla u_{0,xx}  \Vert_{L^2}^2 \right) e^{C \int_0^t \Vert \nabla u \Vert_{L^{\infty}} d\tau + \sqrt{M_1 (T^{\ast})} t } -\sqrt{2E(u_0)}.	\label{Rod-KP blow-up criterion eq-6}
	\end{align}
	Hence, we get, for all $t \in [0,T_{u_0}^{\ast})$,
	\begin{align}
		\Vert \nabla u(\tau) \Vert_{L^2}^2 + \Vert \nabla u_x (\tau)\Vert_{L^2}^2  + \Vert \nabla u_{xx} (\tau)  \Vert_{L^2}^2  \leq M_2 (T^{\ast}) ,	\label{Rod-KP blow-up criterion eq-7}
	\end{align}
	where 
	\begin{align}
		M_2(T^{\ast}) := \left(	\sqrt{E(u_0)} + \Vert \nabla u_0 \Vert_{L^2}^2 + \Vert \nabla u_{0,x}\Vert_{L^2}^2  + \Vert \nabla u_{0,xx}  \Vert_{L^2}^2 \right) e^{C \int_0^{T^{\ast}} \Vert \nabla u \Vert_{L^{\infty}} d\tau + \sqrt{M_1 (T^{\ast})} T^{\ast} } -\sqrt{E(u_0)}.\nonumber 
	\end{align}
	Similarly, thanks to (\ref{Rod-KP LWP case s=2 eq-22}), we have 
	\begin{align}
		& \frac{d}{dt} \left( \Vert (u,u_x ) \Vert_{H^1}^2 + \Vert (\nabla u_{xx}, u_{yy}, u_{xyy} ) \Vert_{L^2}^2 + \Vert (u,\partial_x^{-1} u) \Vert_{H^2}^2 \right)  \nonumber \\
		\leq & C \left( \sqrt{E(u_0)}  + \sqrt{M_2(T^{\ast})}+ \Vert \nabla  u \Vert_{L^{\infty}}\right)  \nonumber \\
		& \quad \quad \quad \quad  \times  \left( \Vert (u,u_x ) \Vert_{H^1}^2 + \Vert (\nabla u_{xx}, u_{yy}, u_{xyy} ) \Vert_{L^2}^2 + \Vert (u,\partial_x^{-1} u) \Vert_{H^2}^2 \right) .
			\label{Rod-KP blow-up criterion eq-8}
	\end{align}
	By Gronwall's inequality, we have for all $t\in[0,T^{\ast})$
	\begin{align}
		& \sup\limits_{\tau\in [0,t]} \left( \Vert (u,u_x ) \Vert_{H^1}^2 + \Vert (\nabla u_{xx}, u_{yy}, u_{xyy} ) \Vert_{L^2}^2 + \Vert(u, \partial_x^{-1} u)\Vert_{H^2}^2  \right)  \nonumber \\
		\leq & \left( \Vert (u_0 ,u_{0,x} ) \Vert_{H^1}^2 + \Vert (\nabla u_{0,xx}, u_{0,yy}, u_{0,xyy} ) \Vert_{L^2}^2+ \Vert(u_0, \partial_x^{-1} u_0)\Vert_{H^2}^2  \right)  \nonumber \\
		& \quad \quad \quad \times e^{ C \left( \int_0^t  \Vert\nabla u (\tau) \Vert_{L^{\infty}} d\tau + (\sqrt{E(u_0)}+ \sqrt{M_2(T^{\ast})}) t\right) }. \nonumber 
	\end{align}
	Hence, we get 
	\begin{align}
		\Vert (u,u_x ) \Vert_{H^1}^2 + \Vert (\nabla u_{xx}, u_{yy}, u_{xyy} ) \Vert_{L^2}^2 + \Vert (u, \partial_x^{-1} u) \Vert_{H^s}^2  \leq M_3(T^{\ast}), \nonumber
	\end{align}
	where
	\begin{align}
		M_3(T^{\ast}):=  & \left( \Vert (u_0 ,u_{0,x} ) \Vert_{H^1}^2 + \Vert (\nabla u_{0,xx}, u_{0,yy}, u_{0,xyy} ) \Vert_{L^2}^2  + \Vert (u_0, \partial_x^{-1} u_0) \Vert_{H^2}^2 \right)  \nonumber \\
		& \quad \quad \quad \times e^{ C \left( \int_0^{T^{\ast}}  \Vert\nabla u (\tau) \Vert_{L^{\infty}} d\tau + (\sqrt{E(u_0)}+ \sqrt{M_2(T^{\ast})}) T^{\ast}  \right) }. \nonumber 
	\end{align}
	It is then follows that 
	\begin{align}
		\Vert u( t ) \Vert_{X^2}^2 \leq2 M_3 (T^{\ast}),\quad  \forall \, t \in [0,T^{\ast}).  \label{Rod-KP blow-up criterion eq-89}
	\end{align}
This is a contradiction, which completes the proof of Theorem \ref{Blow-up criterion} for $s=2$. \\
{\bf Step 2: $\mathbf{s=3}$.} Differentiating (\ref{Rod-KP LWP case s=2 eq-3}) with respect to $x$, we have 
	\begin{align}
	& \partial_t \nabla u_{xxx} +\sigma \left( u\nabla u_{xxxx} + 3u_x\nabla u_{xxx} +4u_{xx} \nabla u_{xx} +3u_{xxx} \nabla u_x  + u_{xxxx} \nabla u\right) 
	\nonumber \\
	&\quad \quad \quad \quad -(3-\sigma) \left( u\nabla u+ u_{xx} \nabla u + u_x \nabla u_x + u \nabla u_{xx} \right) -\kappa (\nabla u_{xx}  +\nabla u_x )- 3u_x \nabla u_x \nonumber \\
	& \quad \quad \quad \quad +p \ast \left(\left( 3-\sigma \right) u\nabla u + \sigma u_x \nabla u_x+\kappa \nabla u_x  \right) +p \ast \nabla u_{xxyy} =0.\label{Rod-KP blow-up criterion eq-9}
	\end{align} 
Multiplying (\ref{Rod-KP blow-up criterion eq-9}) by $\nabla u_{xxx}$ and integrating  it on $\mathbb{R}^2$, we get 
\begin{align}
\frac{1}{2} 	\frac{d}{dt}  \Vert \nabla u_{xxx} \Vert_{L^2}^2 = \sum\limits_{i=1}^4 K_i ,\nonumber 
	\end{align} 
where 
\begin{align}
	& K_1= - \sigma \int_{\mathbb{R}^2} \left(  u\nabla u_{xxxx} + 3u_x\nabla u_{xxx} +4u_{xx} \nabla u_{xx} +3u_{xxx} \nabla u_x  + u_{xxxx} \nabla u \right) \nabla u_{xxx} dxdy ,\nonumber \\
	& K_2 = \int_{\mathbb{R}^2}  (3-\sigma) \left( u\nabla u+ u_{xx} \nabla u + u_x \nabla u_x + u \nabla u_{xx}  \right) \nabla u_{xxx}  + \left(  \kappa (\nabla u_{xx}  +\nabla u_x )  + 3u_x \nabla u_x\right) \nabla u_{xxx} dxdy, \nonumber \\
	& K_3 = -\int_{\mathbb{R}^2}  p \ast \left(\left( 3-\sigma \right) u\nabla u + \sigma u_x \nabla u_x+\kappa \nabla u_x  \right) 
	\nabla u_{xxx} dxdy, \nonumber \\
	& K_4 = -\int_{\mathbb{R}^2} p \ast \nabla u_{xxyy} \nabla u_{xxx} dxdy. \nonumber 
\end{align}
Thanks to  intergation by parts and Lemma \ref{Sobolev embedding remark-2}, then we have $K_4=0$, and 
\begin{align}
	\left| K_1 \right|  \leq  & 2 |\sigma| \Vert u_x \Vert_{L^{\infty}}  \Vert \nabla u_{xxx} \Vert_{L^2}^2 + 4 |\sigma | \Vert  u_{xx} \nabla u_{xx} \Vert_{L^2} \Vert \nabla u_{xxx} \Vert_{L^2}   \nonumber \\
	& + 3|\sigma| \Vert u_{xxx} \nabla u_x \Vert_{L^2} \Vert \nabla u_{xxx}\Vert_{L^2}  + |\sigma| 
	\Vert \nabla  u \Vert_{L^{\infty}} \Vert \nabla u_{xxx} \Vert_{L^2}^2  \nonumber \\
	\leq & 2 |\sigma| \Vert u_x \Vert_{L^{\infty}}  \Vert \nabla u_{xxx} \Vert_{L^2}^2 + 4 |\sigma | \Vert  u_{xx} \Vert_{L^2}^{\frac{1}{2}} \Vert u_{xxy} \Vert_{L^2}^{\frac{1}{2}} \Vert \nabla u_{xx} \Vert_{L^2}^{\frac{1}{2}} \Vert \nabla u_{xxx} \Vert_{L^2}^{\frac{1}{2}} \Vert \nabla u_{xxx} \Vert_{L^2}   \nonumber \\
	& + 3|\sigma| \Vert u_{xxx}  \Vert_{L^2}^{\frac{1}{2}} \Vert u_{xxxy} \Vert_{L^2}^{\frac{1}{2}}  \Vert \nabla u_x \Vert_{L^2}^{\frac{1}{2}} \Vert \nabla u_{xx} \Vert_{L^2}^{\frac{1}{2}} \Vert \nabla u_{xxx} \Vert_{L^2} + |\sigma| 
	\Vert \nabla  u \Vert_{L^{\infty}} \Vert \nabla u_{xxx} \Vert_{L^2}^2  \nonumber \\
	\leq & C(\Vert \nabla u \Vert_{L^{\infty}} + \Vert \nabla u_{xx} \Vert_{L^2} ) \left(  \Vert \nabla u_x\Vert_{L^2}^2 + \Vert \nabla u_{xxx} \Vert_{L^2}^2 \right) ,\nonumber 
\end{align}
\begin{align}
	| K_2 |  \leq & |3-\sigma | (\Vert \nabla u \Vert_{L^{\infty}} \Vert u \Vert_{L^2} + \Vert \nabla u \Vert_{L^{\infty}} \Vert u_{xx} \Vert_{L^2} + \Vert u_x\Vert_{L^{\infty}} \Vert \nabla u_x \Vert_{L^2} )  \Vert \nabla u_{xxx} \Vert_{L^2}  \nonumber \\
	& + |3-\sigma | \Vert u_x\Vert_{L^{\infty}} \Vert \nabla u_{xx} \Vert_{L^2}^2 +|\kappa| \Vert \nabla u_{xx} \Vert_{L^2}^2 + 3 \Vert u_x\Vert_{L^{\infty}} \Vert \nabla u_x \Vert_{L^2} \Vert \nabla u_{xxx} \Vert_{L^2}  \nonumber \\
\leq & C(\Vert \nabla u \Vert_{L^{\infty}} + 1) \left(  \Vert u \Vert_{L^2}^2 + \Vert \nabla u_x \Vert_{L^2}^2 + \Vert \nabla u_{xx} \Vert_{L^2}^2 + \Vert \nabla u_{xxx} \Vert_{L^2}^2 \right) ,\nonumber 
\end{align}
\begin{align}
 |K_3|\leq & \frac{1}{2} ( |3-\sigma|\Vert \nabla u \Vert_{L^{\infty}} \Vert u \Vert_{L^2}  + |\sigma| \Vert u_x\Vert_{L^{\infty}}  \Vert \nabla u_x \Vert_{L^2} +|\kappa| \Vert \nabla u_x \Vert_{L^2} )\Vert \nabla u_{xxx} \Vert_{L^2} \nonumber \\
 \leq & C ( \Vert \nabla u \Vert_{L^{\infty}} + 1) \left(  \Vert u \Vert_{L^2}^2 + \Vert \nabla u_x \Vert_{L^2}^2 + \Vert \nabla u_{xxx} \Vert_{L^2}^2 \right) . \nonumber 
\end{align}
Therefore, we obtain
\begin{align}
	\frac{d}{dt} \Vert \nabla u_{xxx} \Vert_{L^2}^2 \leq  & C  \left( \Vert \nabla u \Vert_{L^{\infty}} + \Vert \nabla u_x \Vert_{L^2} +1  \right) \nonumber \\
	& \times \left(  \Vert  u \Vert_{L^2}^2 + \Vert \nabla u_x \Vert_{L^2}^2 + \Vert \nabla u_{xx} \Vert_{L^2}^2 + \Vert \nabla u_{xxx} \Vert_{L^2}^2  \right) , \nonumber 
\end{align}
which follows 
\begin{align}
	&	\frac{d}{dt} \left( \Vert u\Vert_{L^2}^2 + \Vert \nabla u \Vert_{L^2}^2 + \Vert \nabla u_x \Vert_{L^2}^2 + \Vert \nabla u_{xx} \Vert_{L^2}^2 + \Vert \nabla u_{xxx} \Vert_{L^2}^2 \right)\nonumber \\
	\leq & C \left( \Vert \nabla u \Vert_{L^{\infty}} + \Vert \nabla u_x \Vert_{L^2} + 1 \right)  \nonumber \\
	&	\times \left(\Vert u \Vert_{L^2}^2 + \Vert \nabla u \Vert_{L^2}^2 + \Vert \nabla u_x \Vert_{L^2}^2 + \Vert \nabla u_{xx} \Vert_{L^2}^2 + \Vert \nabla u_{xxx} \Vert_{L^2}^2\right) . \label{Rod-KP blow-up criterion eq-12}
\end{align}
Applying Gronwall's inequality to (\ref{Rod-KP blow-up criterion eq-12}) yields that for all $t \in [0,T^{\ast})$
\begin{align}
	& \sup\limits_{\tau \in [0,t] } \left( \Vert u \Vert_{L^2}^2 + \Vert \nabla u \Vert_{L^2}^2 + \Vert \nabla u_x \Vert_{L^2}^2 + \Vert \nabla u_{xx} \Vert_{L^2}^2 + \Vert \nabla u_{xxx} \Vert_{L^2}^2\right) \nonumber \\
	\leq &\left( \Vert u_0 \Vert_{L^2}^2 + \Vert \nabla u_0 \Vert_{L^2}^2 + \Vert \nabla u_{0,x} \Vert_{L^2}^2 + \Vert \nabla u_{0,xx} \Vert_{L^2}^2 + \Vert \nabla u_{0,xxx} \Vert_{L^2}^2\right) 
	\times e^{C\left( \int_0^t \Vert \nabla u \Vert_{L^{\infty}} d\tau + t\right) }. \nonumber 
\end{align}
Hence, we get, for all $t \in [0,T^{\ast})$, we have 
\begin{align}
	\sup\limits_{\tau \in [0,t] } \left( \Vert u \Vert_{L^2}^2 + \Vert \nabla u \Vert_{L^2}^2 + \Vert \nabla u_x \Vert_{L^2}^2 + \Vert \nabla u_{xx} \Vert_{L^2}^2 + \Vert \nabla u_{xxx} \Vert_{L^2}^2\right) \leq  M_4(T^{\ast}),\nonumber 
\end{align}
where 
\begin{align}
	M_4(T^{\ast}) := \left( \Vert u_0 \Vert_{L^2}^2 + \Vert \nabla u_0 \Vert_{L^2}^2 + \Vert \nabla u_{0,x} \Vert_{L^2}^2 + \Vert \nabla u_{0,xx} \Vert_{L^2}^2 + \Vert \nabla u_{0,xxx} \Vert_{L^2}^2\right) 
	e^{C\left( \int_0^{T^{\ast}} \Vert \nabla u \Vert_{L^{\infty}} d\tau + T^{\ast} \right)} . \nonumber 
\end{align}

Then, acting the operators $\nabla $ to (\ref{Rod-KP LWP case s=2 eq-13}) and (\ref{Rod-KP LWP case s=2 eq-14}) respectively gives rise to 
\begin{align}
	&\partial_t \nabla u_{yy} + \sigma \left( \nabla u u_{xyy}+ u\nabla u_{xyy} + 2\nabla u_y u_{xy} + 2u_y \nabla u_{xy} + \nabla u_x u_{xy} +u_x \nabla u_{xy}  \right) \nonumber \\
	& \quad\quad + p_x \ast \left((3-\sigma ) \left( \nabla u u_{yy} + u \nabla u_{yy} + 2u_y \nabla u_y \right) +\sigma \left( \nabla u_x u_{xyy} + u_x \nabla u_{xyy} + 2u_{xy} \nabla u_{xy} \right) + \kappa \nabla u_{yy}   \right) \nonumber \\
	& \quad \quad + p \ast \partial_x^{-1} \nabla u_{yyyy} =0,\label{Rod-KP blow-up criterion eq-10}\\
	& \partial_t \nabla u_{xyy} + \sigma ( \nabla u u_{xxyy} + u \nabla u_{xxyy} + 2\nabla u_{xy} u_{xy} + 2\nabla u_y u_{xy} + 2 u_y \nabla u_{xy} + \nabla u_{xx} u_{xy} \nonumber \\
	& \quad \quad + u_{xx} \nabla u_{xy} + \nabla u_x u_{xxy}+ u_x\nabla u_{xxy}  ) -(3-\sigma) (\nabla u u_{yy} + u \nabla u_{yy} + 2\nabla u_y u_y) \nonumber \\
	& \quad\quad +p\ast ((3-\sigma) (\nabla u u_{yy} + u \nabla u_{yy} +2\nabla u_y u_y ) +\sigma \left(\nabla u_x u_{xyy} + u_x \nabla u_{xyy} + 2\nabla u_{xy} u_{xy} +\kappa \nabla u_{yy}  \right) ) \nonumber \\
	& \quad\quad + p \ast \nabla u_{yyyy} =0. \label{Rod-KP blow-up criterion eq-11}
\end{align}
Taking $L^2$ inner-products between  two equations in (\ref{Rod-KP blow-up criterion eq-10})-(\ref{Rod-KP blow-up criterion eq-11}) and $\nabla u_{yy}$, $\nabla u_{xyy}$, and summating them, we get 
\begin{align}
	\frac{1}{2} \frac{d}{dt} \left( \Vert \nabla u_{yy} \Vert_{L^2}^2 + \Vert \nabla u_{xyy} \Vert_{L^2}^2 \right)  = \sum\limits_{i=1}^4 L_i, \nonumber 
\end{align}
where
\begin{align}
& L_1 = -\int_{\mathbb{R}^2} \sigma \left( \nabla u u_{xyy}+ u\nabla u_{xyy} + 2\nabla u_y u_{xy} + 2u_y \nabla u_{xy} + \nabla u_x u_{xy} +u_x \nabla u_{xy}  \right)  \nabla u_{yy} dxdy, \nonumber \\
& L_2 =-\int_{\mathbb{R}^2} \sigma ( \nabla u u_{xxyy} + u \nabla u_{xxyy} + 2\nabla u_{xy} u_{xy} + 2\nabla u_y u_{xy} + 2 u_y \nabla u_{xy}  + \nabla u_{xx} u_{xy} + u_{xx} \nabla u_{xy} \nonumber \\
&   \quad \quad\quad\quad + \nabla u_x u_{xxy}+ u_x\nabla u_{xxy}  ) \nabla u_{xyy} -(3-\sigma) (\nabla u u_{yy} + u \nabla u_{yy} + 2\nabla u_y u_y)  \nabla u_{xyy} dxdy ,\nonumber \\
& L_3= - \int_{\mathbb{R}^2} p_x\ast ((3-\sigma) (\nabla u u_{yy} + u \nabla u_{yy} +2\nabla u_y u_y ) \nonumber \\
&\quad \quad\quad\quad +\sigma \left(\nabla u_x u_{xyy} + u_x \nabla u_{xyy} + 2\nabla u_{xy} u_{xy} +\kappa \nabla u_{yy}  \right) )  \nabla u_{yy} dxdy \nonumber \\
& \quad \quad\quad\quad -\int_{\mathbb{R}^2} p\ast ((3-\sigma) (\nabla u u_{yy} + u \nabla u_{yy} +2\nabla u_y u_y ) \nonumber \\
& \quad \quad\quad\quad +\sigma \left(\nabla u_x u_{xyy} + u_x \nabla u_{xyy} + 2\nabla u_{xy} u_{xy} +\kappa \nabla u_{yy}  \right) )  \nabla u_{xyy} dxdy,\nonumber \\
&L_4=-\int_{\mathbb{R}^2} p \ast \partial_x^{-1} \nabla u_{yyyy}  \nabla u_{yy} dxdy -\int_{\mathbb{R}^2} p \ast \nabla u_{yyyy} \nabla u_{xyy}dxdy. \nonumber 
\end{align}
Applying integration by parts yields $L_3 =L_4=0$, and applying Lemma \ref{Sobolev embedding remark-2} again, we have 
\begin{align}
	|L_1| \leq  &|\sigma| (\Vert \nabla u\Vert_{L^{\infty}} \Vert u_{xyy} \Vert_{L^2}  + \Vert u_x\Vert_{L^{\infty}} \Vert \nabla u_{yy}\Vert_{L^2}  +2 \Vert \nabla u_y  u_{xy} \Vert_{L^2} + 2\Vert u_y \Vert_{L^{\infty}}\Vert \nabla u_{xy} \Vert_{L^2} \nonumber \\
	& + \Vert \nabla u_x  u_{xy} \Vert_{L^2} + \Vert u_x \Vert_{L^{\infty}} \Vert \nabla u_{xy} \Vert_{L^2} ) \Vert \nabla u_{yy} \Vert_{L^2} \nonumber \\ 
	\leq  &|\sigma| (\Vert \nabla u\Vert_{L^{\infty}} \Vert u_{xyy} \Vert_{L^2}  + \Vert u_x\Vert_{L^{\infty}} \Vert \nabla u_{yy}\Vert_{L^2}  +2 \Vert \nabla u_y\Vert_{L^2}^{\frac{1}{2}} \Vert \nabla u_{yy}\Vert_{L^2}^{\frac{1}{2}}  \Vert u_{xy} \Vert_{L^2}^{\frac{1}{2}} \Vert u_{xxy}\Vert_{L^2}^{\frac{1}{2}} \nonumber \\
	& + 2\Vert u_y \Vert_{L^{\infty}}\Vert \nabla u_{xy} \Vert_{L^2}+ \Vert \nabla u_x \Vert_{L^2}^{\frac{1}{2}} \Vert \nabla u_{xx} \Vert_{L^2}^{\frac{1}{2}} \Vert u_{xy} \Vert_{L^2}^{\frac{1}{2}} \Vert u_{xyy} \Vert_{L^2}^{\frac{1}{2}}  + \Vert u_x \Vert_{L^{\infty}} \Vert \nabla u_{xy} \Vert_{L^2} ) \Vert \nabla u_{yy} \Vert_{L^2} \nonumber \\ 
	\leq & C (\Vert \nabla u \Vert_{L^{\infty}} + \Vert u_{yy} \Vert_{L^2} + \Vert \nabla u_x \Vert_{L^2}) \left( \Vert u_{xyy} \Vert_{L^2}^2 + \Vert \nabla u_{xx} \Vert_{L^2}^2 +\Vert \nabla u_{yy} \Vert_{L^2}^2 \right) , \label{Rod-KP blow-up criterion eq-11.1}
\end{align}
and
\begin{align}
	|L_2| \leq & |\sigma | ( \Vert \nabla u \Vert_{L^{\infty}} \Vert \nabla u_{xyy} \Vert_{L^2} + \Vert u_x \Vert_{L^{\infty}}\Vert\nabla u_{xyy} \Vert_{L^2} + 2 \Vert \nabla u_{xy} u_{xy} \Vert_{L^2} + 2 \Vert \nabla u_y u_{xy} \Vert_{L^2} \nonumber \\
	&  + 2\Vert u_y \Vert_{L^{\infty}} \Vert \nabla u_{xy} \Vert_{L^2} + \Vert \nabla u_{xx} u_{xy}\Vert_{L^2} + \Vert u_{xx} \nabla u_{xy} \Vert_{L^2}+ \Vert \nabla u_x u_{xxy} \Vert_{L^2} \nonumber \\
	&  + \Vert u_x \Vert_{L^{\infty}} \Vert \nabla u_{xxy} \Vert_{L^2}   ) \Vert \nabla u_{xyy} \Vert_{L^2}  + |3-\sigma|( \Vert \nabla u \Vert_{L^{\infty}} \Vert u_{yy} \Vert_{L^2} +2\Vert u_{y} \Vert_{L^{\infty}}\Vert \nabla u_y \Vert_{L^2} )  \Vert \nabla u_{xyy} \Vert_{L^2} \nonumber \\
	&+ \Vert u_x \Vert_{L^{\infty}} \Vert \nabla u_{yy} \Vert_{L^2}^2 \nonumber \\
	\leq & |\sigma | ( \Vert \nabla u \Vert_{L^{\infty}} \Vert \nabla u_{xyy} \Vert_{L^2} + \Vert u_x \Vert_{L^{\infty}}\Vert\nabla u_{xyy} \Vert_{L^2} + 2 \Vert \nabla u_{xy}\Vert_{L^2}^{\frac{1}{2}} \Vert \nabla u_{xyy} \Vert_{L^2}^{\frac{1}{2}} \Vert u_{xy} \Vert_{L^2}^{\frac{1}{2}} \Vert u_{xxy} \Vert_{L^2}^{\frac{1}{2}} \nonumber \\
	& + 2 \Vert \nabla u_y\Vert_{L^2}^{\frac{1}{2}} \Vert \nabla u_{yy} \Vert_{L^2}^{\frac{1}{2}} \Vert  u_{xy} \Vert_{L^2}^{\frac{1}{2}} \Vert u_{xxy}\Vert^{\frac{1}{2}}  + 2\Vert u_y \Vert_{L^{\infty}} \Vert \nabla u_{xy} \Vert_{L^2} + \Vert \nabla u_{xx}\Vert_{L^2}^{\frac{1}{2}} \Vert \nabla u_{xxy} \Vert_{L^2}^{\frac{1}{2}} \Vert u_{xy}\Vert_{L^2}^{\frac{1}{2}} \Vert u_{xxy} \Vert_{L^2}^{\frac{1}{2}}  \nonumber \\
	&+ \Vert u_{xx}\Vert_{L^2}^{\frac{1}{2}} \Vert u_{xxy} \Vert_{L^2}^{\frac{1}{2}} \Vert \nabla u_{xy} \Vert_{L^2}^{\frac{1}{2}} \Vert \nabla u_{xxy} \Vert_{L^2}^{\frac{1}{2}}  + \Vert \nabla u_x\Vert_{L^2}^{\frac{1}{2}} \Vert \nabla u_{xy}\Vert_{L^2}^{\frac{1}{2}} \Vert u_{xxy} \Vert_{L^2}^{\frac{1}{2}} \Vert u_{xxxy} \Vert_{L^2}^{\frac{1}{2}}
	\nonumber \\
	&+ \Vert u_x \Vert_{L^{\infty}} \Vert \nabla u_{xxy} \Vert_{L^2}   ) \Vert \nabla u_{xyy} \Vert_{L^2}+ |3-\sigma|( \Vert \nabla u \Vert_{L^{\infty}} \Vert u_{yy} \Vert_{L^2}2\Vert u_{y} \Vert_{L^{\infty}}\Vert \nabla u_y \Vert_{L^2} )  \Vert \nabla u_{xyy} \Vert_{L^2}  
	\nonumber \\
	& + \Vert u_x\Vert_{L^{\infty}} \Vert \nabla u_{yy} \Vert_{L^2}^2  \nonumber \\
	\leq & C \left( \Vert \nabla u  \Vert_{L^{\infty}} +\Vert \nabla u_x \Vert_{L^2}+ \Vert \nabla u_{xx}\Vert_{L^2} +\Vert u_{yy} \Vert_{L^2} + \Vert u_{xyy} \Vert_{L^2}\right) \nonumber \\
	& \times \left( \Vert \nabla u_{xx} \Vert_{L^2}^2 + \Vert  \nabla u_{xyy} \Vert_{L^2}^2 + \Vert \nabla u_{xxx} \vert_{L^2}^2 + \Vert \nabla u_{xyy} \Vert_{L^2}^2  \right) . \label{Rod-KP blow-up criterion eq-11.2}
\end{align}
It then follows form (\ref{Rod-KP blow-up criterion eq-11.1}) and (\ref{Rod-KP blow-up criterion eq-11.2}) that 
\begin{align}
& 	\frac{d}{dt} \left(  \Vert \nabla u_{yy} \Vert_{L^2}^2 + \Vert \nabla u_{xyy} \Vert_{L^2}^2\right)\nonumber \\
	\leq & C \left(\Vert \nabla u  \Vert_{L^{\infty}} +\Vert \nabla u_x \Vert_{L^2}+ \Vert \nabla u_{xx}\Vert_{L^2} +\Vert u_{yy} \Vert_{L^2} + \Vert u_{xyy} \Vert_{L^2} \right) \nonumber \\
	& \times \left( \Vert \nabla u_{xx} \Vert_{L^2}^2 +\Vert \nabla u_{yy} \Vert_{L^2}^2 +  \Vert  u_{xyy} \Vert_{L^2}^2 + \Vert \nabla u_{xxx} \vert_{L^2}^2 + \Vert \nabla u_{xyy} \Vert_{L^2}^2 \right) .	\label{Rod-KP blow-up criterion eq-13}
\end{align}
Combining (\ref{Rod-KP blow-up criterion eq-12}), (\ref{Rod-KP blow-up criterion eq-13}) with  (\ref{Rod-KP LWP eq-8}) implies 
\begin{align}
	& \frac{d}{dt} \left( \Vert (u,u_x) \Vert_{H^1}^2+ 	\Vert (\nabla u_{xx} ,\nabla u_{xxx}) \Vert_{L^2}^2 + \Vert ( u_{yy}, u_{xyy})\Vert_{L^2}^2+ \Vert (\nabla u_{yy}, \nabla u_{xyy} )   \Vert_{L^2}^2 + \Vert (u, \partial_x^{-1} u ) \Vert_{H^3}^2  \right) \nonumber \\
	\leq &  C \left(\sqrt{E(u_0)}+  \Vert \nabla u \Vert_{L^{\infty}} +\Vert \nabla u_x \Vert_{L^2}+ \Vert \nabla u_{xx}\Vert_{L^2} +\Vert u_{yy} \Vert_{L^2} + \Vert u_{xyy} \Vert_{L^2} \right) \nonumber \\
	& \times \left( \Vert (u,u_x) \Vert_{H^1}^2+ 	\Vert (\nabla u_{xx} ,\nabla u_{xxx}) \Vert_{L^2}^2 + \Vert ( u_{yy}, u_{xyy})\Vert_{L^2}^2+ \Vert (\nabla u_{yy}, \nabla u_{xyy} )   \Vert_{L^2}^2  \right) \nonumber \\
	\leq &C \left(  \Vert \nabla u \Vert_{L^{\infty}} + \sqrt{M_3(T^{\ast})}  +\sqrt{E(u_0)} \right) \nonumber \\
	& \times \left( \Vert (u,u_x) \Vert_{H^1}^2+ 	\Vert (\nabla u_{xx} ,\nabla u_{xxx}) \Vert_{L^2}^2 + \Vert ( u_{yy}, u_{xyy})\Vert_{L^2}^2+ \Vert (\nabla u_{yy}, \nabla u_{xyy} )   \Vert_{L^2}^2 + \Vert (u, \partial_x^{-1} u ) \Vert_{H^3}^2 \right) . \nonumber 
\end{align}
Thus, for all $t \in [0,T^{\ast})$, we have 
\begin{align}
	& \sup\limits_{\tau \in [0,t]} \left(  \Vert (u,u_x) \Vert_{H^1}^2+ 	\Vert (\nabla u_{xx} ,\nabla u_{xxx}) \Vert_{L^2}^2 + \Vert ( u_{yy}, u_{xyy})\Vert_{L^2}^2+ \Vert (\nabla u_{yy}, \nabla u_{xyy} )   \Vert_{L^2}^2  + \Vert (u, \partial_x^{-1} u) \Vert_{H^3}^2\right) \nonumber \\
	\leq &C  (  \Vert (u_0,u_{0,x}) \Vert_{H^1}^2+ 	\Vert (\nabla u_{0, xx} ,\nabla u_{0, xxx}) \Vert_{L^2}^2 + \Vert ( u_{0, yy}, u_{0, xyy})\Vert_{L^2}^2+ \Vert (\nabla u_{0, yy}, \nabla u_{0, xyy} )   \Vert_{L^2}^2 \nonumber \\
	& + \Vert (u_0, \partial_x^{-1} u_0) \Vert_{H^3}^2 ) \times e^{C \left( \int_0^t  \Vert \nabla u \Vert_{L^{\infty}} d\tau + \sqrt{M_3(T^{\ast})} t+\sqrt{E(u_0)} t \right) } .\nonumber 
\end{align}
If $\int_0^{T^{\ast}} \Vert \nabla u(\tau) \Vert_{L^{\infty}} d\tau$ is finite, then  $t \in [0,T^{\ast})$
\begin{align}
	\left(  \Vert (u,u_x) \Vert_{H^1}^2+ 	\Vert (\nabla u_{xx} ,\nabla u_{xxx}) \Vert_{L^2}^2 + \Vert ( u_{yy}, u_{xyy})\Vert_{L^2}^2+ \Vert (\nabla u_{yy}, \nabla u_{xyy} )   \Vert_{L^2}^2  \right) + \Vert(u, \partial_x^{-1} u ) \Vert_{H^3}^2 \leq M_5(T^{\ast}), \nonumber 
\end{align}
with 
\begin{align}
M_5 := 	& (  \Vert (u_0,u_{0,x}) \Vert_{H^1}^2+ 	\Vert (\nabla u_{0, xx} ,\nabla u_{0, xxx}) \Vert_{L^2}^2 + \Vert ( u_{0, yy}, u_{0, xyy})\Vert_{L^2}^2+ \Vert (\nabla u_{0, yy}, \nabla u_{0, xyy} )   \Vert_{L^2}^2 \nonumber \\
&+ \Vert (u_0, \partial_x^{-1} u_0) \Vert_{H^3}^2 ) \times e^{C \left( \int_0^{T^{\ast}}  \Vert \nabla u \Vert_{L^{\infty}} d\tau + \sqrt{M_3(T^{\ast})} T^{\ast} + \sqrt{E(u_0)} T^{\ast} \right) }.  \nonumber 
\end{align}
By means of  the fact that 
\begin{align}
	 \Vert u \Vert_{X^3}^2 \leq  3 \left( \Vert (u,u_x) \Vert_{H^1}^2+ 	\Vert (\nabla u_{xx} ,\nabla u_{xxx}) \Vert_{L^2}^2 + \Vert ( u_{yy}, u_{xyy})\Vert_{L^2}^2+ \Vert (\nabla u_{yy}, \nabla u_{xyy} )   \Vert_{L^2}^2 + \Vert (u, \partial_x^{-1} u)\Vert_{H^3}^2 \right), \nonumber 
 \end{align} 
we obtain 
\begin{align}
	\Vert u \Vert_{X^3}^2 \leq 3M_5(T^{\ast}).\label{Rod-KP blow-up criterion eq-14}
\end{align}
Similar to Step 1, we can easily prove  Theorem \ref{Blow-up criterion} for $s=3$. \\
{\bf Step 3: $\mathbf{s>3}$.} In view of (\ref{Rod-KP LWP eq-14}), we have 
\begin{equation}
	\begin{aligned}
	\frac{d}{dt} \left(  2 \Vert u \Vert_{H^s}^2+ \Vert \partial_x^{-1} u  \Vert_{H^s}^2 + \Vert u_x\Vert_{H^s}^2 \right) 
	\leq & C \left(  \Vert u \Vert_{X^2} + \Vert u \Vert_{X^3} + \Vert \nabla u \Vert_{L^{\infty}} \right) \left( 2 \Vert u \Vert_{H^s}^2 + \Vert u_x\Vert_{H^s}^2 + \Vert \partial_x^{-1} u \Vert_{H^s}^2\right)  \nonumber \\
	 \leq & C \left(  \sqrt{2M_3(T^{\ast})}+ \sqrt{3M_5(T^{\ast})} + \Vert \nabla u \Vert_{L^{\infty}} \right)\nonumber \\
	 & \times  \left( 2 \Vert u \Vert_{H^s}^2 + \Vert u_x\Vert_{H^s}^2 + \Vert \partial_x^{-1} u \Vert_{H^s}^2\right) . 
	 \end{aligned} 
\end{equation}
By Gronwall's inequality, we have for all $t\in [0,T^{\ast})$
\begin{align}
	\sup\limits_{\tau\in [0,t]} \left( \Vert u \Vert_{H^s}^2+\Vert \partial_x^{-1} u \Vert_{H^s}^2  + \Vert u_x\Vert_{H^s}^2 \right) \leq  & \left( 2  \Vert u_0\Vert_{H^s}^2 + \Vert \partial_x^{-1}  u_0  \Vert_{H^s}^2+ \Vert u_{0,x} \Vert_{H^s}^2\right) \nonumber \\
	& \times e^{C\left(\int_0^t \Vert \nabla u(\tau) \Vert_{L^{\infty}} d\tau + \sqrt{2M_2(T^{\ast})t} + \sqrt{3M_3(T^{\ast})}t  \right) }. \nonumber 
\end{align}
 If $T^{\ast}<\infty$ satisfies $\int_0^{T^{\ast}} \Vert \nabla u (\tau) \Vert_{L^{\infty}} d\tau <\infty$, we have 
\begin{equation}
	\Vert u \Vert_{X^s}^2 \leq 2 \Vert u \Vert_{H^s}^2 + \Vert \partial_x^{-1} u \Vert_{H^s}^2 + \Vert u_x\Vert_{H^s}^2  \leq M_6(T^{\ast}) , \quad \forall t \in [0,T^{\ast}) ,\nonumber 
\end{equation}
where 
\begin{equation}
	M_6(T^{\ast}) :=	\left(  2 \Vert u_0\Vert_{H^s}^2+ \Vert \partial_x^{-1} u_0 \Vert_{H^s}^2  + \Vert u_{0,x} \Vert_{H^s}^2\right)  e^{C\left(  \int_0^{T^{\ast}} \Vert \nabla u (\tau ) \Vert_{L^{\infty}} d\tau  + \sqrt{M_3(T^{\ast})}T^{\ast}  +\sqrt{ M_4(T^{\ast})} T^{\ast} \right) }. \nonumber 
\end{equation}
which contradicts  the assumption that $T^{\ast} <\infty$ is the maximal existence time. This proves Theorem \ref{Blow-up criterion} for $s>3$. \\
{\bf Step 4: $\mathbf{2<s<3}$.} By the interpolation inequality,  for all $u \in X^3$, we have 
\begin{equation}
	\begin{aligned} 
	\Vert u \Vert_{X^s}^2  =  &  \Vert u\Vert_{H^s}^2 + \Vert u_x \Vert_{H^s}^2 + \Vert \partial_x^{-1} u \Vert_{H^s}^2 \nonumber \\ 
	  \leq  &  \Vert u \Vert_{H^2}^{2(s-2)} \Vert u \Vert_{H^3}^{2(3-s)}  + \Vert u_x \Vert_{H^2}^{2(s-2)} \Vert u_x  \Vert_{H^3}^{2(3-s)}  +\Vert \partial_x^{-1} u\Vert_{H^2}^{2(s-2)} \Vert \partial_x^{-1} u \Vert_{H^3}^{2(3-s)}  \nonumber \\
	  \leq&   3 \Vert u \Vert_{X^2}^{2(s-2)}  \Vert u\Vert_{X^3}^{2(3-s)}.
	\end{aligned} 
\end{equation}
If $\int_0^{T^{\ast}} \Vert \nabla u(\tau ) \Vert_{L^{\infty}} d\tau <\infty$, from (\ref{Rod-KP blow-up criterion eq-89}) and (\ref{Rod-KP blow-up criterion eq-14}), we get 
$\Vert u \Vert_{X^s}^2 \leq 18 M_3(T^{\ast})^{s-2} M_5(T^{\ast})^{3-s} $.
Note that $X^3$ is dense in $X^s$ $(2<s<3)$. A standard density argument yields  
\begin{align}
	\Vert u \Vert_{X^s}^2 \leq 19M_3(T^{\ast})^{s-2} M_5(T^{\ast})^{3-s} , \quad \forall u \in X^s. \nonumber 
\end{align}
This leads to a contradiction, which completes the proof of Theorem \ref{Blow-up criterion} for $2<s<3$. 

Consequently, we have completed the proof of Theorem \ref{Blow-up criterion} from Step 1 to Step 4. 
\end{proof}
%\begin{remark1}\label{blow-up criterion remark}	Under the same assumption in Theorem \ref{Blow-up criterion}, the solution $u$ blows up in finite time $T^{\ast} <+\infty$ if and only if 	\begin{align}	\limsup\limits_{t \to T^{\ast}} \Vert  \sigma \nabla u(t) \Vert_{L^{\infty}(\mathbb{R}^2)} =+\infty. \nonumber \end{align}\end{remark1}

%\begin{theorem1}%\label{BBM blow-up criterion} Suppose that $\sigma =0$, $s\geq 2$, let $u_0 \in X^s(\mathbb{R}^2)$ and $u$ be the corresponding solution to (\ref{the nonlocal weak form}). Assume $T^{\ast}$ is the maximal time of existence. Then  	\begin{align} 	T^{\ast} < \infty \quad \Rightarrow \quad \int_0^{T^{\ast}} \Vert u(\tau ) \Vert_{L^{\infty}(\mathbb{R}^2)} d\tau =\infty.  	\end{align}\end{theorem1}
We now proceed to prove Theorem \ref{BBM blow-up criterion}, which gives the blow-up criterion for $\sigma =0$. %This result is distinct from the one obtained for $\sigma \neq 0$. 
\begin{proof}[Proof of Theorem \ref{BBM blow-up criterion}]
	 Firstly, applying the operator $\Lambda^s \partial_x^{-2}$ to equation (\ref{BBM-KP equation}), and taking the $L^2$ inner-product with $\Lambda^s \partial_x^{-1} u$, we have 
	\begin{align}
		\frac{1}{2} \frac{d}{dt} \left( \Vert \partial_x^{-1} u \Vert_{H^s}^2 + \Vert u\Vert_{H^s}^2\right)  = & \int_{\mathbb{R}^2}  \Lambda^s \partial_x^{-1} u \left(  \partial_t \Lambda^s\partial_x^{-1} u - \partial_t  \Lambda^s u_x \right) dxdy \nonumber \\
		=& -\frac{3}{2}\int_{\mathbb{R}^2} \Lambda^s(u^2) \Lambda^s\partial_x^{-1} u dxdy -\kappa \int_{\mathbb{R}^2} \Lambda^s u_x \Lambda^s \partial_x^{-1} udxdy  \nonumber \\
		& -\int_{\mathbb{R}^2}  \Lambda^s \partial_x^{-2} u_{yy} \Lambda^s  \partial_x^{-1} udxdy  \nonumber \\
		= &-\frac{3}{2}\int_{\mathbb{R}^2} \Lambda^s(u^2) \Lambda^s\partial_x^{-1} u dxdy  . \label{BBM-KP LWP eq-1}
	\end{align}
	Then, applying the operator $\Lambda^s \partial_x^{-1}$ to equation (\ref{BBM-KP equation}), and taking the $L^2$ inner-product with $\Lambda^s u$, we have 
	\begin{align}
		\frac{1}{2} \frac{d}{dt} \left(\Vert u \Vert_{H^s}^2 + \Vert u_x \Vert_{H^s}^2 \right) = & \int_{\mathbb{R}^2} \Lambda^s u (\partial_t \Lambda^s u - \partial_t \Lambda^s u_{xx})dxdy \nonumber \\
		= &- \int_{\mathbb{R}^2}\Lambda^s (uu_x) \Lambda^s u dxdy-\int_{\mathbb{R}^2} \kappa \Lambda^s u_x \Lambda^s u dxdy \nonumber \\
		& -\int_{\mathbb{R}^2} \Lambda^s \partial_x^{-1} u_{yy} \Lambda^s u dxdy \nonumber \\
		= &\frac{3}{2} \int_{\mathbb{R}^2} \Lambda^s (u^2) \Lambda^s u_x dxdy . \label{BBM-KP LWP eq-2}
	\end{align}
	Thanks to the Cauchy-Schwarz inequality, we have 
	\begin{align}
		\left| \int_{\mathbb{R}^2} \Lambda^s(u^2) \Lambda^s\partial_x^{-1} u dxdy \right| \leq C \Vert u^2 \Vert_{H^s} \Vert \partial_x^{-1} u \Vert_{H^s} \leq C \Vert u \Vert_{L^{\infty}} \Vert u \Vert_{H^s} \Vert \partial_x^{-1} u \Vert_{H^s},\label{BBM-KP LWP eq-3}
	\end{align}
	\begin{align}
		\left| \int_{\mathbb{R}^2} \Lambda^s (u^2) \Lambda^s u_x dxdy \right| \leq C \Vert u^2 \Vert_{H^s} \Vert u_x \Vert_{H^s} \leq C \Vert u \Vert_{L^{\infty}} \Vert u \Vert_{H^s} \Vert u_x \Vert_{H^s}. \label{BBM-KP LWP eq-4}
	\end{align}
     Thus, by virtue of (\ref{BBM-KP LWP eq-1})-(\ref{BBM-KP LWP eq-4}), we have 
	\begin{align}
		\frac{d}{dt} \left( \Vert \partial_x^{-1} u\Vert_{H^s}^2 + 2\Vert u \Vert_{H^s}^2 + \Vert u_x \Vert_{H^s}^2 \right) \leq C \Vert u \Vert_{L^{\infty}} \left( \Vert \partial_x^{-1} u\Vert_{H^s}^2 + 2\Vert u \Vert_{H^s}^2 + \Vert u_x \Vert_{H^s}^2\right) . \label{BBM-KP LWP eq-5}
	\end{align}
	Applying Gronwall's inequality to (\ref{BBM-KP LWP eq-5}) yields that for all $t \in [0,T^{\ast})$
	\begin{align}
		\sup\limits_{\tau \in [0,t]}\left( \Vert u (\tau) \Vert_{X^s}^2 \right) \leq & \sup\limits_{\tau\in [0,t]} \left( \Vert \partial_x^{-1} u (\tau)\Vert_{H^s}^2 + 2\Vert u (\tau)\Vert_{H^s}^2 + \Vert u_x(\tau) \Vert_{H^s}^2 \right) \nonumber \\
		\leq & \left( \Vert \partial_x^{-1} u_0\Vert_{H^s}^2 + 2\Vert u_0\Vert_{H^s}^2 + \Vert u_{0,x} \Vert_{H^s}^2 \right) e^{\int_0^t \Vert u (\tau)\Vert_{L^{\infty}} }. \nonumber
	\end{align}
	If $\int_0^{T^{\ast}} \Vert u (\tau ) \Vert_{L^{\infty}}d\tau $ is finite, then for all $t\in [0,T^{\ast})$, we have 
	\begin{align}
		\Vert u (t)\Vert_{X^s}^2 \leq \left( \Vert \partial_x^{-1} u_0\Vert_{H^s}^2 + 2\Vert u_0\Vert_{H^s}^2 + \Vert u_{0,x} \Vert_{H^s}^2 \right) e^{\int_0^{T^{\ast}}  \Vert u (\tau)\Vert_{L^{\infty}} d\tau }.\nonumber 
	\end{align}
This completes the proof of Theorem \ref{BBM blow-up criterion}.
\end{proof}

\section{Global Existence in the case $\sigma =0$}\label{Global BBM}
In this section,  we investigate global existence  for the BBM-KP equation (\ref{BBM-KP equation}). %When $\sigma =0$, the conservation law $F_{\sigma}(u_0)$ in (\ref{Rod-KP conserved quantity-1}) enables us to prove the existence of a global solution for all $u_0\in X^s$, $s\geq 2$. 
%\begin{theorem2}%\label{BBM global result}
%	Assume $\sigma =0$, let $u_0 \in X^s(\mathbb{R}^2)$, $s\geq 2$. Then there exists a unique global strong solution $u\in C([0,T]; X^s(\mathbb{R}^2)) \cap C^1([0,T]; X^{s-2}(\mathbb{R}^2))$.  \end{theorem2}
\begin{proof}[Proof of Theorem \ref{BBM global result}]
	Applying the operator $\partial_x^{-2} \partial_y$ to (\ref{BBM-KP equation}), we have 
	\begin{align}
		\partial_x^{-1} u_{ty} - u_{txy} + 3uu_y + \kappa u_y + \partial_x^{-2} u_{yyy} =0. \label{BBM-KP global eq-1}
	\end{align}
Then taking the $L^2$ inner-product between (\ref{BBM-KP global eq-1}) and $\partial_x^{-1} u_y$ to give 
	\begin{align}
		\frac{1}{2}\frac{d}{dt} \left(  \Vert \partial_x^{-1} u_y \Vert_{L^2}^2 + \Vert u_y \Vert_{L^2}^2\right) =&\int_{\mathbb{R}^2} \partial_x^{-1} u_y \left( \partial_x^{-1} u_{ty} - u_{txy}  \right) dxdy \nonumber \\
		= & -\int_{\mathbb{R}^2} \left( 3uu_y \partial_x^{-1} u_y  + \kappa u_y \partial_x^{-1} u_y + \partial_x^{-2} u_{yyy} \partial_x^{-1} u_y \right) dxdy  \nonumber \\
		= & -3\int_{\mathbb{R}^2} uu_y \partial_x^{-1} u_y dxdy . \label{BBM global result proof eq-01}
	\end{align}
	Thanks to the Cauchy-Schwartz inequality and Lemma \ref{Sobolev embedding remark-2}, we have 
	\begin{align}
		\left| \int_{\mathbb{R}^2} uu_y \partial_x^{-1} u_y dxdy\right| \leq &\Vert u\partial_x^{-1} u_y\Vert_{L^2} \Vert u_y \Vert_{L^2}  \nonumber \\
		\leq & C\Vert u \Vert_{L^2}^{\frac{1}{2}}\Vert \partial_x^{-1} u_y \Vert_{L^2}^{\frac{1}{2}} \Vert u_y \Vert_{L^2}^2.\label{BBM global result proof eq-02}  
	\end{align}
In view of Lemma \ref{Sobolev inequality lemma}-(\ref{Sobolev inequality lemma-4}), we have 
\begin{align}
	\Vert \partial_x^{-1} u_y \Vert_{L^2}^2 =& 2F_0(u_0) -\kappa \int_{\mathbb{R}^2} u^2 dxdy -\int_{\mathbb{R}^2} u^3 dxdy  \nonumber \\
	\leq & 2F_0(u_0) +2 |\kappa| E(u_0)+ \left| \int_{\mathbb{R}^2} u^3dxdy\right|  \nonumber \\
	\leq & 2F_0(u_0)+ 2|\kappa| E(u_0)+C \Vert u \Vert_{L^2}^{\frac{3}{2}} \Vert u_x \Vert_{L^2}\Vert \partial_x^{-1} u_y \Vert_{L^2}^{\frac{1}{2}} \nonumber \\
	\leq &2F_0(u_0)+2 |\kappa| E(u_0)+2C E(u_0)\Vert \partial_x^{-1} u_y \Vert_{L^2}^{\frac{1}{2}}  \nonumber \\
	\leq & 2F_0(u_0)+2 |\kappa| E(u_0)+CE^2(u_0) +C \Vert \partial_x^{-1} u_y \Vert_{L^2} .\nonumber 
	\end{align}
Thus, we have
\begin{align}
	\Vert \partial_x^{-1} \partial_y \Vert_{L^2} \leq  \sqrt{2F_0(u_0)+ 2|\kappa|E(u_0) + CE^2(u_0)+\frac{C^2}{4}}-\frac{C}{2} . \label{BBM global result proof eq-03}
\end{align}
	Therefore, combining (\ref{BBM global result proof eq-01}), (\ref{BBM global result proof eq-02}) and (\ref{BBM global result proof eq-03}) yields
	\begin{align}
		\frac{d}{dt} \left(  \Vert \partial_x^{-1} u_y \Vert_{L^2}^2 + \Vert u_y \Vert_{L^2}^2\right)  \leq C E^{\frac{1}{4}} (u_0)C_1^{\frac{1}{2}}\left(  \Vert \partial_x^{-1} u_y \Vert_{L^2}^2 + \Vert u_y \Vert_{L^2}^2\right) , \nonumber 
	\end{align}
where $C_1 := \sqrt{2F_0(u_0)+ 2|\kappa|E(u_0) + CE^2(u_0)+\frac{C^2}{4}}-\frac{C}{2} >0$. 
	Applying the Gronwall inequality, one may find, for all $t \in [0,T^{\ast})$
	\begin{align}
		\Vert u_y (t)\Vert_{L^2}^2 \leq & \Vert \partial_x^{-1} u_y (t)\Vert_{L^2}^2 + \Vert u_y(t) \Vert_{L^2} \nonumber \\
		\leq & \left( \Vert \partial_x^{-1} u_{0,y}\Vert_{L^2}^2 + \Vert u_{0,y}\Vert_{L^2}^2 \right) e^{CE^{\frac{1}{4}}(u_0) C_1^{\frac{1}{2}} t}. \label{BBM global result proof eq-1}
	\end{align}
Applying the operator $\partial_x^{-1} \partial_y^2$ to (\ref{BBM-KP equation}), we have 
\begin{align}
	u_{ty}-u_{txxy} + 3u_y u_x + 3uu_{xy} + \kappa u_{xy} +\partial_x^{-1} u_{yyy} =0. \label{BBM-KP global eq-2}
\end{align}
	Taking the $L^2$ inner-product between (\ref{BBM-KP global eq-2}) and $u_y$, we get 
	\begin{align}
		\frac{1}{2} \frac{d}{dt} \left(\Vert u_y \Vert_{L^2}^2 + \Vert u_{xy} \Vert_{L^2}^2 \right)  = & \int_{\mathbb{R}^2}u_y (u_{ty} - u_{txxy}) dxdy \nonumber \\
		= & -\int_{\mathbb{R}^2} \left( 3u_y^2 u_x +3 u u_y u_{xy} +\kappa u_{xy} u_y + \partial_x^{-1} u_{yyy} u_y \right) dxdy \nonumber \\
		= & 3\int_{\mathbb{R}^2} uu_y u_{xy} dxdy . \label{BBM global result proof eq-21}
	\end{align}
	Thanks to the Cauchy-Schwartz inequality and Lemma \ref{Sobolev embedding remark-2} again, we have 
	\begin{align}
		\left| \int_{\mathbb{R}^2} uu_y u_{xy} dxdy \right| \leq  & \Vert uu_y \Vert_{L^2}\Vert u_{xy} \Vert_{L^2}\nonumber \\
		\leq & C\Vert u \Vert_{L^2}^{\frac{1}{2}} \Vert u_y \Vert_{L^2} \Vert u_{xy} \Vert_{L^2}^{\frac{3}{2}}. \label{BBM global result proof eq-22}
	\end{align}
	Thus, combinbing  (\ref{BBM global result proof eq-1}), (\ref{BBM global result proof eq-21}) and (\ref{BBM global result proof eq-22}), it follows that 
	\begin{align}
		&\frac{d}{dt} (\Vert u_y \Vert_{L^2}^2 + \Vert u_{xy} \Vert_{L^2}^2 )\nonumber \\
		\leq & C \Vert u \Vert_{L^2}^{\frac{1}{2}} \Vert u_y \Vert_{L^2}^{\frac{1}{2}} \left(  \Vert u_y \Vert_{L^2}^2+ \Vert u_{xy} \Vert_{L^2}^2\right)\nonumber  \\
		 \leq & CE^{\frac{1}{4}}(u_0) \left( \Vert \partial_x^{-1} u_{0,y} \Vert_{L^2}^2 + \Vert u_{0,y} \Vert_{L^2}^2 \right)^{\frac{1}{4}} e^{CE^{\frac{1}{4}}(u_0)  C_1^{\frac{1}{2}} t}(\Vert u_y \Vert_{L^2}^2 + \Vert u_{xy} \Vert_{L^2}^2 ) . \label{BBM global eq-3}
	\end{align}
	Using Gronwall's inequality to (\ref{BBM global eq-3}) yields that 
	\begin{align}
		&\sup\limits_{t\in [0,T^{\ast})} \left( \Vert u_y (t) \Vert_{L^2}^2 + \Vert u_{xy} (t) \Vert_{L^2}^2\right) \nonumber \\
		\leq &\left( \Vert u_{0,y} \Vert_{L^2}^2 + \Vert u_{0,xy}  \Vert_{L^2}^2\right)\left( \Vert \partial_x^{-1} u_{0,y} \Vert_{L^2}^2 + \Vert u_{0,y} \Vert_{L^2}^2 \right)^{\frac{1}{4}} CE^{\frac{1}{4}}(u_0)T^{\ast}e^{CE^{\frac{1}{4}} (u_0) C_1^{\frac{1}{2}} T^{\ast}}. \nonumber 
	\end{align}
From  Lemma \ref{Sobolev inequality lemma}-(\ref{Sobolev inequality lemma-1}), the following estimate is derived 
\begin{align}
	\sup\limits_{t\in [0,T^{\ast})} \Vert u(t)  \Vert_{L^{\infty}}^2 \leq C \Vert u(t) \Vert_{L^2}^{\frac{1}{2}}  \Vert u_x(t) \Vert_{L^2}^{\frac{1}{2}} \Vert u_y(t) \Vert_{L^2}^{\frac{1}{2}} \Vert u_{xy} (t) \Vert_{L^2}^{\frac{1}{2}} \leq CE^{\frac{1}{2}} (u_0) C_2^{\frac{1}{2}}, \nonumber 
\end{align}
where 
\begin{align}
	C_2: =\left( \Vert u_{0,y} \Vert_{L^2}^2 + \Vert u_{0,xy}  \Vert_{L^2}^2\right)\left( \Vert \partial_x^{-1} u_{0,y} \Vert_{L^2}^2 + \Vert u_{0,y} \Vert_{L^2}^2 \right)^{\frac{1}{4}} CE^{\frac{1}{4}}(u_0)T^{\ast}e^{CE^{\frac{1}{4}}(u_0)  C_1^{\frac{1}{2}} T^{\ast}}. \nonumber 
\end{align}
This completes the proof of Theorem \ref{BBM global result}.
\end{proof}

%\begin{theorem2}%\label{blow-up data-1}	Fix $\phi \in H^2(\mathbb{R})$ with $\phi \geq 0$ and $\int_{\mathbb{R}} \phi dy=1$. Suppose that $\sigma \neq 0$, $u_0 \in X^s(\mathbb{R}^2)$ with $s\geq 2$. If  the corresponding solution $u$ satisfies that $u(t,x,y) =-u(t,-x,y)$ and the initial value $u_0$ is subject to the constraint  \begin{align} 	\int_{\mathbb{R}} \sigma u_{0,x}(0,y) \phi(y) dy < -\sqrt{2C_2}, \nonumber 	\end{align} where  \begin{equation}	C_3:= \left\{ 	\begin{aligned} 	& \frac{\sigma(\sigma-3)}{2}  E(u_0) \Vert \phi\Vert_{L^{\infty}},  && \quad\sigma >4, \\ 	&0, & &\quad 0\leq \sigma \leq 4,\\ 	& \frac{\sigma(\sigma-3)}{2}E(u_0) \Vert \phi\Vert_{L^{\infty}},  & & \quad \sigma <0.	\end{aligned}	\right.\end{equation}
%Then the corresponding strong solution $u$ blows up in finite time $T^{\ast}$ such that  \begin{equation}	T^{\ast} \leq T_1 := \left\{ 	\begin{aligned} 	& \frac{1}{\sqrt{2C_3}}\ln \left( \frac{\int_{\mathbb{R}} \sigma u_{0,x}(0,y) \phi (y) dy -\sqrt{2C_3} }{\int_{\mathbb{R}} \sigma u_{0,x}(0,y) \phi (y) dy +\sqrt{2C_3}} \right) , && \quad \text{$\sigma>4$ or $\sigma<0$}, \\ 	& -\frac{2}{\int_{\mathbb{R}} u_{0,x} (0,y) dy}, && \quad 0<\sigma \leq 4. \end{aligned}	\right.\nonumber \end{equation}\end{theorem2}
\section{Blow-up  data for the case $\sigma \neq 0$}\label{Blow-up rod}
We now in this section focus our attention on blow-up data for the rod-KP equation (\ref{the nonlocal weak form}) when $\sigma \neq 0$. The lemma stated below is useful in the proof of Theorem \ref{blow-up data-1}. 
\begin{lemma2}[\cite{Brandloese CMP 2014}]\label{Brandloese CMP 2014 lemma}
	Let $u \in H^1(\mathbb{R})$ and $ 0 \leq \sigma \leq 4$. Then 
	\begin{align}
		p \ast \left( \frac{3-\sigma}{2} u^2 + \frac{\sigma}{2} u_x^2\right) \geq J(\sigma) u^2, \nonumber
	\end{align}
	where $J(\sigma) =\frac{\sqrt{\sigma}}{2} \left( \sqrt{12-3\sigma} -\sqrt{\sigma} \right) \leq \frac{3-\sigma}{2}$ and $p(x) = \frac{1}{2} e^{-|x|}$ is the fundmental solution of $(1-\partial_x^2)^{-1}$ on the real line $\mathbb{R}$.
\end{lemma2}
\begin{proof}[Proof of Theorem \ref{blow-up data-1}]
	In view of the assumption $u(t,x,y)=-u(t,-x,y)$, $\forall t\geq 0$, $(x,y)\in \mathbb{R}^2$, we have $u(t,0,y)=0$, $(p \ast u) (t,0,y) =0$ and $(p \ast u_{yy} )(t,0,y)=0$. Multiplying by $\phi(y)$ on the both sides of (\ref{the nonlocal weak form}), differentiating it with respect to $x$ and integrating it over $\mathbb{R}$ with $y$, it yields that 
	\begin{align}
		\int_{\mathbb{R}} \left( u_{tx} +\sigma uu_{xx}\right) \phi dy =& \int_{\mathbb{R}} \left( -\frac{\sigma}{2} u_x^2 + \frac{3-\sigma}{2} u^2 + \kappa u \right) \phi dy \nonumber \\
		& -\int_{\mathbb{R}} p \ast \left( \frac{3-\sigma}{2} u^2 +\frac{\sigma}{2} u_x^2 + \kappa u \right) \phi dy -\int_{\mathbb{R}} p \ast u_{yy} \phi dy. \nonumber 
	\end{align}
At the point $(t,0,y)$, we have 
\begin{align}
	\frac{d}{dt} \int_{\mathbb{R}} u_x (t,0,y) \phi (y) dy = &- \frac{\sigma}{2} \int_{\mathbb{R}} u_x^2(t,0,y) \phi (y) dy -\int_{\mathbb{R}} p \ast \left( \frac{3-\sigma}{2} u^2 + \frac{\sigma}{2} u_x^2  \right) (t,0,y) \phi(y) dy . \nonumber 
\end{align}
Taking advantage of H\"{o}lder inequality and the conservation law $E(u_0)$ in (\ref{Rod-KP conserved quantity}), we get 
\begin{align}
	\left( \int_{\mathbb{R}} u_x(t,0,y) \phi(y) dy\right)^2 \leq \int_{\mathbb{R}} u_x^2(t,0,y) \phi(y) dy , \nonumber 
\end{align}
\begin{align}
	\left| \int_{\mathbb{R}} (p \ast u^2) (t,0,y) \phi(y) dy \right| \leq \int_{\mathbb{R}} \Vert p \ast u^2 \Vert_{L^{\infty}_x} dy \Vert \phi \Vert_{L^{\infty}}\leq E(u_0)\Vert \phi \Vert_{L^{\infty}}, \nonumber 
\end{align}
\begin{align}
		\left| \int_{\mathbb{R}} (p \ast u_x^2) (t,0,y) \phi(y) dy \right| \leq \int_{\mathbb{R}} \Vert p \ast u_x^2 \Vert_{L^{\infty}_x} dy \Vert \phi \Vert_{L^{\infty}}\leq E(u_0)\Vert \phi \Vert_{L^{\infty}}, \nonumber 
\end{align} 
%\begin{align} 	\left| \int_{\mathbb{R}} \left( p\ast u\right) (t,0,y) \phi(y) dy \right| \leq \int_{\mathbb{R}} \Vert p \ast u \Vert_{L^{\infty}_x} \phi(y) dy \leq H_1 \Vert \phi \Vert_{L^{\infty}}. \nonumber  \end{align}
If $\sigma >4$, at the point $(t,0,y)$, we have 
\begin{align}
	\frac{d}{dt} \int_{\mathbb{R}} u_x(t,0,y) \phi(y) dy \leq &-\frac{\sigma}{2} \left( \int_{\mathbb{R}} u_x (t,0,y) \phi(y) dy \right)^2 +\frac{\sigma -3}{2} \int_{\mathbb{R}} \left( p \ast u^2 \right) (t,0,y) \phi(y) dy\nonumber \\
	\leq & -\frac{\sigma}{2} \left( \int_{\mathbb{R}} u_x (t,0,y) \phi(y) dy \right)^2 + \frac{\sigma-3}{2} E(u_0)\Vert \phi\Vert_{L^{\infty}}. \nonumber 
	\end{align} 
If $ 0 <\sigma \leq 4$, in view of Lemma \ref{Brandloese CMP 2014 lemma} and the fact that $u(t,0,y)=0$, we have 
\begin{align}
		\frac{d}{dt} \int_{\mathbb{R}} u_x(t,0,y) \phi(y) dy \leq &-\frac{\sigma}{2} \left( \int_{\mathbb{R}} u_x (t,0,y) \phi(y) dy \right)^2 +J(\sigma) \int_{\mathbb{R}} u^2(t,0,y) \phi(y) dy \nonumber \\
		\leq &-\frac{\sigma}{2} \left( \int_{\mathbb{R}} u_x (t,0,y) \phi(y) dy \right)^2 . \nonumber   
\end{align}
If $\sigma <0$, we have 
\begin{align}
	\frac{d}{dt} \int_{\mathbb{R}} u_x(t,0,y) \phi(y) dy \geq & -\frac{\sigma}{2} \left( \int_{\mathbb{R}} u_x(t,0,y) \phi(y)\right)^2 -\frac{3-\sigma}{2} \int_{\mathbb{R}} (p \ast u^2)(t,0,y) \phi(y) \nonumber \\
	\geq & -\frac{\sigma}{2} \left( \int_{\mathbb{R}} u_x(t,0,y) \phi(y)\right)^2 - \frac{3-\sigma}{2}E(u_0)\Vert \phi \Vert_{L^{\infty}}. \nonumber 
\end{align}
Then, we denote
\begin{equation}
	C_3:= \left\{ 
	\begin{aligned}
	& \frac{\sigma(\sigma-3)}{2} E(u_0)\Vert \phi\Vert_{L^{\infty}},  && \quad\text{$\sigma >4$ or $\sigma <0$}, \\
	&0,  & &\quad 0< \sigma \leq 4.
	\end{aligned}
	\right.
\end{equation}
Let $ B(t) : =\int_{\mathbb{R}} \sigma u_x(t,0,y) \phi(y)dy$,  we have 
\begin{align}
	\frac{d}{dt} B(t)\leq -\frac{1}{2} B^2(t) + C_3. \label{blow-up data eq-1}
\end{align}
Since the initial data satisfies $\int_{\mathbb{R}} \sigma u_{0,x} (0,y) \phi (y) dy < -\sqrt{2C_3}$, we deduce that 
\begin{align}
	\lim\limits_{t \to T^{\ast}} \int_{\mathbb{R}} \sigma u_x(t,0,y) \phi(y) dy =-\infty, \nonumber 
\end{align}
with 
\begin{equation}
	T^{\ast} \leq T_1 := \left\{ 
	\begin{aligned}
		& \frac{1}{\sqrt{2C_3}}\ln \left( \frac{\int_{\mathbb{R}} \sigma u_{0,x}(0,y) \phi (y) dy -\sqrt{2C_3 }}{\int_{\mathbb{R}} \sigma u_{0,x}(0,y) \phi (y) dy +\sqrt{2C_3}} \right) , && \quad \text{$\sigma>4$ or $\sigma<0$}, \\
		& -\frac{2}{\int_{\mathbb{R}}\sigma u_{0,x} (0,y)\phi(y)  dy}, && \quad 0<\sigma \leq 4. 
	\end{aligned}
\right.\nonumber 
\end{equation}
This completes the proof of Theorem \ref{blow-up data-1}. 
\end{proof}
%\begin{theorem2}Suppose that $\sigma \neq 0$, $u_0 \in X^s(\mathbb{R})$ with $s\geq 2$. Let $T^{\ast}$ be the maximal existence time of  the corresponding solution $u$ to (\ref{the nonlocal weak form}) with the initial value $u_0$. Assume that \end{theorem2}
\begin{remark2}
	In the case $(x,y)\in \mathbb{R}\times \mathbb{T}$, we have a similar blow-up result. Suppose that $ \sigma \neq 0$,  and $u_0 \in X^s(\mathbb{R} \times \mathbb{T})$ with $s\geq 2$. Let $T^{\ast}$ be the maximal existence time of the corresponding solution $u$ to (\ref{the nonlocal weak form}) with initial data $u_0$. If the solution $u$ satisfies that $u(t,x,y) =-u(t,-x,y)$, $\forall t\geq 0$, $(x,y) \in \mathbb{R} \times \mathbb{T}$, and the initial data $u_0$ satisfies 
	\begin{align}
		\int_{\mathbb{T}} \sigma u_{0,x} (0,y) dy < -\sqrt{2C_4}, \nonumber 
	\end{align}
where 
\begin{equation}
	C_4:= \left\{ 
	\begin{aligned}
		& \frac{\sigma(\sigma -3)}{2}E(u_0), && \quad \text{$\sigma >4$ or $\sigma <0$},\\
		& 0, && \quad 0< \sigma \leq 4.
	\end{aligned}
\right. \nonumber 
\end{equation}
Then  the corresponding strong solution $u$ blows up in finite time $T^{\ast}$ such that 
\begin{equation}
	T^{\ast} \leq T_2:= \left\{ 
	\begin{aligned}
		&\frac{1}{\sqrt{2C_4}} \ln \left(\frac{\int_{\mathbb{T}} \sigma u_{0,x}(0,y) dy -\sqrt{2C_4}}{\int_{\mathbb{T}} \sigma u_{0,x} (0,y) dy +\sqrt{2C_4}}\right) ,&&\quad \text{$\sigma>4$ or $\sigma <0$},\\
		& -\frac{2}{\int_{\mathbb{T}}\sigma u_{0,x}(0,y) dy}, && \quad 0< \sigma \leq 4. 
	\end{aligned}
\right.\nonumber 
\end{equation}
\end{remark2}
\section{Uniqueness continuation property of the solutions}\label{Liouville-type property of the solution}
\newtheorem{remark4}{Remark}[section]
\newtheorem{theorem4}{Theorem}[section]
\newtheorem{lemma4}{Lemma}[section]
\newtheorem{corollary4}{Corollary}[section]
\newtheorem{proposition4}{Proposition}[section]
This section is devoted to the study of the uniqueness continuation properity of  solutions to  rod-KP equation (\ref{the nonlocal weak form}). 
%Now we complete the proof of Theorem \ref{UCP Rod}. 
\begin{proof}[Proof of Theorem \ref{UCP Rod}]
	%The proof is given by a contradiction argument. We claim that if there is such a domain $\Omega$, then the solution is trivial. 
	
	Firstly, integrating (\ref{the nonlocal weak form}) over $\mathbb{R}$  with $y$, it yields that 
	\begin{align}
		\int_{\mathbb{R}} \left(  u_t + \sigma uu_x \right) dy = -\int_{\mathbb{R}} p_x \ast \left( \frac{3-\sigma}{2} u^2 +\frac{\sigma}{2} u_x^2 \right)  dy .
		\label{sec:6 theorem eq-1}
	\end{align}
	By virtue of the assumption and (\ref{sec:6 theorem eq-1}), we see that 
	\begin{align}
		\frac{3-\sigma}{2}u^2 +\frac{\sigma}{2} u_x^2 =0 , \quad \forall (t,x) \in \Omega, \, y \in \mathbb{R},
		\nonumber 
	\end{align}
	and 
	\begin{align}
		\int_{\mathbb{R}}  p_x \ast \left( \frac{3-\sigma}{2}u^2 +\frac{\sigma}{2} u_x^2\right) dy =0 ,\quad \forall (t,x) \in \Omega.\nonumber
	\end{align}
	Since $\Omega$ is an open set, there exists a $t^{\ast} \in (0,T^{\ast})$ and $I=[a,b]$, $a<b$ such that $\left\lbrace t^{\ast} \right\rbrace \times I \subset \Omega$. We denote 
	\begin{align}
		F(x,y) := p_x \ast \left( \frac{3-\sigma}{2}u^2 +\frac{\sigma}{2} u_x^2\right) \left( t^{\ast},x,y\right) ,\nonumber
	\end{align}
	and 
	\begin{align}
		f(x,y) := \left( \frac{3-\sigma}{2}u^2 +\frac{\sigma}{2} u_x^2\right) \left( t^{\ast},x,y\right) \geq 0,\nonumber
	\end{align}
	with $0\leq \sigma\leq 3$. The regularity of $u$ ensures that $F$ and $f$ are smooth functions. From the above argument, we observe that $\int_{\mathbb{R}} F(x,y) dy =0$, $f(x,y)=0$, $\forall x\in [a,b]$. 
	
	On the other hand, for any $y \in \mathbb{R}$, we have 
	\begin{align}
		F(b,y) = & -\frac{e^{-b}}{2}\int_{-\infty}^b e^z f(z,y) dz +\frac{e^b}{2} \int_b^{+\infty} e^{-z} f(z,y) dz \nonumber \\
		= & -\frac{e^{-b}}{2}\int_{-\infty}^a e^z f(z,y) dz +\frac{e^b}{2} \int_a^{+\infty} e^{-z} f(z,y) dz \nonumber \\
		\geq &  -\frac{e^{-a}}{2}\int_{-\infty}^a e^z f(z,y) dz +\frac{e^a}{2} \int_a^{+\infty} e^{-z} f(z,y) dz \nonumber \\
		=&  F(a,y). \nonumber 
	\end{align}
	Since $\int_{\mathbb{R}} \left( F(b,y)  - F(a,y) \right) dy = 0-0=0$ and $F(b,y)- F(a,y) \geq 0$, it follows that $F(b,y) -F(a,y) =0$, which implies that 
	\begin{align}
		\frac{e^{-a} -e^{-b}}{2} \int_{-\infty}^a e^z f(z,y) dz + \frac{e^b -e^a}{2} \int_a^{+\infty} e^{-z} f(z,y) dz =0 ,\nonumber 
	\end{align}
	thus $f(x,y) =0$, $\forall (x,y) \in \mathbb{R}^2$. By the uniqueness of the solution, we obtain the desired result.
\end{proof}
%\noindent {\bf Acknowledgments.}  
%	Let $\varphi (x-ct)$ be a solitary wave of  Rod equation with speed $c\in \mathbb{R}$. Then $\varphi (x)$ satisfy \begin{align}	\left( -c+\kappa \right) \varphi_x + c \varphi_{xxx} + 3 \varphi \varphi_x = \sigma \left( 2\varphi_x \varphi_{xx} + \varphi \varphi_{xxx} \right) . \label{smooth wave solution eq-1} 	\end{align}
% After integration we get  \begin{align} 	\left( -c+\kappa \right) \varphi + c\varphi_{xx} + \frac{3}{2} \varphi^2 = \sigma \left( \varphi \varphi_{xx} +\frac{1}{2} \varphi_x^2 \right) .\label{smooth wave solution eq-2} \end{align}
% Multiplying both sides of (\ref{smooth wave solution eq-2}) by $2 \varphi_x$ and integrating on $(-\infty,x]$, and applying the integartion by parts, it follows that  \begin{align} 	\varphi_x^2 \left( c-\sigma \varphi \right)  = \varphi^2\left( c-\kappa -\varphi \right) . \nonumber  \end{align}

%Notice that   the Rod-KP equation (\ref{Rod-KP equation form})  recovers  the CH-KP equation  (\ref{CH-KP eqaution}) in the special case $\sigma=0$.                                                                                                                                                                
\section{Travelling wave  solutions} 
\label{Travelling wave  solutions}
In this section, we focus on the travelling wave solutions of the rod-KP equation (\ref{Rod-KP equation form-1}). 
Let us consider  the travelling wave solution  in the form of $u(t,x,y) = u(x-ct,y)$.   %By Theorem \ref{Travelling wave solution rigidity theorem}, we obtain that the solution  to (\ref{travelling wave solution eq-1}) is trivial in the case $c>\min \left\lbrace0, \kappa  \right\rbrace $.  The detailed proof of Theorem \ref{Travelling wave solution rigidity theorem} as follows 
\begin{proof}[Proof of Theorem \ref{Travelling wave solution rigidity theorem}]
	Using the fact that $2u_xu_{xx} + uu_{xxx} = \left(  \frac{1}{2} u_x^2 + uu_{xx}\right)_x$, we can rewrite (\ref{travelling wave solution eq-1}) as follows: 
	\begin{align}
		\left( \left( -c+ \kappa \right) u + cu_{xx} + \frac{3}{2} u^2 - \sigma \left( \frac{1}{2} u_x^2 + uu_{xx}\right)  \right)_x + \partial_x^{-1} u_{yy} =0 . \label{travelling wave solution eq-2}
	\end{align} 
	Multiply (\ref{travelling wave solution eq-2}) by $ y \partial_x^{-1} u_y$ and integrate on $\mathbb{R}^2$, we have 
	\begin{align}
		\int_{\mathbb{R}^2} (- c+\kappa )yuu_y +cy u_{xx} u_y +\frac{3}{2} y u^2 u_y + \sigma y\left( \frac{1}{2} u_x^2 + uu_{xx}\right) u_y +\frac{1}{2} y \left( \partial_x^{-1} u_y \right)^2 dxdy =0 .  \label{travelling wave solution eq-3}
	\end{align}
	It follows from the integration by parts that 
	\begin{align}
		\int_{\mathbb{R}^2} (-c+\kappa)y uu_y  + c y u_{xx} u_y   =\frac{1}{2} \int_{\mathbb{R}^2} \left( (c-\kappa) u^2+ c  u_x^2  \right)  dxdy , 
		\label{travelling wave solution eq-4}
	\end{align}
	\begin{align}
		\int_{\mathbb{R}^2} \frac{3}{2} y u^2 u_y dxdy =-\frac{1}{2}  \int_{\mathbb{R}^2}  u^3dxdy  , 	\label{travelling wave solution eq-4.5}
	\end{align}
	\begin{align}
		\int_{\mathbb{R}^2} y \left( \frac{1}{2} u_x^2 + uu_{xx} \right) u_y dxdy =&  -\int_{\mathbb{R}^2} \frac{1}{2} yu_x^2 u_y + y uu_x u_{xy} dxdy  \nonumber \\
		=& -\frac{1}{2} \int_{\mathbb{R}^2} y (uu_x^2)_y dxdy =\frac{1}{2}\int_{\mathbb{R}^2} uu_x^2dxdy.\label{travelling wave solution eq-5}
	\end{align}
	Combining (\ref{travelling wave solution eq-3})-(\ref{travelling wave solution eq-5}), one gets 
	\begin{align}
		\frac{1}{2} \int_{\mathbb{R}^2} \left(  -( c -\kappa)u^2 -c  u_x^2 +u^3 + \sigma uu_x^2 - \left( \partial_x^{-1} u_y \right)^2  \right) dxdy=0  . \label{travelling wave solution eq-6}
	\end{align}
	Applying the operator $\partial_x^{-1}$ to (\ref{travelling wave solution eq-2}), we have 
	\begin{align}
		(-c+\kappa )u + cu_{xx} +\frac{3}{2} u^2 -\sigma \left(  \frac{1}{2} u_x^2 + uu_{xx} \right) 
		+ \partial_x^{-2} u_{yy} =0. \label{travelling wave solution eq-7}
	\end{align}
	Multiply (\ref{travelling wave solution eq-7}) by $u$ and integrate on $\mathbb{R}^2$, we have 
	\begin{align}
		\int_{\mathbb{R}^2} \left( (-c+\kappa)u^2 -cu_x^2 +\frac{3}{2} u^3-\sigma \left(  u^2 u_{xx}  +\frac{1}{2} uu_x^2 \right) + \left(\partial_x^{-1} u_y  \right)^2 \right) dxdy =0 . \label{travelling wave solution eq-8}
	\end{align}
	Thanks to integration by parts again, we get 
	\begin{align}
		\int_{\mathbb{R}^2} u^2 u_{xx} dxdy = -2\int_{\mathbb{R}^2} uu_x^2 dxdy . \label{travelling wave solution eq-9}
	\end{align}
	Taking (\ref{travelling wave solution eq-9}) into (\ref{travelling wave solution eq-8}), one yields
	\begin{align}
		\int_{\mathbb{R}^2} \left( (-c+\kappa)u^2 -cu_x^2 +\frac{3}{2} u^3+ \frac{3\sigma}{2}   uu_x^2  + \left(\partial_x^{-1} u_y  \right)^2 \right) dxdy =0 . \label{travelling wave solution eq-10}
	\end{align}
	From (\ref{travelling wave solution eq-6}) and (\ref{travelling wave solution eq-10}), we obtain 
	\begin{align}
		\int_{\mathbb{R}^2}2(c-\kappa) u^2+ 2cu_x^2 + 4\left( \partial_x^{-1} u_y \right)^2 dxdy =0 . \nonumber 
	\end{align}
	Since $c>\max \left\lbrace0,  \kappa\right\rbrace $, it then follows that $u=0$. This completes the proof of Theorem \ref{Travelling wave solution rigidity theorem}. 
\end{proof}
Next, we  present the  symmetry property of the solitary wave solutions to (\ref{travelling wave solution eq-1}).
%  of the rotation-modified Kadomtsev-Petviashvili equation in \cite{Esfahani-Levandosky PAMS 2019}, we are able to present  the proof of Theorem \ref{Symmetry result}. 
\begin{proof}[Proof of Theorem \ref{Symmetry result}]
	Using the fact that $p \ast f = (1-\partial_x^2)^{-1} f$,  we rewritten (\ref{travelling wave solution eq-1}) as the following form: 
	\begin{align}
		-c u_x + \sigma uu_x +p_x \ast \left(\frac{3-\sigma}{2} u^2 + \frac{\sigma}{2} u_x^2 + \kappa u  \right) + p\ast \partial_x^{-1} u_{yy} =0.\label{symmetry eq-1}
	\end{align}
	Then applying the operator $\partial_x^{-1}$ to  (\ref{symmetry eq-1}), one has
	\begin{align}
		cu - \partial_x^{-2} p \ast u_{yy} - \kappa p \ast u=\frac{\sigma}{2} u^2 + p \ast \left( \frac{3-\sigma}{2} u^2 +\frac{\sigma}{2} u_x^2 \right) . \label{symmetry eq-2}
	\end{align}
	Now  we put $ u = \partial_x^2 (1-\partial_x^2) \varphi= \varphi_{xx} - \varphi_{xxxx}$, then (\ref{symmetry eq-2}) is transformed to 
	\begin{align}
		& \left( c-\kappa \right) \varphi_{xx} - c\varphi_{xxxx} - \varphi_{yy} \nonumber \\
		=& \frac{\sigma}{2} \left( \varphi_{xx} - \varphi_{xxxx} \right)^2 + p \ast \left(  \frac{3-\sigma}{2} \left(  \varphi_{xx}-\varphi_{xxxx} \right)^2 +\frac{\sigma}{2} \left(  \varphi_{xxx} -\varphi_{xxxxx} \right)^2  \right)  , \nonumber 
	\end{align}
	which can be written in the integral form  as follows
	\begin{align}
		\varphi  = K \ast f (\varphi), \label{symmetry eq-4}
	\end{align}
	where $f(\varphi) = \frac{\sigma}{2} \left( \varphi_{xx} - \varphi_{xxxx} \right)^2 + p \ast \left(  \frac{3-\sigma}{2} \left(  \varphi_{xx}-\varphi_{xxxx} \right)^2 +\frac{\sigma}{2} \left(  \varphi_{xxx} -\varphi_{xxxxx} \right)^2  \right) $ and the Fourier transform of $K$ is given by 
	\begin{align}
		\widehat{K} = \frac{1}{ -(c-\kappa) \xi^2 -c \xi^4 + \eta^2} . \nonumber 
	\end{align}
	In view of    $0 \leq \sigma \leq 3 $ and $0< c\leq \kappa $, we deduce that  any solution $\varphi$ of (\ref{symmetry eq-4}) is nonnegative. Then by going along the similar lines in  \cite{Esfahani-Levandosky PAMS 2019}, one can readily complete the proof of Theorem \ref{Symmetry result}.
\end{proof} 
Finally, we  investigate the existence of peaked solitary-wave solutions to the rod-KP equation (\ref{Rod-KP equation form-1}) in the form $u(t,x,y) = ce^{-|x+\beta y -ct|}$, $c\in \mathbb{R}$.  Such type of solutions are understood  as weak solutions in the following sense. 
\begin{definition1}
	Let $u_0 \in H^1(\mathbb{R}^2)$. If  the function $u \in C([0,T); H^1_{loc} (\mathbb{R}^2))$ satisfies the following identity:
	\begin{align}
		\int_0^T \int_{\mathbb{R}^2} & \left(-\varphi_{t,x} u + \varphi_x \left( \sigma uu_x  + p_x \ast  \left( \frac{3-\sigma}{2}u^2 +\frac{\sigma}{2}u_x^2 - \kappa u \right) \right) -\varphi_{yy} p \ast u \right) dxdydt \nonumber \\
		& + \int_{\mathbb{R}^2} u_0 (x,y) \varphi_x (0,x,y)  dxdy =0 , \nonumber 
	\end{align}
	for any smooth test function $ \varphi(t,x) \in C_c^{\infty} ([0,T) \times \mathbb{R}^2)$, then $u$ is  called a weak solution to the rod-KP equation (\ref{the nonlocal weak form}). And $u$ is said to be a global weak solution, if it is a weak solution on $[0,T)$ for every $T>0$. 
\end{definition1}

\begin{proof}[Proof of Theorem \ref{Peaked solitary wave}]
	Assume that $u_c(t,x,y) = ce^{-|x+\beta  y -ct|}$. Then we have 
	\begin{align}
		\partial_t u_c= c \cdot{\rm sign} (x+\beta y -ct)\cdot u_c, \quad \partial_x u_c  =-{\rm sign} (x+ \beta y -ct) \cdot u_c, \nonumber 
	\end{align}
	which implies $(\partial_x u_c)^2 =u_c^2$ and $\frac{3-\sigma}{2} u_c^2 +\frac{\sigma}{2} u_c^2 = \frac{3}{2} u_c^2$.  Note that $p_x(x) =-\frac{1}{2} {\rm sign} (x) e^{-|x|}$ for $x\in \mathbb{R}$, we have 
	\begin{align}
		p_x \ast \left( \frac{3-\sigma}{2} u^2 + \frac{\sigma}{2} u_x^2 \right) =& \int_{-\infty}^{+\infty} \left( -\frac{1}{2} {\rm sign } (x-z) e^{-|x-z|} \right) \frac{3}{2} c^2 \cdot e^{-2|z+ \beta y -ct|} dz \nonumber \\
		= & \int_{-\infty}^{+\infty} \left( -\frac{1}{2} {\rm sign } (x+\beta y -z) e^{-|x+\beta y-z|} \right) \frac{3}{2} c^2 \cdot e^{-2|z-ct|} dz . \label{peaked solitary-wave solution proof eq-1}
	\end{align}
	When $x+ \beta y >ct$, we split the right hand side of (\ref{peaked solitary-wave solution proof eq-1}) into the following three parts
	\begin{align}
		I=&  p_x \ast \left( \frac{3-\sigma}{2} u_c^2 +\frac{\sigma}{2} (\partial_x u_c)^2 \right) (t,x,y) \nonumber \\
		= & \left( \int_{-\infty}^{ct} + \int_{ct}^{x+\beta y} + \int_{x+\beta y}^{+\infty} \right)  \left( -\frac{1}{2} {\rm sign } (x+\beta y -z) e^{-|x+\beta y-z|} \right) \frac{3}{2} c^2 \cdot e^{-2|z-ct|} dz \nonumber \\
		= : & I_1 + I_2 +I_3. \nonumber 
	\end{align}
	We directly compute $I_1$ as follows:
	\begin{align}
		I_1= -\frac{3}{4} c^2 \int_{-\infty}^{ct} e^{-(x+\beta y -z)} e^{-2(ct-z)} dz = -\frac{1}{4} c^2 e^{-(x+\beta y -ct)}. \nonumber 
	\end{align}
	In a similar manner, we have 
	\begin{align}
		I_2= -\frac{3}{4} c^2 \int_{ct }^{x+\beta y} e^{-(x+\beta y -2ct)-z}  dz =\frac{3}{4} c^2 \left(  e^{-2(x+\beta y -z) } - e^{-(x+\beta y -ct)} \right) , \nonumber 
	\end{align}
	and 
	\begin{align}
		I_3 = \frac{3}{4}c^2 \int_{x+\beta y}^{+\infty} e^{x+\beta y +2ct -3z} dz = \frac{1}{4} c^2 e^{_2(x+\beta y -ct)} .  \nonumber 
	\end{align} 
	Thus, we deduce that for $x+\beta y >ct$, 
	\begin{align}
		I= c^2 \left(  e^{-2(x+\beta y -z) } - e^{-(x+\beta y -ct)} \right) . 
		\label{peaked solitary-wave solution proof eq-2}
	\end{align}
	While for the case $x+ \beta y < ct$, we split the right side of (\ref{peaked solitary-wave solution proof eq-1}) into the following three parts
	\begin{align}
		II=&  p_x \ast \left( \frac{3-\sigma}{2} u_c^2 +\frac{\sigma}{2} (\partial_x u_c)^2 \right) (t,x,y) \nonumber \\
		= & \left( \int_{-\infty}^{x+\beta y} + \int_{x+\beta y }^{ct} + \int_{ct }^{+\infty} \right)  \left( -\frac{1}{2} {\rm sign } (x+\beta y -z) e^{-|x+\beta y-z|} \right) \frac{3}{2} c^2 \cdot e^{-2|z-ct|} dz \nonumber \\
		= : & II_1 + II_2 +II_3. \nonumber 
	\end{align}
	A direct computation gives rise to 
	\begin{align}
		II_1 =-\frac{3}{4}c^2\int_{-\infty}^{x+\beta y} e^{-(x+\beta y +2ct)} e^{3z} dz = -\frac{1}{4} c^2e^{2(x+\beta y -ct)}, \nonumber 
	\end{align}
	\begin{align}
		II_2 = \frac{3}{4}c^2 \int_{x+\beta y }^{ct} e^{x+\beta y -2ct}e^{z} dz = \frac{3}{4}c^2 \left( e^{x+\beta y -ct} - e^{2(x+\beta y -ct)} \right) ,\nonumber 
	\end{align}
	\begin{align}
		II_3 = \frac{3}{4}c^2 \int_{ct}^{+\infty} e^{x+\beta y +2ct} e^{-3z} dz = \frac{1}{4} c^2 e^{x+\beta y -ct}. \nonumber 
	\end{align}
	We deduce that for $x+\beta y \leq ct$, 
	\begin{align}
		II= c^2 \left( e^{x+\beta y -ct} - e^{2(x+\beta y -ct)} \right) , \nonumber 
	\end{align}
	which along with (\ref{peaked solitary-wave solution proof eq-1}) yields 
	\begin{align}
		p_x \ast \left( \frac{3-\sigma}{2} u_c^2 +\frac{\sigma}{2} (\partial_x u_c)^2 \right)  = \left\{ 
		\begin{aligned}
			& c^2 \left(  e^{-2(x+\beta y -z) } - e^{-(x+\beta y -ct)} \right) , \quad &&x+\beta y >ct, \\ 
			& c^2 \left( e^{x+\beta y -ct} - e^{2(x+\beta y -ct)} \right) , \quad && x+ \beta y \leq ct. 
		\end{aligned}
		\right.	\label{peaked solitary-wave solution proof eq-3}
	\end{align}
	On the other hand, one can deduce 
	\begin{align}
		\partial_t u_c + \sigma u_c \partial_x u_c = c\left(c-c\sigma e^{-|x+\beta y -ct|} \right) e^{-|x+\beta y -ct|},\nonumber 
	\end{align}
	which implies that
	\begin{align}
		\partial_t u_c + \sigma u_c \partial_x u_c = \left\{ 
		\begin{aligned}
			&c^2 \left( e^{-(x+\beta y -ct)} - \sigma e^{-2(x+\beta y -ct)} \right) , \quad &&x+\beta y >ct, \\
			&c^2 \left( \sigma e^{2(x+\beta y -ct)} - e^{x+\beta y-ct}\right) , \quad && x+\beta y \leq ct. 
		\end{aligned}
		\right. 	\label{peaked solitary-wave solution proof eq-4}
	\end{align}
	Thanks to (\ref{peaked solitary-wave solution proof eq-3}) and (\ref{peaked solitary-wave solution proof eq-4}), we get 
	\begin{align}
		\partial_t u_c + \sigma u_c \partial_x u_c + p_x \ast \left( \frac{3-\sigma}{2} u_c^2 + \frac{\sigma}{2} (\partial_x u_c)^2 \right) = \left\{ 
		\begin{aligned}
			&\left( 1- \sigma \right) c^2 e^{-2(x+\beta y -ct)}  , \quad &&x+\beta y >ct, \\
			&\left( \sigma-1\right)  c^2e^{2(x+\beta y -ct)} , \quad && x+\beta y \leq ct. 
		\end{aligned}
		\right. 	\label{peaked solitary-wave solution proof eq-5}
	\end{align}
	Sincce $(u_c)_y = \beta  (u_c)_x$, we have $p \ast\left(  \partial_x^{-1}  \partial_y^2 u_c\right) = \beta^2 p \ast (\partial_x u_c)$. Next, we consider 
	\begin{align}
		p_x \ast u_c = -\frac{1}{2}c  \int_{-\infty}^{+\infty} {\rm sign} \left( x+\beta y-z\right) e^{-|x+\beta y -z|} e^{-|z-ct|} dz . 
		\label{peaked solitary-wave solution proof eq-6}
	\end{align}
	When $x+\beta y >ct$, we split the right side of (\ref{peaked solitary-wave solution proof eq-6}) into the following three parts:
	\begin{align}
		III= & p_x \ast u_c (t,x,y) \nonumber \\
		= & -\frac{1}{2} c \left(  \int_{-\infty}^{ct} + \int_{ct}^{x+\beta y} + \int_{x+\beta y}^{+\infty} \right) {\rm sign} (x+\beta y -z) e^{-|x+\beta y -z|-|z-ct|} dz \nonumber \\
		= : & III_1 + III_2 + III_3.	\label{peaked solitary-wave solution proof eq-7}
	\end{align}
	We directly compute $III_1$, $III_2$ and $III_3$ as follows
	\begin{align}
		III_1 = -\frac{1}{2} c \int_{-\infty}^{ct} e^{-(x+\beta y +ct)} e^{2z} dz = -\frac{1}{4} c e^{-(x+\beta y -ct)},	\label{peaked solitary-wave solution proof eq-8}
	\end{align}
	\begin{align}
		III_2 =-\frac{1}{2} \int_{ct}^{x+\beta y} e^{-(x+\beta y-ct)} dz= -\frac{1}{2} c \left(x+\beta y -ct  \right) e^{-(x+\beta y -ct)} , 	\label{peaked solitary-wave solution proof eq-9}
	\end{align}
	\begin{align}
		III_3= \frac{1}{2}c \int_{x+\beta y}^{+\infty} e^{x+\beta y +ct} e^{-2z} dz =\frac{1}{4} c e^{-(x+\beta y -ct)}. 	\label{peaked solitary-wave solution proof eq-10}                                                                                                                            
	\end{align}
	Plugging (\ref{peaked solitary-wave solution proof eq-8})-(\ref{peaked solitary-wave solution proof eq-10}), we deduce that for $x+\beta y >ct$, 
	\begin{align}
		p_x \ast u_c = -\frac{1}{2} c (x+\beta y -ct) e^{-(x+\beta y -ct)}. \label{peaked solitary-wave solution proof eq-10.5}
	\end{align}
	While for the case $x+\beta y \leq ct$, we split the right hand side of (\ref{peaked solitary-wave solution proof eq-6}) into the following three parts
	\begin{align}
		IV= & p_x \ast u_c (t,x,y) \nonumber \\
		= & -\frac{1}{2} c \left(  \int_{-\infty}^{x+\beta y} + \int_{x+\beta y}^{ct} + \int_{ct}^{+\infty} \right) {\rm sign} (x+\beta y -z) e^{-|x+\beta y -z|-|z-ct|} dz \nonumber \\
		= : & IV_1 + IV_2 + IV_3.	\label{peaked solitary-wave solution proof eq-11}
	\end{align}
	For $IV$, a direct computation gives rise to
	\begin{align}
		IV_1 = -\frac{1}{2} c\int_{-\infty}^{x+\beta y} e^{-(x+\beta y +ct)} e^{2z} dz = -\frac{1}{4} e^{x+\beta y -ct}, 
		\label{peaked solitary-wave solution proof eq-12}
	\end{align}
	\begin{align}
		IV_2 = \frac{1}{2} c \int_{x+\beta y}^{ct} e^{x+\beta y -ct} dz= -\frac{1}{2} (x+\beta y -ct) e^{x+\beta y -ct}, \label{peaked solitary-wave solution proof eq-13}
	\end{align}
	\begin{align}
		IV_3= \frac{1}{2} c \int_{ct}^{+\infty} e^{x+\beta y +ct} e^{-2z} dz = \frac{1}{4}ce^{x+\beta y -ct} .  \label{peaked solitary-wave solution proof eq-14}
	\end{align}
	Plugging (\ref{peaked solitary-wave solution proof eq-12})-(\ref{peaked solitary-wave solution proof eq-14}) into (\ref{peaked solitary-wave solution proof eq-11}), we deduce that for $x+\beta y \leq ct$, 
	\begin{align}
		p_x \ast u_c =-\frac{1}{2} c (x+\beta y -ct) e^{x+\beta y -ct} ,\nonumber 
	\end{align}
	which along with (\ref{peaked solitary-wave solution proof eq-10.5}) yields
	\begin{align}
		p_x \ast u_c = \left\{   
		\begin{aligned}
			& -\frac{1}{2} c (x+\beta y -ct) e^{-(x+\beta y -ct)}, \quad && x+\beta y >ct,\\
			& -\frac{1}{2} c (x+\beta y -ct) e^{x+\beta y -ct}, \quad && x+\beta y \leq ct. 
		\end{aligned} \right. 
		\label{peaked solitary-wave solution proof eq-15}
	\end{align}
	
	Therefore, thanks to (\ref{peaked solitary-wave solution proof eq-5}) and (\ref{peaked solitary-wave solution proof eq-15}), we see that the  rod-KP equation (\ref{Rod-KP equation form-1}) has a global weak solutin in the peak form of  $u(t,x,y) = ce^{-|x+\beta y -ct|}$ for some constant $\beta \in \mathbb{R}$ if  and only if there holds $\sigma =1$ and $\beta^2 + \kappa  =0$. 
\end{proof}

\noindent {\bf Acknowledgments.}  This work was  partially supported by the National Natural Science Foundation of China under grant 11971188.

\end{document}